%%%%%%%%%%%%%%%%THIS IS AN AMS-TEX DOCUMENT%%%%%%%%%%%%%%%%%%%%%%%%%%%%%%

\magnification=\magstep1
\vsize=22truecm
\input amstex
\documentstyle{amsppt}
\leftheadtext{K. B\"or\"oczky, A. Heppes, E. Makai, Jr.}
\rightheadtext{Densest packings of translates of strings and layers of balls}
\topmatter
\title Densest packings of translates of strings and layers of balls
\endtitle
\author K. B\"or\"oczky$^*$, A. Heppes$^{**}$, E. Makai, Jr.$^{***}$ 
\vskip.5cm
{\rm{L. E\"otv\"os University, Faculty of Science, Institute of
Mathematics, Chair of Geometry,
1117 Budapest, P\'azm\'any P. S\'et\'any 1/a, Hungary (K. B.)}}
\vskip.1cm
{\centerline{\rm{1124 Budapest, V\'ercse u. 24/A, Hungary (A. H.)}}}
\vskip.1cm
\centerline{
{\rm{MTA Alfr\'ed R\'enyi Institute of Mathematics,}}} 
\centerline{
{\rm{H-1364 Budapest, Pf. 127, Hungary (E. M., Jr.)}}}
\vskip.1cm
{\rm{http://www.renyi.mta.hu/\~{}makai}}
\vskip.1cm
{\rm{E-mail: boroc\@ludens.elte.hu, h9202hep\@ella.hu,
makai.endre\@renyi.mta.hu}}
\endauthor
\thanks 
$^*$ 
Research (partially) supported by Hungarian National Foundation for 
Scientific Research, grant nos. T04520, K81146, ...???
\newline
$^{**}$ 
Research (partially) supported by Hungarian National Foundation for 
Scientific Research, grant nos. T04520, K81146, ...???
\newline
$^{***}$ Research (partially) supported by Hungarian National Foundation for 
Scientific Research, grant nos. T04520, K68398, K75016, K81146\endthanks
\subjclass {\it 2000 Mathematics Subject Classification.} Primary: 52C17;
Secondary: 52C07, 52A45, 52A43\endsubjclass
\abstract Let $L \subset {\Bbb R}^3$ be the union of unit balls, whose centres 
lie on the $z$-axis, and are equidistant with distance 
$2d \in [2, 2\sqrt{2}]$. Then a packing of
unit balls in ${\Bbb R}^3$ consisting of translates of $L$ has a density at
most $\pi /(3d\sqrt{3-d^2})$, with equality for a certain lattice packing of
unit balls.
Let $L \subset {\Bbb R}^4$ be the union of unit balls, whose centres lie on
the $x_3x_4$ coordinate plane, 
and form either a square lattice or a regular triangular
lattice, of edge length $2$. Then a packing of
unit balls in ${\Bbb R}^4$ consisting of translates of $L$ has a density at
most $\pi ^2/16$, with equality for the densest lattice packing of unit 
balls in ${\Bbb R}^4$. This is the first class of non-lattice packings of unit
balls in ${\Bbb R}^4$, for which this conjectured upper bound for the packing
density of balls is proved.
Our main tool for the proof is a theorem on
$(r,R)$-systems in ${\Bbb R}^2$. If $R/r \le 2 \sqrt{2}$, then the Delone
triangulation associated to this $(r,R)$-system has the following property.
The average area of a Delone triangle is at least $\min \{ V_0, 2r^2 \} $,
where $V_0$ is the infimum of the areas of the non-obtuse Delone triangles.
This general theorem has applications also in other problems about packings:
namely for $2r^2 \ge V_0$ it is sufficient to deal only
with the non-obtuse Delone triangles, which is in general a much easier task.
Still we give a proof of an unpublished theorem of L. Fejes T\'oth and E (=J.)
Sz\'ekely: for the $2$-dimensional analogue of our question about equidistant
strings of unit
balls, we determine the densest packing of translates of an equidistant
string of unit circles with distance $2d$, 
for the first non-trivial interval $2d \in (2{\sqrt{3}},4)$.
\endabstract
\endtopmatter\document

\head{\S 1 Introduction}\endhead

The well known Kepler conjecture states that in ${\Bbb R}^3$ a packing of unit
balls has density at most the density of the densest lattice packing of balls
in ${\Bbb R}^3$, i.e., $\pi / \sqrt{18}$ (cf. \cite{FT72}, \cite{Ro},
\cite{FT64}, \cite{GL}, \cite{BKJ}). 

As a special case, the first named author posed in \cite{Bo}, in 1975, 
the following problem. Prove that
the density of a packing of unit balls in ${\Bbb R}^3$, consisting of parallel
strings of unit balls, whose centers are equally spaced on a straight line,
at distances $2$, is at most $\pi / \sqrt{18}$.

\cite{BKM91}, in 1991, solved this question in the affirmative. For some time
their result was the most general result about a class of 
packings of unit balls in ${\Bbb R}^3$, for whose density the sharp upper bound
$\pi / \sqrt{18}$ was proved.  

Now of course this result is superceded by T. C. Hales, S. Ferguson 
\cite{HF}, who proved Kepler's
conjecture in full generality. A revision of their original proof is given in
\cite{HHMNOZ}. A complete formal proof of the Kepler conjecture that can be
verified by automated proof checking software (cited from \cite{V}) is given
in \cite{HABD...}.

L. Fejes T\'oth \cite{FT} generalized the question of \cite{Bo}, 
requiring that the
distance of the neighbouring centres of the balls in a string should be $2d$
rather than $2$, where $d \in [1, {\sqrt{2}}]$, 
and asked for the maximal density under this more general
condition. We will prove a general theorem, which contains as a special case
the solution of L. Fejes T\'oth's question. We will
obtain that among the extremal packings there are lattice packings, 
which we will concretely determine. We note that the authors of \cite{BKM91}
were aware that their paper gave possibility to prove L. Fejes T\'oth's
conjecture for $d$ sufficientky close to $1$, but this was not included in
their paper.

Observe that the result of \cite{HF} does
not solve this problem, except for those values of $d$ for which 
$2d$ is the distance of some centres of unit balls in a packing of unit balls
consisting of closely packed regular hexagonal layers
which is periodic with some period $k$ (i.e., some fixed
translation carries the $n$'th hexagonal layer to the $(n+k)$'th layer, for
each integer $n$) and $2d$ is the distance of the centres of some balls from
the first and $(k+1)$'th hexagonal layers.
Such values are, e.g., $d=1$, $d=\sqrt{2}$ (for the densest lattice packing)
and $d=\sqrt{8/3}$ (for the densest regular non-lattice packing).

The analogous question for the plane is that of the maximum 
density of a packing of unit circles in ${\Bbb R}^2$,
consisting of parallel
strings of unit circles, whose centers are equally spaced on a straight line,
at distances $2d$ (where $d \ge 1$). This question
is solved for $1 \le d \le \sqrt{3}$ in
\cite{FT62} (in the form of packing geodesic unit circles on the surface of
a circular cylinder with base of perimeter $2d$): a densest packing is 
lattice-like, with the corresponding point lattice spanned by the vertices of
a triangle of sides $2$, $2$, $2d$.
Further, this question is
solved for $\sqrt{3} \le d \le 2$ \cite{FT62} (he did not give the densest
packing) and
\cite{Sz}, oral communication. In the second part 

\newpage

of the introduction we will
give a proof for this case. This
question is unsolved for all values $d > 2$, except those for which $2d$ is
the distance of two centres of circles in the densest lattice packing of unit
circles in ${\Bbb R}^2$. 

\cite{Sz} conjectured that for any $d$ one of the
densest packings for the planar case 
is obtained in the following way. We choose a natural number
$k$, and consider $k$ neighbouring strings (of touching unit circles)
of a densest lattice packing of unit
circles. We translate this block of $k$ strings periodically so that any two
neighbouring blocks touch each other and for a circle of the first string of a
block and some circle of the first string 
of the following block the distance of the
centres is $2d$. This number $k$ equals $1$ for $1 < d < \sqrt{3}$
(cf. [FT62]), and
$2$ for $\sqrt{3} < d < 2$ (cf. [Sz], oral communication). 
By this definition we get in fact a packing of unit circles, and also a
packing of equidistant strings of unit circles with distance $2d$.
The next interesting interval is $2 < d  < {\sqrt{7}}$. Then still we should
have $k=2$ (for $k \ge 4$ the symmetry axes of the $k$'th neighbour strings
have a distance at least $2{\sqrt{12}}$, and for $k=3$ the centre of any circle
touching the $3$'rd neighbour string from the ``other side'' has a distance at
least $2{\sqrt{7}}$ from the centre of any circle in the original string. The
proof of the theorem of L. Fejes T\'oth and E. (=J.) Sz\'ekely in the second
half of the introduction does not seem to be suitable for the case 
$2 < d  < {\sqrt{7}}$.

A related problem is investigated in J. Moln\'ar \cite{M78}. Namely the
maximal packing density of unit balls in a parallel slab in ${\Bbb{R}}^3$,
of width $w \in [2, 2 + {\sqrt{2}}]$ was determined, as a function of $w$.
The analogous problem of the
maximal packing density of unit circles in a parallel strip in ${\Bbb{R}}^2$,
of width $w \in [2, 2 + 2{\sqrt{2}}]$ was determined, as a function of $w$, by
G. Kert\'esz \cite{Ke}, unpublished, and for width $w \in [2, 2 + 2{\sqrt{3}}]$
by Z. F\"uredi \cite{Fu}.

The idea of our proofs is taken from \cite{BKM}. Our problem is reduced
essentially to a planar problem, cf. Theorem 2.9. There we have a nice system
of points on the plane (an $(r,R)$-system), 
and consider the Delone triangulation associated to
it. The acute or right
Delone triangles present no problem, the problem is only with obtuse Delone
triangles. However, an obtuse Delone triangle has a neighbour at its longest
side, for which the common side is longer than guaranteed by the problem, and 
its vertex opposite to the common side is outside of the circumcircle of the
obtuse triangle. These two properties will have the consequence that ``the
smaller area of the obtuse triangle is compensated by the surplus of the
area of its neighbouring Delone triangle at its longest side''. Analogously: if
some Delone triangle shares the longest sides of two or three obtuse
triangles, then
``its area surplus compensates the small areas of the mentioned two
or three obtuse triangles''. Details cf. in \S 4.

Still we note that the analogue of the Kepler conjecture was most recently
solved in
${\Bbb{R}}^8$ by M. Viazovska \cite{V} 
(the solution is the densest lattice packing of unit balls in ${\Bbb{R}}^8$) 
and by H. Cohn, A. Kumar, S.D. Miller, D. Radchenko, M. Viazovska
\cite{CKMRV}

\newpage

in ${\Bbb{R}}^{24}$ (the solution is the Leech lattice, which was proved to be
the densest lattice packing by H. Cohn, A. Kumar \cite{CK} in 2009).

\vskip.1cm

{\centerline{***}}

Since for the problem of the densest packing of
translates of equidistant strings of unit circles
with equal 
distances $2d$, for $d \in ({\sqrt{3}},2)$, no proof seems to have been
published, not even the densest packing or the packing density seems to have
been published, we give the not complicated proof of the respective result.

\vskip.1cm

{\bf{Theorem.}} (L. Fejes T\'oth and E. (=J.) Sz\'ekely).
{\it{For $d \in ({\sqrt{3}},2)$,
the maximum density of packings of translates of equidistant strings of unit
circles, with equal distances $2d$,
is $2 \pi /[d({\sqrt{4-d^2}}+d{\sqrt{3}})]$. This density is
attained for the packing of translates of our equidistant string of unit
circles described in the next paragraph.}}
 
\vskip.1cm

We consider the axes of symmetry of the strings (later called axes of the
strings) as horizontal, and then the
strings constituting our packing have a natural order: namely according to the
$y$-coordinates of the centres of the circles in the strings. If the $0$'th
string has two neighbourly centres of circles $A_0,A_1$, with $|A_0A_1|=2d$,
then let the $1$'st string lie above the $0$'th string and have a circle with
centre $B_0$ where $|A_0B_0|=|A_1B_0|=2$. Similarly, let the $2$'nd string lie
above the $1$'st string, and have a circle with centre $C_1$, where
$|A_1C_1|=|B_0C_1|=2$. Then let 
the translation carrying the $0$'th string to the
$2$'nd string carry the $i$'th string to the $(i+2)$'nd string for each
integer $i$. 
(This system is identical with the one conjectured by E. (=J.)
Sz\'ekely, with $k=2$. In particular, the above construction gives
a packing of unit circles, and also a
packing of equidistant strings of unit circles with distance $2d$.)

\vskip.1cm

{\it{Proof of the theorem of L. Fejes T\'oth and E. (=J.) Sz\'ekely.}}
It will be sufficient to prove that the distance of the axes of the $0$'th and
$2$'nd strings is at least the distance in the above construction, i.e.,
$({\sqrt{4-d^2}}+{\sqrt{3}}d)/2$. {\it{Indirectly 
suppose that this distance is less than 
$({\sqrt{4-d^2}}+{\sqrt{3}}d)/2$.}} Then {\it{this distance is also at most
$({\sqrt{4-d^2}}+{\sqrt{3}}d)/2$,}} and we are going to show that 
$$
\cases
{\text{under this new (weaker) assumption the distance of the axes}} \\
{\text{of the }}
0{\text{'th and }} 2 {\text{'nd strings is equal to }} 
({\sqrt{4-d^2}}+{\sqrt{3}}d)/2 \,.
\endcases
\tag 0.1
$$
This will prove our indirect statement.

We take $A_0$ as the origin, and then $A_1=(2d,0)$, $B_0=(d,{\sqrt{4-d^2}})$,
and $C_1=\left( (3d+{\sqrt{3}}{\sqrt{4-d^2}})/2, 
({\sqrt{3}}d+{\sqrt{4-d^2}})/2 \right) $. 
Let $u$ and $v$ be the vectors of translations carrying the $0$'th string to
the $1$'st string and the $1$'st string to the $2$'nd string, where $u$ and
$v$ lie in the strip $0 \le x \le 2d$. Then by the packing property and our
hypothesis 

\newpage

$$
\cases
u,v \in K:= \{ (x,y) \in {\Bbb{R}}^2 \mid 0 \le x \le 2d, \,\,0 \le y \le \\
({\sqrt{4-d^2}}+{\sqrt{3}}d)/2, \,\,
\| (x,y) \| \ge 2, \,\,\|(x,y)-(2d,0) \| \ge 2 \} \,.
\endcases
\tag 0.2
$$
Then translation by 
$w:=u+v$ carries the $0$'th string to the $2$'nd string, hence
analogously to \thetag{0.2}, by the packing property and our hypothesis
$$
w \in K+(2id,0) {\text{ for some integer }} i\,.
\tag 0.3
$$
Since $w=u+v$ lies in the strip $0 \le x \le 4d$,
so actually 
$$
w \in K {\text{ or }} w \in K+(2d,0) \,. 
\tag 0.4
$$
Also $w=u+v \in K+K$. Hence by \thetag{0.4} we
have
$$
w = u + v \in K \cap (K+K) \ne \emptyset {\text{ or }} 
w = u + v \in \left( K +(2d,0) \right)
\cap (K+K) \ne \emptyset \,.
\tag 0.5
$$

The set $K$ is bounded by the segment $[C_0,C_1]$, where $C_0$ is the
symmetric image of $C_1$ with respect to the line $x=d$, and by two concave
circular arcs ${\widehat{B_0C_0}}$ and ${\widehat{B_0C_1}}$ of radii $2$,
central angles $\pi /3$ and centres $A_0$ and $A_1$. Thus $K$ is a ``concave
arc triangle'', consisting of the Jordan curve $[C_0,C_1] \cup 
{\widehat{B_0C_0}} \cup {\widehat{B_0C_1}}$ and its interior. Its convex hull
${\text{conv}}K$ is the triangle $\Delta B_0C_0C_1$. 

Since $B_0,C_0,C_1 \in K$, we have $2B_0,2C_0,2C_1,B_0+C_0,B_0+C_1,C_0+C_1 \in
K+K$. Since $K \subset \Delta B_0C_0C_1$, we have $K+K \subset 2 \cdot \Delta
B_0C_0C_1$. We are going to bound $K+K$ from outside better.

It will be more convenient to deal with $(K+K)/2$. From above we have 
$$
\cases
\{ B_0,C_0,C_1,(B_0+C_0)/2, (B_0+C_1)/2, (C_0+C_1)/2 \} \\
\subset (K+K)/2 \subset \Delta B_0C_0C_1 \,.
\endcases
\tag 0.6
$$

Diminishing $K$ from any of its vertices in ratio $1/2$, we see that $(K+K)/2$
contains $[C_0,(C_0+C_1)/2]$, $[(C_0+C_1)/2,C_1]$ 
and and four (concave) circular arcs of radii $1$ and
central angles $\pi /3$, say, ${\widehat{B_0 \left( (B_0+C_0) /2 \right) }}$, 
${\widehat{\left( (B_0+C_0) /2 \right)  C_0}}$,
\newline
${\widehat{B_0 \left( (B_0+C_1) /2 \right) }}$,
${\widehat{\left( (B_0+C_1) /2 \right) C_1}}$, which six arcs altogether 
form a Jordan
curve $J$, such that $({\text{int}}\,J) \cup J \supset K$, where
${\text{int}}\,J$ denotes the interior of $J$. Of course all
three diminished
copies of ${\text{int}}\,K$, from any of its vertices in ratio $1/2$,
are contained in ${\text{int}} \left( (K+K)/2 \right) $.

\newpage

Let $x,y \in K$, and let one of $x$ and $y$ be an interior point of $K$. Then
$(x+y)/2 \in {\text{int}} \left( (K+K)/2 \right)$. The case when one of $x$
and $y$ is a vertex of $K$ was investigated in the last paragraph.

There remains the case when
$x,y \in K$ and both are relative inner points of some
arc sides of $K$. The slope of ${\text{bd}}K$ at relative inner points of
$[C_0,C_1]$, ${\widehat{B_0C_0}}$ and ${\widehat{B_0C_1}}$ is $0$, negative
and positive, respectively. 
Therefore if $x$ and $y$ are relative inner points of different
arc sides of $K$, then these slopes are different, hence 
$(x+y)/2 \in
{\text{int}} \left( (K+K)/2 \right) $. If $x,y \in [C_0,C_1]$, then
$(x+y)/2 \in [C_0,C_1]$. If $x,y$ are different and both belong to the relative
interior of either ${\widehat{B_0C_0}}$ or ${\widehat{B_0C_1}}$, then the
slopes of ${\text{bd}}\,K$ at $x.y$ are different, hence, like above,
$(x+y)/2 \in {\text{int}} \left( (K+K)/2 \right) $. If $x=y$ belongs to 
the relative
interior of either ${\widehat{B_0C_0}}$ or ${\widehat{B_0C_1}}$, then
$x=y=(x'+y')/2$ for some $x',y' \in {\text{int}}\,K$, hence $x=y \in
{\text{int}} \left( (K+K)/2 \right) $.

Summing up: boundary points of $(K+K)/2$ lie either on $J$, or on the
relative interiors of the diminished opposite arc sides of $K$, from its three
vertices in ratio $1/2$ (which bound a
convex arc-triangle $T$ with vertices the side midpoints of ${\text{conv}}\,K$).

Next we show that $J \subset {\text{bd}} \left( (K+K)/2 \right) $.
We show this for the subarcs
\newline
${\widehat{B_0 \left( (B_0+C_0)/2 \right) }}$ and 
${\widehat{B_0 \left( (B_0+C_1)/2 \right) }}$ of $J$. (For the other subarcs the
argument is the same.) These bound (partially) 
a diminished copy of $K$ from $B_0$, in ratio $1/2$,
whose interior is a subset of ${\text{int}} \left( (K+K)/2 \right) $. If one
of our considered two arcs were not contained in ${\text{bd}} \left(
(K+K)/2 \right)$, then some inner point of $(K+K)/2$ (namely in the diminished
copy of ${\text{int}}\,K$, from $B_0$, in ratio $1/2$)
could be connected by
an arc to infinity, avoiding ${\text{bd}}
\left( (K+K)/2 \right) $,
contradicting $(K+K)/2 \subset {\text{conv}}\,K$. 

Then by $J \subset {\text{bd}} \left( (K+K)/2 \right) $
the relative
interiors of the diminished copies of the opposite arc sides 
(i.e., the relative interiors
of the sides of the above arc triangle $T$) cannot have
points on the {\it{outer boundary of $(K+K)/2$, i.e., on the boundary of the
unbounded connected component of ${\Bbb{R}}^2 \setminus \left( (K+K)/2 \right)
$}}. Therefore the outer boundary of $(K+K)/2$ equals $J$. 

If $(K+K)/2$ had
some point outside $J$, then it would have also some outer boundary point
outside $J$, a contradiction. Therefore
$(K+K)/2 \subset ({\text{int}}\,J) \cup
J$. Therefore $K+K \subset L:=2\left( ({\text{int}}\,J) \cup J \right) $. 
(Actually here equality holds --- but we do not need this. 
Namely, $\emptyset \ne ({\text{int}}\,K) \cap ({\text{int}}\,T)
\subset \left( {\text{int}}\,\left( (K+K)/2 \right) \right)
\cap ({\text{int}}\,T)$. Further, since
$({\text{int}}\,T) \cap \left( {\text{bd}}\,\left( (K+K)/2 \right) \right)
= \emptyset $,
therefore ${\text{int}}\,T \subset {\text{int}}\,\left(
(K+K)/2 \right) \subset (K+K)/2$. Also the 
diminished copies of $K$ from its three vertices, in ratio $1/2$, lie in
$(K+K)/2$ as well. These together prove our claim.)
Then \thetag{0.5} implies
$$
w \in K \cap L \ne \emptyset {\text{ or }} w \in
\left( K+(2d,0) \right) \cap L \ne \emptyset \,.
\tag 0.7
$$
We are going to show that 
$$
\cases
{\text{any point (e.g., }} w) {\text{ of }} K \cap L {\text{ or of }} 
\left( K+(2d,0) \right) \cap L \\
{\text{has }} y {\text{-coordinate }} ({\sqrt{4-d^2}} + {\sqrt{3}}d)/2 \,, 
\endcases
\tag 0.8
$$

\newpage

as promised in \thetag{0.1} (observe that the distance of the axes of the
$0$'th and $2$'nd strings is the $y$-coordinate of $w$).

Observe that we have a symmetric trapezoid $A_0A_1C_1C_0$, with $A_0$ the
origin, which contains $K$. Both $A_0A_1C_1C_0$ and $K$ have the line $x=d$ as
symmetry axis. The vertices of $K$ are $B_0,C_0,C_1$, and the 
boundary of $K$ is 
${\widehat{B_0C_0}} \cup {\widehat{B_0C_1}} \cup [C_0,C_1]$. Then $L$ is
bounded by $[2C_0,C_0+C_1]$, $[C_0+C_1,2C_1]$ 
and four concave circular arcs of radius $2$ and
central angles $\pi /3$, say, ${\widehat{(2B_0)(B_0+C_0)}}$, ${\widehat{
(B_0+C_0)(2C_0)}}$ --- which are translates of the boundary arc
${\widehat{B_0C_0}}$ of $K$ --- and
${\widehat{(2B_0)(B_0+C_1)}}$, ${\widehat{(B_0+C_1)(2C_1)}}$ 
--- which are translates of the boundary arc
${\widehat{B_0C_1}}$ of $K$. Let the images of $A_1,C_0,C_1$ by the translation
through $(2d,0)$ be $A_2,C_2,C_3$. Then the line $x=2d$ is an axis of symmetry
of $L$, and the axially symmetric images of $A_0A_1C_1C_0$ and $K$ with
respect to this line are $A_1A_2C_3C_2$ and $K+(2d,0)$. Therefore also $K \cap
L$ and $\left( K+(2d,0) \right) \cap L$ are axially symmetric images of each
other with respect to
the line $x=2d$. Hence, rather than {\thetag{0.8}}, it suffices to show that 
$$
{\text{any point (e.g., }} w) {\text{ of }} 
K \cap L {\text{ has }} y{\text{-coordinate }} 
({\sqrt{4-d^2}} + {\sqrt{3}}d)/2 \,.
\tag 0.9
$$

The line $(B_0+C_0)(B_0+C_1)$ cuts $L$ into two closed parts, say $L_1$ is the
lower part, and $L_2$ is the upper part. Since the $y$-coordinate of $B_0$ is
positive, the $y$-coordinate of $B_0+C_0$ is larger than the $y$-coordinate
of $C_0$, which equals 
$({\sqrt{4-d^2}}+{\sqrt{3}}d)/2$. Also $K$ lies (not strictly)
below the line $y = ({\sqrt{4-d^2}}+{\sqrt{3}}d)/2$. Thus $K \cap L_2
= \emptyset $, and
$$
K \cap L = K \cap L_1
\,.
\tag 0.10 
$$
The set $L_1$ is a translate of $K$, through the vector $B_0$. The leftmost
point of $L_1$ is $B_0+C_0$, whose $x$-coordinate is greater than the
$x$-coordinate of $B_0$ which equals 
${\sqrt{4-d^2}}$. The part $K_1$ of $K$ (not
strictly) to the left hand side
of the line $x={\sqrt{4-d^2}}$ is therefore disjoint to
$L_1$. Therefore only the part $K_2$ of $K$ (not strictly) to the right 
hand side of this
line can intersect $L_1$, i.e., 
$$
K \cap L_1 = K_2 \cap L_1 \,.
\tag 0.11
$$
Hence 
$$
\cases
K_2 \subset K \subset {\text{conv}}\,K \subset {\text{ the circle with centre
}} B_0 {\text{ and radius }} 2, 
{\text{ and}} \\
{\text{the only point of }} K_2  {\text{ lying on the boundary of this 
circle is }} C_1 \,.
\endcases
\tag 0.12
$$

On the other hand, 
$$
\cases
L_1 {\text{ lies (not strictly) to the right hand side}} \\
{\text{of the arc }} {\widehat{(2B_0)(B_0+C_0)}} {\text{ of the circle
in \thetag{0.12}.}}
\endcases
\tag 0.13
$$

\newpage

Then \thetag{0.10}, \thetag{0.11}, \thetag{0.12} and \thetag{0.13} imply
$$
K \cap L = K \cap L_1 = K_2 \cap L_1 = \{ C_1 \} \,,
\tag 0.14
$$
and the $y$-coordinate of $C_1$ is $({\sqrt{4-d^2}}+{\sqrt{3}}d)/2$, as was
promised in \thetag{0.9}. Hence also \thetag{0.8} and \thetag{0.1} hold, thus
the theorem is proved.
$\blacksquare $

%%%%%%%%%%%%%%%%%%%%%%%%%%%%%%%%%%%%%%%%%%%%%%%%%%%%%%%%%%%%%%%%%%%%%%%%%%%
%%%%%%%%%%%%%%%%%%%%%%%%%%%%%%%%%%%%%%%%%%%%%%%%%%%%%%%%%%%%%%%%%%%%%%%%%%%

\head{\S 2 Results}\endhead

We begin with some notations.
We write $B^n \subset {{\Bbb R}}^n$ for the closed unit ball centred at $0$. 
$B\left( (x_1, \ldots , x_n) ,R \right) \subset {\Bbb{R}}^n$ 
is the closed ball of centre $(x_1, \ldots , x_n)$ and radius $R$. An analogous
notation will be applied for closed balls in lower dimensional subspaces (the
coordinates of the centre of the ball will indicate the subspace).  
For points $A,B \in {{\Bbb R}}^n$ we write $[A,B]$ for the segment with
endpoints $A$ and $B$, $|AB|$ for the length of $[A,B]$,
and for a vector $v \in {\Bbb{R}}^n$ we write $|v|$ for its norm.

By a {\it{point lattice in ${\Bbb{R}}^n$}} we mean an
inhomogeneous lattice, i.e., a translate of a homogeneous lattice (i.e., of a
discrete subgroup of ${\Bbb{R}}^n$, possibly not full dimensional). 
We denote point lattices by $\Lambda $.
By a {\it{lattice vector}} we mean a vector (point) in the
corresponding homogeneous lattice.
{\it{For $L \subset {\Bbb{R}}^n$, a two-dimensional lattice packing of
translates of $L$}} is a packing of
the form $ \{ L+\lambda \mid \lambda \in \Lambda \} $, where 
$\Lambda \subset {\Bbb R}^n$ is a two-dimensional point lattice.

Our Theorem 2.1 about ${\Bbb{R}}^3$ 
generalizes the theorem of \cite{BKM91}, which is the special
case $d=1$ of our Theorem 2.1.

%%%%%%%%%%%%%%%%%%%%%%%%%%%%%%%%%%%%%%%%%%%%%%%%%%%%%%%%%%%%%%%%%%%%%%%%%%%
%%%%%%%%%%%%%%%%%%%%%%%%%%%%%%%%%%%%%%%%%%%%%%%%%%%%%%%%%%%%%%%%%%%%%%%%%%%

\proclaim{Theorem 2.1}
Let $L \subset {\Bbb R}^3$ be the union of a string of unit balls in ${\Bbb
R}^3$, whose centres lie on the $z$-axis, and are equidistant with distance
$2d$, where $1 \le d \le {\sqrt{2}}$. Then any packing of translates of $L$
in ${\Bbb R}^3$
(i.e., any packing of unit balls in ${\Bbb R}^3$,
consisting of entire translates of $L$)
has a density at most $\pi /(3d\sqrt{3-d^2})$, with equality e.g. for the
lattice generated by the
vertices of a tetrahedron with five edges of length $2$ and one edge of length
$2d$. 
Another way of giving this lattice is the following: it is 
generated by the vertices of a rectangular right pyramid with 
base edges of lengths $2$ and $2d$, and lateral edges of length $2$. 
Among
lattice packings of translates of $L$ the above lattice is the unique lattice 
with this maximal density. 
\endproclaim

%%%%%%%%%%%%%%%%%%%%%%%%%%%%%%%%%%%%%%%%%%%%%%%%%%%%%%%%%%%%%%%%%%%%%%%%%%%
%%%%%%%%%%%%%%%%%%%%%%%%%%%%%%%%%%%%%%%%%%%%%%%%%%%%%%%%%%%%%%%%%%%%%%%%%%%

\definition{Remark 2.2}
Consider this packing as one consisting of layers generated by the balls
corresponding to a face of the tetrahedron in the Theorem with sides $2,2,2d$.
Then we obtain packings of translates of $L$ of the same density when these
layers are packed parallelly and closely, two neighbouring layers joining in
any of the
two congruent ways (as in case of the densest packing of unit balls in
${\Bbb{R}}^3$ we
may put on a regular hexagonal layer a next one in two different ways).
As particular cases, we obtain the lattice packing
from Theorem 2.1 and the non-lattice-like regular packing (in analogy with
the densest lattice packing of unit
balls and the densest non-lattice-like regular packing of unit balls in
${\Bbb{R}}^3$).
\enddefinition

%%%%%%%%%%%%%%%%%%%%%%%%%%%%%%%%%%%%%%%%%%%%%%%%%%%%%%%%%%%%%%%%%%%%%%%%%%%
%%%%%%%%%%%%%%%%%%%%%%%%%%%%%%%%%%%%%%%%%%%%%%%%%%%%%%%%%%%%%%%%%%%%%%%%%%%

\newpage

The next theorem is an analogue of Theorem 2.1 for ${\Bbb R}^4$.
It is a well known conjecture 
that in ${\Bbb{R}}^4$ a densest packing of unit balls is
lattice like, with packing lattice the space-centred cubic lattice, and this
ball packing has a density $\pi ^2 / 16$, cf. \cite{Ro}, \cite{GL}, \cite{BKJ}.
Our next Theorem 2.2 about ${\Bbb{R}}^4$ gives the first class of non-lattice 
packings of unit balls in ${\Bbb{R}}^4$, for which the conjectured density
estimate is proved. 
We hope that this result will be soon superceded by the solution of the
densest ball packing in ${\Bbb R}^4$, like it happened with
\cite{BKM91} and \cite{HF} in ${\Bbb{R}}^3$.

\proclaim{Theorem 2.3}
Let $L \subset {\Bbb R}^4$ be the union of the two-dimensional lattice packing
of translates of the unit ball in ${\Bbb R}^4$, where the corresponding
two-dimensional point lattice is either the square lattice with edge length 
$2$, or the regular triangular lattice with edge length $2$. Then any packing of
translates of $L$ in ${\Bbb R}^4$
(i.e., any packing of unit balls in ${\Bbb R}^4$,
consisting of entire translates of $L$)
has a density at most $\pi ^2/16$, i.e., the density of a
densest lattice packing of unit balls in $ {\Bbb R}^4$.
\endproclaim

%%%%%%%%%%%%%%%%%%%%%%%%%%%%%%%%%%%%%%%%%%%%%%%%%%%%%%%%%%%%%%%%%%%%%%%%%%%
%%%%%%%%%%%%%%%%%%%%%%%%%%%%%%%%%%%%%%%%%%%%%%%%%%%%%%%%%%%%%%%%%%%%%%%%%%%

Returning to ${\Bbb R}^3$, we need some notations. 
Let $L$ be the union of a string of unit balls in ${\Bbb R}^3$, 
whose centres lie on the $z$-axis. Let the distance of the
neighbouring centres of
balls in $L$ be always at least $2$, further let them
satisfy that the average distance of the neighbouring balls is $2d$, and this
holds uniformly. By this we mean the following.
$$
\cases
{\text{for any segment }}[(0,0,z-R),(0,0,z+R)]
{\text{ in the }}z{\text{-axis}} \\
{\text{ the quotient of }} 2R {\text{ and the number of ball centres from }}
L \\
{\text{ in }} [(0,0,z-R),(0,0,z+R)] {\text{ tends for }} R \to \infty 
{\text{ to }} 2d, \\
{\text{ and this convergence is uniform for all }} z \in
{\Bbb R}.
\endcases
\tag 2.1
$$

%%%%%%%%%%%%%%%%%%%%%%%%%%%%%%%%%%%%%%%%%%%%%%%%%%%%%%%%%%%%%%%%%%%%%%%%%
%%%%%%%%%%%%%%%%%%%%%%%%%%%%%%%%%%%%%%%%%%%%%%%%%%%%%%%%%%%%%%%%%%%%%%%%%

The next Proposition has a weaker hypothesis (more general sets $L$)
and also a weaker conclusion (the extremal lattices are not described) 
than Theorem 2.1. 
We say that a {\it{two-dimensional point lattice $\Lambda \subset {\Bbb{R}}^n$
projects orthogonally on the $x_1x_2$-coordinate plane injectively}}, 
if for distinct $\lambda _1, \lambda _2 \in \Lambda $ their images by
this projection are also distinct. 

%%%%%%%%%%%%%%%%%%%%%%%%%%%%%%%%%%%%%%%%%%%%%%%%%%%%%%%%%%%%%%%%%%%%%%%%%
%%%%%%%%%%%%%%%%%%%%%%%%%%%%%%%%%%%%%%%%%%%%%%%%%%%%%%%%%%%%%%%%%%%%%%%%%

\proclaim{Proposition 2.4}
Let $L \subset {\Bbb R}^3$ be the union of a
string of unit balls in ${\Bbb R}^3$, whose centres lie on the $z$-axis.
Let the distance of the neighbouring centres of
balls in $L$ be always at least $2$, and let the
average distance of neighbouring centres of unit balls on the $z$-axis exist
uniformly and equal $2d$ in the sense of \thetag{2.1}, where $1 \le
d \le {\sqrt{2}}$. 
Then any packing of
translates of $L$ in ${\Bbb{R}}^3$
(i.e., any packing of unit balls in ${\Bbb{R}}^3$,
consisting of entire translates of $L$)
has a density at most the supremum of the densities of those two-dimensional 
lattice packings $\{ L + \lambda \mid \lambda \in \lambda \} $
of translates of $L$ in ${{\Bbb R}}^3$, for which the following holds.
The point lattice $\Lambda $ projects
orthogonally to the $xy$-plane injectively, onto a two-dimensional point
lattice in the $xy$-plane.
\endproclaim

%%%%%%%%%%%%%%%%%%%%%%%%%%%%%%%%%%%%%%%%%%%%%%%%%%%%%%%%%%%%%%%%%%%%%%%%%%%%
%%%%%%%%%%%%%%%%%%%%%%%%%%%%%%%%%%%%%%%%%%%%%%%%%%%%%%%%%%%%%%%%%%%%%%%%%%%%

Now we turn to ${\Bbb{R}}^n$. We again need some notations.
Let $L$ be the union of a packing of unit balls in ${\Bbb R}^n$, 
whose centres lie on the 

\newpage

$x_3 \ldots x_n$-coordinate plane.
Further (as a generalization of (2.1)), let us suppose
$$
\cases
{\text{for any }} (0,0,x_3, \ldots ,x_n) {\text{ in the }}
x_3 \ldots x_n{\text{-coordinate hyperplane}} \\
{\text{ and for }} R \to \infty
{\text{ the quotient of the number of ball centres from }}
L \\
{\text{ in }} B\left( (0,0,x_3, \ldots ,x_n),R \right) {\text{ and of }} 
R^{n-2} 
{\text{ tends to some positive}}\\ 
{\text{number, and this convergence is uniform for all }} 
(0,0,x_3, \ldots, x_n) \\
{\text{in the }} x_3 \ldots x_n{\text{-coordinate hyperplane}}
\endcases
\tag 2.2
$$

The inequalities in both Theorems 2.1 and 2.3 follow from the next Proposition.

%%%%%%%%%%%%%%%%%%%%%%%%%%%%%%%%%%%%%%%%%%%%%%%%%%%%%%%%%%%%%%%%%%%%%%%%%%%%%
%%%%%%%%%%%%%%%%%%%%%%%%%%%%%%%%%%%%%%%%%%%%%%%%%%%%%%%%%%%%%%%%%%%%%%%%%%%%%

\proclaim{Proposition 2.5}
Let $n \ge 2$, and let $L$ be the union of a 
packing of translates of $B^n$, with centres in the $x_3 \ldots
x_n$-coordinate plane, such that 
the concentric balls of radius $\sqrt{2}$ form a
covering of the $x_3 \ldots x_n$-coordinate plane. Let us suppose that 
{\rm{(2.2)}} holds.
Then any packing of translates of $L$ in ${\Bbb{R}}^n$ 
(i.e., any packing of unit balls in ${\Bbb{R}}^n$,
consisting of entire translates of $L$)
has a density at most the supremum of the densities of those two-dimensional 
lattice packings $\{ L + \lambda \mid \lambda \in \Lambda \} $
of translates of $L$ in ${{\Bbb R}}^n$,
for which the following holds. The point lattice $\Lambda $ projects
orthogonally to the $x_1x_2$-coordinate 
plane injectively, onto a two-dimensional point
lattice in the $x_1x_2$-coordinate plane.
\endproclaim

%%%%%%%%%%%%%%%%%%%%%%%%%%%%%%%%%%%%%%%%%%%%%%%%%%%%%%%%%%%%%%%%%%%%%%%%%%
%%%%%%%%%%%%%%%%%%%%%%%%%%%%%%%%%%%%%%%%%%%%%%%%%%%%%%%%%%%%%%%%%%%%%%%%%%

We say that $X \subset {\Bbb R}^n$ has {\it{rotational symmetry about the 
$x_3 \ldots x_n$-coordinate plane}} (of dimension $n-2$), if $(x_1,x_2,x_3,
\ldots , x_n) \in X$ implies 
$(x_1 \cos \varphi + x_2 \sin \varphi , - x_1 \sin \varphi $
\newline
$ + x_2 \cos \varphi , x_3, \ldots , x_n) \in X$ for every $\varphi 
\in [0, 2 \pi ]$. 

{\bf{Notation.}} If $X \subset {\Bbb{R}}^n$ has rotational symmetry about the 
$x_3 \ldots x_n$-coordinate plane, and also is open, then let 
$$
\cases
m(X):=
\inf \{ \sqrt{x_1^2+x_2^2}/2 \mid (x_1,x_2,x_3, \ldots , x_n) \in
{\Bbb{R}}^n, \\
X \cap \left( X + (x_1,x_2,x_3, \ldots , x_n) \right) 
= \emptyset \} \,,
\endcases
\tag 2.3
$$ 
and 
$$
M(X):= \sup \{ \sqrt{x_1^2+x_2^2} \mid (x_1,x_2,x_3,
 \ldots , x_n) \in X \} \,.
\tag 2.4
$$
Clearly, $m(X) \le M(X)$. 

%%%%%%%%%%%%%%%%%%%%%%%%%%%%%%%%%%%%%%%%%%%%%%%%%%%%%%%%%%%%%%%%%%%%%%%%%%%
%%%%%%%%%%%%%%%%%%%%%%%%%%%%%%%%%%%%%%%%%%%%%%%%%%%%%%%%%%%%%%%%%%%%%%%%%%%

\proclaim{Proposition 2.6}
Let $L$ be the union of an $(n-2)$-dimensional 
lattice packing of translates of some convex body $K
\subset {\Bbb R}^n$, where $K$ has rotational symmetry about the $x_3 \ldots
x_n$-coordinate plane\, and the corresponding $(n-2)$-dimensional point
lattice lies in the $x_3 \ldots x_n$-coordinate plane. With the 
notations \thetag{2.3} and \thetag{2.4}, let $\left( M(L) =
\right) M(K)=1$ and let $1/{\sqrt{2}} \le m(L)$. 
Then any packing of translates of $L$ in ${\Bbb{R}}^n$
(i.e., any packing of translates of $K$ in ${\Bbb{R}}^n$, consisting of entire

\newpage

translates of $L$) has a density at most the density of the
densest lattice packing of translates of $K$ in ${\Bbb{R}}^n$.
\endproclaim

%%%%%%%%%%%%%%%%%%%%%%%%%%%%%%%%%%%%%%%%%%%%%%%%%%%%%%%%%%%%%%%%%%%%%%%%%%%
%%%%%%%%%%%%%%%%%%%%%%%%%%%%%%%%%%%%%%%%%%%%%%%%%%%%%%%%%%%%%%%%%%%%%%%%%%%

{\definition{Remark 2.7}} To get many examples where Proposition 2.6
can be applied, we consider the following examples.
We suppose $M(K)=1$, and ensure $m(L) \ge 1/{\sqrt{2}}$ by 
letting $ \{ (x_i) \in {\Bbb R}^n \mid \sqrt{x_1^2+x_2^2} \le 1/\sqrt{2} \}
\subset L$.
Let us consider a lattice-tiling on the $x_3 \ldots x_n$-coordinate
plane, by translates of some convex $(n-2)$-polytope $P$, by the vectors of
some $(n-2)$-dimensional lattice $\Lambda $.
Let $K:= \{ (x_1,x_2,x_3, \ldots ,x_n) 
\in {\Bbb R}^n \mid (x_3, \ldots , x_n) \in P,\,\,\,\,
\sqrt{x_1^2+x_2^2} \le f(x_3, \ldots , x_n) \} $, where $f:P \to [1/\sqrt{2},
1]$ is a concave function with maximum $1$. 
Then $K$ is a convex body, rotationally symmetric
about the $x_3 \ldots x_n$-coordinate-plane, 
whose translates by the
vectors in $\Lambda $ have a union satisfying the required properties for $L$.

%%%%%%%%%%%%%%%%%%%%%%%%%%%%%%%%%%%%%%%%%%%%%%%%%%%%%%%%%%%%%%%%%%%%%%%%%%%%%
%%%%%%%%%%%%%%%%%%%%%%%%%%%%%%%%%%%%%%%%%%%%%%%%%%%%%%%%%%%%%%%%%%%%%%%%%%%%%

Observe that in Theorems 2.1, 2.3, and Proposition
2.6 we had packings of lower dimensional
lattices of some convex bodies, where the existence of density was automatic.
On the other hand, in Propositions 2.4, 2.5 we had to take care for the 
densities. 

All the above  statements will turn out to be rather 
simple consequences of the 
following Theorem 2.8. 
Of course, also in this theorem we have to take care for the density. 

For this theorem we have to introduce some
notations.
Let $L \subset {\Bbb{R}}^n$ be rotationally symmetric with
respect to the $x_3 \ldots
x_n$-coordinate plane and be open, with $1/{\sqrt{2}} \le m(L) $ and
$M(L)=1$. 
$$
\cases
{\text{Two translates of }} L, {\text{ say, }} L+(x_1,x_2,x_3, \ldots ,x_n)
{\text{ and }} L+ \\
(y_1,y_2,y_3, \ldots ,y_n) {\text{ are said to }} touch\,\, 
each\,\, other {\text{ if they are disjoint,}} \\
{\text{but for every }}
(x'_1,x'_2) {\text{ and }}
(y'_1,y'_2) {\text{ with }} (y'_1-x'_1)^2 + (y'_2-x'_2)^2 < \\
(y_1-x_1)^2 +
(y_2-x_2)^2 {\text{ we have }}
\left( L+(x'_1,x'_2,x_3, \ldots ,x_n)\right) 
\cap \\
\left( L+(y'_1,y'_2,y_3, \ldots ,y_n) \right) \ne \emptyset \,.
\endcases
\tag 2.5
$$ 
$$
\cases
{\text{We define the function }} g(x_3, \ldots , x_n) {\text{ as the minimal 
value of }} |x_1|/2 \\
{\text{(equivalently: of }} {\sqrt{x_1^2+x_2^2}}/2) {\text{ such that }} L \cap 
\left( L+(x_1,0,x_3, \ldots , x_n) \right) \\
= \emptyset {\text{ (equivalently, such that }} L \cap \left( L + 
(x_1,x_2,x_3, \ldots ,x_n) \right) = \emptyset {\text{).}}
\endcases
\tag 2.6
$$
(This minimum exists for each $(x_3, \ldots ,x_n) \in {\Bbb{R}}^{n-2}$, 
by openness of $L$.)
Clearly for ``touching'' translates of $L$ 
the distance of their axes of rotation (translates of the $x_3 \ldots
x_n$-coordinate plane) is at most $2$. This is sharp: we have 

\newpage

$$
m(L) = \inf g \le \sup g = g(0, \ldots , 0) = M(L) = 1
\tag 2.7
$$

For $L \subset {\Bbb{R}}^n$ with the above properties we consider
some three translates of $L$ mutually ``touching'' each other in the sense of
(2.5) (provided these exist).
Observe that here the $3$'rd, $\ldots $, $n$'th coordinates of the translation
vectors $(x_1,x_2,x_3, \ldots ,$
\newline
$x_n)$ and $(x_1,x_2,x_3, \ldots ,x_n)$ are fixed.
This condition determines the distances of the axes of rotation of 
these three translates of $L$. 
$$
\cases
{\text{Let }} V_0(L) {\text{ denote the infimum of
the areas of the triangles (provided}} \\
{\text{these exist), 
with vertices the points of intersection of the axes of}} \\
{\text{rotation of
these three mutually touching translates of }} L 
{\text{ with the}} \\ 
x_1x_2{\text{-coordinate
plane. Here the
infimum is taken for all three}} \\
{\text{translates of }} L {\text{ with }}
3{\text{'rd}}, 
\ldots , d{\text{'th
coordinates of the respective}} \\
{\text{translation vectors chosen arbitrarily 
in the}} 
x_3 \ldots x_d{\text{-coordinate plane.}} 
\endcases
\tag 2.8
$$
(Later it will turn out that in our
paper the triangles in question will always exist.)
We bound the  densities of packings of translates of $L$ in ${\Bbb R}^n$ from
above. 

For $m(L) \ge 1/{\sqrt{2}} > 0$ and a two-dimensional lattice packing $\{ L
+ \lambda \mid \lambda \in \Lambda \} $ of translates of $L$,
naturally no non-zero lattice vector of $\Lambda $ is 
in the $x_3 \ldots x_n$-coordinate plane, and even in its open
$m(L)$-neighbourhood, i.e., the
point lattice $ \Lambda $ projects orthogonally
to the $x_1x_2$-plane injectively onto a
two-dimensional point lattice in the $x_1x_2$-plane.

We need a generalization of hypotheses (2.1) and (2.2) about density.
$$
\cases 
{\text{Let }} L \subset {\Bbb{R}}^n {\text{ be rotationally symmetric with
respect}} \\
{\text{to the }} 
x_3 \ldots x_n{\text{-coordinate plane and be open, with }} \\
1/{\sqrt{2}} \le m(L) {\text{ and }} M(L)=1 .
\endcases
\tag 2.9
$$
Then $L$ lies in the cylinder $L_0$ 
with base a circle of radius $1$ and centre the origin in the
$x_1x_2$-cooordinate plane, and with axis the $x_3 \ldots x_n$-coordinate
plane (and here the radius $1$ cannot be replaced by any smaller number).
Then it makes sense to speak about the density of $L$ with 
respect to $L_0$ (if it exists). {\it{We write $V$ for volume}} 
(Lebesgue measure).
Let us suppose

\newpage

$$
\cases
{\text{for any }} (0,0,x_3, \ldots ,x_n) {\text{ in the }}
x_3 \ldots x_n{\text{-coordinate hyperplane,}} \\
{\text{ and for }} R \to \infty ,{\text{ the quotient }}
V \left( \{ (\xi _1,\xi _2, \xi _3, \ldots ,\xi _n) \in L \right. \\
\left. \mid {\sqrt{
(\xi _3 -x_3)^2 + \ldots (\xi _n -x_n)^2}} \le R \} \right) \\
/
V \left( \{ (\xi _1,\xi _2, \xi _3, \ldots ,\xi _n) \in L_0 \mid {\sqrt{
(\xi _3 -x_3)^2 + \ldots (\xi _n -x_n)^2}} \le R \} \right) \\
{\text{ tends to some positive 
number }} d(L), {\text{ and this convergence is uniform}} \\
{\text{for all }} 
(0,0,x_3, \ldots, x_n) 
{\text{ in the }} x_3 \ldots x_n{\text{-coordinate hyperplane.}}
\endcases
\tag 2.10
$$

%%%%%%%%%%%%%%%%%%%%%%%%%%%%%%%%%%%%%%%%%%%%%%%%%%%%%%%%%%%%%%%%%%%%%%%%%%%%%
%%%%%%%%%%%%%%%%%%%%%%%%%%%%%%%%%%%%%%%%%%%%%%%%%%%%%%%%%%%%%%%%%%%%%%%%%%%%%

\proclaim{Theorem 2.8}
With the above notations, suppose the above hypotheses \thetag{2.9}
and \thetag{2.10} about $L \subset {\Bbb{R}}^n$.
Then the supremum of the densities of all packings of translates of $L$ in
${\Bbb{R}}^n$ is the supremum of the densities of the two-dimensional lattice
packings of translates of $L$ --- the corresponding point lattice orthogonally
projecting to the $x_1x_2$-coordinate plane injectively onto a
two-dimensional point lattice in the $x_1x_2$-coordinate plane --- spanned by
three mutually touching translates of $L$ (in the
sense of \thetag{2.5}). The maximal density of
two-dimensional lattice packings of translates of $L$, with corresponding
point lattice as decribed above, is at least
$\pi d(L) /(2\sqrt{4m(L)^2-1})$.
\endproclaim

%%%%%%%%%%%%%%%%%%%%%%%%%%%%%%%%%%%%%%%%%%%%%%%%%%%%%%%%%%%%%%%%%%%%%%%%%%%%%
%%%%%%%%%%%%%%%%%%%%%%%%%%%%%%%%%%%%%%%%%%%%%%%%%%%%%%%%%%%%%%%%%%%%%%%%%%%%%

Theorem 2.8 is itself a rather straightforward consequence of the following
general
theorem. Before stating it, we recall the definition of an 
$(r,R)$-system of points in ${\Bbb R}^n$, where $0 < r < R < \infty $. 
We say that $P = \{ p_1, p_2, \ldots \} \subset {\Bbb
R}^n$ {\it{is an $(r,R)$-system}}, 
if the balls with these centres and radii $r$ form
a packing in ${\Bbb R}^n$, 
and the balls with these centres and radii $R$ form a covering
of ${\Bbb R}^n$.

%%%%%%%%%%%%%%%%%%%%%%%%%%%%%%%%%%%%%%%%%%%%%%%%%%%%%%%%%%%%%%%%%%%%%%%%%%%%%
%%%%%%%%%%%%%%%%%%%%%%%%%%%%%%%%%%%%%%%%%%%%%%%%%%%%%%%%%%%%%%%%%%%%%%%%%%%%%

\proclaim{Theorem 2.9}
Let $P \subset {\Bbb R}^2$ be an $(r,R)$-system in ${\Bbb
R}^2$. Let $R/r \le 2{\sqrt{2}}$. Let us consider a Delone (Delaunay)
triangulation of ${\Bbb R}^2$ associated to $P$ (if there are Delone-polygons
with more than three sides, we triangulate them arbitrarily). Then the Delone
triangles can be grouped in groups, such that a group can be
\newline
(1) a single non-obtuse triangle $T$, or
\newline
(2) an obtuse triangle $T$ of area at least $2r^2$, or
\newline
(3) can consist of one non-obtuse triangle $T$ and one, two or three obtuse
triangles $T_i$, whose longest sides coincide with some sides of $T$,
these two, three or four triangles having 
an average area at least $2r^2$, or
\newline
(4) can consist of an obtuse triangle $T$ and one or two further obtuse
triangles $T_i$, whose longest sides coincide with some sides of $T$ but not
with the longest side of $T$, these two or three triangles having an average
area at least $2r^2$.

In particular, if the infimum of the areas of the non-obtuse Delone 
triangles is $V_0$, then
the average area of the Delone triangles is at least $\min \{ V_0, 2r^2 \} $.
This bound is sharp.
\endproclaim

%%%%%%%%%%%%%%%%%%%%%%%%%%%%%%%%%%%%%%%%%%%%%%%%%%%%%%%%%%%%%%%%%%%%%%%%%%%%
%%%%%%%%%%%%%%%%%%%%%%%%%%%%%%%%%%%%%%%%%%%%%%%%%%%%%%%%%%%%%%%%%%%%%%%%%%%%

\newpage

{\bf{Conjecture 2.10.}} The next distance that can occur after $2$ and
$2\sqrt{2}$ inside one (hexagonal, square, or $2$ by $2\sqrt{2}$ rectangular)
layer or between different such layers is $2\sqrt{8/3}$. 
\newline (A)
One can conjecture that for some $\varepsilon >0$, for 
$d \in (\sqrt{2}, \sqrt{2} + \varepsilon ]$ a densest packing of strings with
distance $2d$ is obtained in the following way. Recall from Theorem 2.1 the
densest packing of strings equidistant with distance $2{\sqrt{2}}$.
This can be obtained as follows. We take a homogeneous
point lattice in ${\Bbb{R}}^3$ with
basis $ \{ (2,0,0), (0, 2, 0), (0,0, 2{\sqrt{2}})$, whose basic cell is a
$2 \times 2 \times 2{\sqrt{2}}$ rectangular box, say $B$. 
We add to the points of this
lattice all centres of lattice translates of $B$, obtaining
this way a new point lattice $\Lambda $. This is the point lattice
corresponding to the
densent lattice packing of unit balls in ${\Bbb{R}}^3$. In particular, the 
centres of the neighbours of the ball with centre $(0,0,0)$ are the following
twelve points: $\{ (\pm 2, 0,0), (0, \pm 2,0), (\pm 1, \pm 1, \pm {\sqrt{2}})$ 
(in the last case the plus-minus signs are independent of each other). Now
let us apply the following linear transformation 
$T$ to our lattice $\Lambda $:
$T(1,0,0):=(1,0,0)$, $T(0,1,0):=(0,1,0)$ and $T(0,0,{\sqrt{2}}):=(0,0,d)$, with
$d \in (\sqrt{2}, \sqrt{2} + \varepsilon ]$ and $\varepsilon >0$ being
sufficiently small.
Then the norms of all lattice
vectors in the $xy$-plane are preserved, and the norms of all lattice vectors 
not in the $xy$-plane become longer. This means that the minimum vectors of
the lattice $T \Lambda $ 
are $\{ (\pm 2, 0,0), (0, \pm 2,0) \} $, and the norms of
$T(\pm 1, \pm 1, \pm {\sqrt{2}})$ 
are just a bit larger than $2$, while the norms of
all other non-zero 
lattice vectors are greater than some number strictly greater than
$2$. Thus 
$$
\cases
T\Lambda {\text{ cannot be the point lattice corresponding to a}} \\
{\text{densest lattice packing of the equidistant string of}} \\
{\text{balls with distances of neighbouring ball centres $2d$.}}
\endcases
\tag 2.11
$$
Namely, let us apply to the point lattice
$T\Lambda $ another linear transformation $S$, defined as follows. We let
$S(0,0,1):=(0,0,1)$, and $S(1,0,0):= (\cos \delta , \sin \delta , 0)$ and 
and $S(0,1,0):= (\sin \delta , \cos \delta , 0)$, where $\delta $ is small
(depending on $\varepsilon $). 
Then $S$ carries the rectangular box $TB$ to a right prism $STB$
with base a rhomb with edge lengths $2$ (which is close to a 
$2 \times 2$ square) and with height $2d$. Then the linear map $S$ carries the
minimum vectors $(\pm 2, 0,0)$ and $(0, \pm 2,0)$ of the point lattice 
$T\Lambda $ to vectors of length $2$,
which we want to be minimum vectors of the lattice $ST \Lambda $. Thus, for
$\delta > 0$ sufficiently small (depending on $\varepsilon $) there arise no new
minimum vactors, and $V(STB) < V(TB)$, showing \thetag{2.11}.
For $\varepsilon > 0$ sufficiently
small, we want to investigate the images of $(\pm 1, \pm 1, \pm {\sqrt{2}})$ 
(said otherwise, the vectors pointing from the centre of $B$ to its vertices)
by $ST$, to see how large $\delta > 0$ 
can be so that still these would have lengths
at least $2$. By symmetry in the $z$ coordinate, it suffices to investigate
the images of $(\pm 1, \pm 1, {\sqrt{2}})$. By central symmetry of the base
rhomb 

\newpage

it suffices to investigate only two neighbouring vertices of this
rhomb. Actually, the vertex of the rhomb with an acute angle (e.g., $(0,0,0)$)
got farther 
from the string containing as one unit ball the one with centre the centre of
$STB$. So there no new touching can occur. However, the vertex of the rhomb
with an obtuse angle (e.g., $2( \cos \delta , \sin \delta , 0)$) 
got closer to the above mentioned string. Namely, its distance to the axis of
rotation of our string, i.e., to the centre of the rhomb, decreased to
$(\cos \delta - \sin \delta ){\sqrt{2}}$. We want to have touching of the
strings containing the unit ball with centre at $(0,0,0)$ and at the centre of
$STB$, i.e.,
$$
4 = 2(\cos \delta - \sin \delta )^2 + d^2\,, \,\,{\text{ i.e., }}\,\, 
2 \sin (2 \delta ) = d^2 -2 \,.
\tag 2.12
$$
In this case, the set of the intersection points of the axes of rotation of
the strings with the $xy$-plane forms the vertices and centres of a lattice of
the above rhombs. This is a $(\cos \delta + \sin \delta ) {\sqrt{2}} \times 
(\cos \delta - \sin \delta ) {\sqrt{2}}$ rectangular lattice. {\it{We
conjecture that for $\varepsilon >0$ sufficiently small, the above constructed
packing is the densest 
packing of translates of our string of unit balls, with
centres equidistant with distances $2d \in (2{\sqrt{2}},
2({\sqrt{2}} + \varepsilon )]$.}} Its density is 
$$
(4 \pi /3) / \left( 4d \cos (2 \delta ) \right) \,.
\tag 2.13
$$
\newline (B)
For $d={\sqrt{8/3}}$ we have another picture. We take the non-lattice-like
regular close packing of regular triangular lattices of unit balls with edge
lengths $2$, consisting of layers which are translates of one such layer
with ball centres 
in the $xy$-plane. Then on each vertical line containing some ball centre the
ball centres on them are equidistant, with distance ${\sqrt{8/3}}$. If we take
only
every second layer, then the intersections of the rotation axes of the strings,
which have a ball centre in these every second layers,
with the $xy$-plane, form a lattice of regular triangles of edge length $2$ in
the $xy$-plane, joining along entire edges. This mosaic $\{ 3,6 \} $ is also
the Delone triangulation corresponding to these points in the $xy$-plane.
If we consider the remaining every
second layers, then the intersection points of the axes of rotation of strings
with one ball centre in these other every second layers with the
$xy$-plane will be the centres of all original regular triangles which are
translates of each other (thus only of a ``half'' of all considered
regular triangles in the $xy$-plane). Adding these new points, the triangles
whose centres are not added have three edge-neighbour triangles whose 
centres are added. The triangle with centre not added remains an empty
triangle, but actually its circumcircle passes also through the centres of its
edge-neighbour triangles. Therefore the Delone triangulation of all these
points consists of the vertices of a mosaic $\{ 6,3 \} $, of regular hexagons
with circumradius $2/{\sqrt{3}}$.
\newline (C)
All this shows that for $d$ in the interval $[\sqrt{2}, \sqrt{8/3}]$ the
solution of the 

\newpage

densest packing of translates of 
our strings conjecturably will behave differently in
several subintervals: it begins with a lattice packing, and ends
with a non-lattice-like regular packing.

%%%%%%%%%%%%%%%%%%%%%%%%%%%%%%%%%%%%%%%%%%%%%%%%%%%%%%%%%%%%%%%%%%%%%%%%%%%%
%%%%%%%%%%%%%%%%%%%%%%%%%%%%%%%%%%%%%%%%%%%%%%%%%%%%%%%%%%%%%%%%%%%%%%%%%%%%

{\definition{Remark 2.11}}
Analogously to the conjecture of \cite{Sz} in \S 1, 
for L. Fejes T\'oth's question for arbitrary $d$ one has also the following
construction.
Let $L$ be a string of unit balls in ${\Bbb{R}}^3$, with centres on a straight
line, which are equidistant with distance $2d$, where $d > {\sqrt{2}}$. (The
case $d \in [1, {\sqrt{2}}]$ is settled in Proposition 2.5.) 
We choose a natural number $k$ and consider $k$ closely packed
parallel hexagonal, or square, or $2$ by $2{\sqrt{2}}$ rectangular
layers of unit balls, in the first case for any layer the following layer
packed in any of the possible two ways
(or possibly some other packings of unit balls with
centres in some lattice plane of the densest lattice packing of balls). 
We translate this block of $k$ layers
periodically so that the centre of a ball of the first layer of a block has
a distance $2d$ from the centre of some ball of the first layer of the
preceding block and has distances $2$ from the centres of two balls 
of the
last layer of the preceding block. 
Of course, this construction does not include the example from
Conjecture 2.10 about $d \in ( {\sqrt{2}}, {\sqrt{2}}+ \varepsilon ]$.
\vskip.1cm

%%%%%%%%%%%%%%%%%%%%%%%%%%%%%%%%%%%%%%%%%%%%%%%%%%%%%%%%%%%%%%%%%%%%%%%%%%%%
%%%%%%%%%%%%%%%%%%%%%%%%%%%%%%%%%%%%%%%%%%%%%%%%%%%%%%%%%%%%%%%%%%%%%%%%%%%%

If $2d$ is the distance of the centres of two balls in some close packing of
parallel hexagonal or square or $2$ by $2\sqrt{2}$ rectangular
layers of unit balls (or possibly of some other packings of unit balls from
the densest lattice packing of unit balls, with
centres in some lattice plane), then evidently some periodic system of
closely packed parallel such layers will be a densest
packing. Of course, the problem for any $d$ 
is hopeful to be proved only for small
values of $d$, since for $d \to \infty $ the solution of 
this problem would imply the
theorem of T. C. Hales-S. Ferguson \cite{HF}
on the densest packing of unit balls in ${\Bbb{R}}^3$.

%%%%%%%%%%%%%%%%%%%%%%%%%%%%%%%%%%%%%%%%%%%%%%%%%%%%%%%%%%%%%%%%%%%%%%%%%%%%
%%%%%%%%%%%%%%%%%%%%%%%%%%%%%%%%%%%%%%%%%%%%%%%%%%%%%%%%%%%%%%%%%%%%%%%%%%%%

{\definition{Remark 2.12}} Also for some other problems,
by estimating in some way
the area of the non-obtuse Delone triangles, that depends on the
particular problem, this theorem can be
applied. As an example, we mention the following problem of L. Fejes T\'oth
\cite{FT}.
A {\it{molecule in ${\Bbb R}^2$}} is the union of a unit circle and of one or
two circles of radius $r<1$, such that they form a packing, and the unit
circle touches the circle(s) of radius $r$. The position of two small circles
is not determined, only the packing and touching properties are prescribed.
For $r \le 2/\sqrt{3} - 1 = 0.1547 \ldots $ this problem is not interesting,
since in the holes of a densest lattice packing of unit circles there is enough
space for the small circles. By using 
some further ideas we can determine the densest
packing of two-atom molecules, for $r \le 0.1899 \ldots $. (The so called
$L^*$-decomposition introduced by J. Moln\'ar \cite{M77}, \cite{M78} also
plays a role there.) This is attained, e.g., in the following case. We
consider packings of unit circles, which consist of horizontal strings of
touching copies of the unit circle. We enumerate these horizontal strings by the
integers, so that the $(i+1)$'st string 
is the upper neighbour of the $i$'th string.
The neighbouring strings of unit circles
touch. The $2i$'th and $(2i+1)$'st strings are closely packed. The $(2i+1)$'st

\newpage

and $(2i+2)$'nd strings are not closely packed, but so that two neighbouring
centres of the $(2i+1)$'st string of circles and some centre of the $(2i+2)$'nd
string of circles form
isosceles triangles with two sides of length $2$,
and with circumradius $1+r$. The small circles have centres exactly in the
circumcentres of these last mentioned isosceles triangles. For three-atom
molecules, for $r - (2/\sqrt{3} - 1)$ sufficiently small, a densest packing
consists of horizontal strings of touching unit circles, where any two 
neighbouring horizontal strings join to each other as the $(2i+1)$'st
and $(2i+2)$'nd strings above.
We hope to return to this problem later.

%%%%%%%%%%%%%%%%%%%%%%%%%%%%%%%%%%%%%%%%%%%%%%%%%%%%%%%%%%%%%%%%%%%%%%%%%%%%
%%%%%%%%%%%%%%%%%%%%%%%%%%%%%%%%%%%%%%%%%%%%%%%%%%%%%%%%%%%%%%%%%%%%%%%%%%%%

\head{\S 3 Preparatory lemmas}\endhead

In this paragraph we will prove statements in ${\Bbb R}^2$.
We write $\Delta ABC$ for the triangle $ABC$ in ${\Bbb R}^2$, 
and denote the area of a set in ${\Bbb R}^2$ by $V( \cdot )$. 
For a triangle $\Delta ABC$ 
in ${\Bbb R}^2$ the sides opposite to $A,B,C$ will be
denoted by $a,b,c$, and 
the angles at $A,B,C$ by $\alpha , \beta , \gamma $.

We will frequently use the following elementary fact. 
$$
\cases
{\text{For\,\,non-obtuse\,\,triangles\,\,the\,\,area\,\,is\,\,a\,\,strictly}} \\
{\text{monotonous\,\,function\,\,of\,\,the\,\,side-lengths.}}
\endcases
\tag 3.0
$$
In fact, if 
the side lengths are $a \le a'$, $b \le b'$, $c \le c'$ then by inflation
we may attain, with the evident notations, that, e.g., $c=c'$ 
(and $A=A'$ and $B=B'$). Then $C'$ has a distance from $A$ and $B$ at least
$b$ and $a$, and projects orthogonally to a point of the
side $AB$, so its distance to the line $AB$ is minimal if and only if $C'=C$.

%%%%%%%%%%%%%%%%%%%%%%%%%%%%%%%%%%%%%%%%%%%%%%%%%%%%%%%%%%%%%%%%%%%%%%%%%%%
%%%%%%%%%%%%%%%%%%%%%%%%%%%%%%%%%%%%%%%%%%%%%%%%%%%%%%%%%%%%%%%%%%%%%%%%%%%

\proclaim{Lemma 3.1}
Let $p,c \in (0, \infty )$ with $c \ge p$ be fixed. Let 
in a triangle $\Delta ABC$ the side $c$ be fixed, the vertex $C$ be variable,
and let $\alpha \le \pi /2$, $\beta \le \pi /2$, $a \ge p$, $b  \ge p$
and $\gamma \le \gamma _0$.
If under these conditions the triangle $\Delta ABC$ 
has a minimal area, then in two of
the above five last inequalities we have equality.

Let us assume, in addition to the above hypotheses, that
also $\gamma \ge \pi /2$.
If under these new conditions the triangle $\Delta ABC$ 
has a minimal area, then in
two of the inequalities $a \ge p$, $b  \ge p$ and $\gamma \le \gamma _0$ we
have equalities.
\endproclaim

%%%%%%%%%%%%%%%%%%%%%%%%%%%%%%%%%%%%%%%%%%%%%%%%%%%%%%%%%%%%%%%%%%%%%%%%%%
%%%%%%%%%%%%%%%%%%%%%%%%%%%%%%%%%%%%%%%%%%%%%%%%%%%%%%%%%%%%%%%%%%%%%%%%%%%

\demo{Proof} 
By $\alpha , \beta \le \pi /2$ we have that our triangle lies in a parallel
strip $S$ 
with boundary lines passing through $A,B$ and orthogonal to the side $AB$.
We may assume $AB$ horizontal, and that $C$ is above the side $AB$ (the
hypotheses imply that $C$ lying on the line $AB$ is impossible). The
inequalities $a \ge p$, $b \ge p$ and $\gamma \le \gamma _0$
express that the vertex $C \in S$ does not
lie strictly below the graphs of certain
concave functions (for $\gamma \le \gamma _0$ we have
to distinguish the cases of obtuse and non-obtuse $\gamma _0$). If $C$ does
not lie on any of these graphs, and also does not lie on the boundary lines of
$S$, then we can move it vertically downwards, decreasing $V(\Delta ABC)$.
If $C$ lies only on one of these graphs, or
boundary lines, then on this graph, or boundary line 
we can move $C$ either increasing or decreasing its $x$-coordinate, or
decreasing its $y$-coordinate, decreasing $V(\Delta ABC)$.

\newpage

Let us have also the additional hypothesis $\gamma \ge \pi /2$. Then
$\alpha , \beta < \pi /2$ holds automatically, which excludes that $C$ lies on
any of the boundary lines of $S$. 
Moreover, we may assume $\gamma _0 \ge \pi /2 $.
As in the proof of the first part, we move $C$
vertically downwards (thus preserving $\gamma \ge \pi /2$), 
until $C$ gets to some of the curves given by $a=p$,
$b=p$ and  $\gamma = \gamma _0$. If this curve ${\Cal{C}}$ is the one given by 
$\gamma = \gamma _0$, then $\gamma \le \pi /2$ is automatically satisfied, and
the minimal area occurs when one of $a$ and $b$ equals $p$. 
If this curve ${\Cal{C}}$ is the one given by, e.g., $b=p$, then rotating $C$
about $A$ towards $B$, thus along ${\Cal{C}}$, 
the area decreases, and $\gamma \ge \pi /2$ will be preserved.
$\blacksquare $
\enddemo

%%%%%%%%%%%%%%%%%%%%%%%%%%%%%%%%%%%%%%%%%%%%%%%%%%%%%%%%%%%%%%%%%%%%%%%%%%
%%%%%%%%%%%%%%%%%%%%%%%%%%%%%%%%%%%%%%%%%%%%%%%%%%%%%%%%%%%%%%%%%%%%%%%%%%%

\proclaim{Lemma 3.2}
Let $a_0,b_0, R_0 \in (0, \infty )$.
Let in the triangle $\Delta ABC$ 
there hold $\gamma \ge \pi /2$, $a \ge a_0$, $b \ge
b_0$, and the circumradius $R$ be at most $R_0$. If under these conditions 
$\Delta ABC$ has a minimal area, then $a=a_0$, $b=b_0$ and $R=R_0$. In
particular, for 
\newline
{\rm{(1)}}
$a_0=p \sqrt{2}$,
$b_0=p$ and $R_0=p \sqrt{2}$ this minimal area is $p^2 
(\sqrt{7} + \sqrt{3})/8 = p^2 \cdot 0.5472 \ldots $\,,
\newline {\rm{(2)}}
$a_0=p \sqrt{2}$, $b_0=p$ and $R_0=p\sqrt{5/2}$ 
this minimal area is $p^2/2$.
\newline {\rm{(3)}} 
$a_0=b_0=p$ and $R_0=p \sqrt{2}$ this minimal area is $p^2\sqrt{7}/8 =
p^2 \cdot 0.3307 \ldots $\,.
\endproclaim

%%%%%%%%%%%%%%%%%%%%%%%%%%%%%%%%%%%%%%%%%%%%%%%%%%%%%%%%%%%%%%%%%%%%%%%%%%
%%%%%%%%%%%%%%%%%%%%%%%%%%%%%%%%%%%%%%%%%%%%%%%%%%%%%%%%%%%%%%%%%%%%%%%%%%%

\demo{Proof} 
Fixing $\gamma $, decrease $a$ to $a_0$, $b$ to $b_0$. Then also $c$
and $R=c/(2 \sin \gamma )$ decrease. After this, fixing $a = a_0$ and $b=b_0$
we increase $\gamma $ till $R=c/(2 \sin \gamma )$ increases (since $c$
increases and $\gamma \ge \pi /2$ increases) to $R_0$. Thus we
have strictly decreased $V(\Delta ABC)$, unless we have equalities in all
three asserted inequalities of the lemma.
$\blacksquare $
\enddemo

%%%%%%%%%%%%%%%%%%%%%%%%%%%%%%%%%%%%%%%%%%%%%%%%%%%%%%%%%%%%%%%%%%%%%%%%%%
%%%%%%%%%%%%%%%%%%%%%%%%%%%%%%%%%%%%%%%%%%%%%%%%%%%%%%%%%%%%%%%%%%%%%%%%%%%

For a triangle $\Delta A_1BC$ we denote the sides opposite to $A_1,B,C$ by
$a_1,b_1,c_1$, and the angles at these vertices by $\alpha _1 , \beta _1,
\gamma _1$. Analogous notations will be used with the indices $2$ and $3$, for
the triangles $\Delta AB_2C$ and $\Delta ABC_3$.

%%%%%%%%%%%%%%%%%%%%%%%%%%%%%%%%%%%%%%%%%%%%%%%%%%%%%%%%%%%%%%%%%%%%%%%%%%
%%%%%%%%%%%%%%%%%%%%%%%%%%%%%%%%%%%%%%%%%%%%%%%%%%%%%%%%%%%%%%%%%%%%%%%%%%%

\proclaim{Lemma 3.3} (Cf. also \cite{Ka}, proof of Theorem 14, formula 10.)
Let the triangles $\Delta ABC$ and $\Delta A_1BC$ 
lie on the opposite sides of side $BC$, and
let $\alpha + \alpha _1 \le \pi$. Fixing $b,c,b_1,c_1$ decrease $\alpha $ and
$\alpha _1$. Then the sum of the areas of $\Delta ABC$ and $\Delta A_1BC$
decreases.
\endproclaim

%%%%%%%%%%%%%%%%%%%%%%%%%%%%%%%%%%%%%%%%%%%%%%%%%%%%%%%%%%%%%%%%%%%%%%%%%%
%%%%%%%%%%%%%%%%%%%%%%%%%%%%%%%%%%%%%%%%%%%%%%%%%%%%%%%%%%%%%%%%%%%%%%%%%%%

\demo{Proof} 
We have $b^2+c^2 -2bc \cos \alpha = b_1^2+c_1^2 -2b_1c_1 \cos \alpha _1$,
hence $d \alpha _1 / d \alpha = (bc \sin \alpha ) / (b_1c_1 \sin \alpha _1)$.
Thus 
$$
{\frac{d}{d\alpha }}( bc \sin \alpha + b_1c_1 \sin \alpha _1) = bc \sin \alpha
( \cot \alpha + \cot \alpha _1 ) \ge 0\,.
\tag 3.3.1
$$
$\blacksquare $
\enddemo

%%%%%%%%%%%%%%%%%%%%%%%%%%%%%%%%%%%%%%%%%%%%%%%%%%%%%%%%%%%%%%%%%%%%%%%%%%
%%%%%%%%%%%%%%%%%%%%%%%%%%%%%%%%%%%%%%%%%%%%%%%%%%%%%%%%%%%%%%%%%%%%%%%%%%%

\proclaim{Lemma 3.4}
Let the triangles $\Delta ABC$ and $\Delta A_1BC$ 
lie on the opposite sides of their common side $BC$, and
let $\alpha = \alpha _1 = \pi /2$. Fixing $b,c,b_1,c_1$ decrease $\alpha $ and
$\alpha _1$ a bit, preserving $c_1 < a$. After this rotate side $A_1B$ about
$B$, decreasing $\beta _1$, until $\alpha _1 $ increases to $\pi /2$. Then 
$V(\Delta ABC) + V(\Delta A_1BC)$ decreases. 
\endproclaim

%%%%%%%%%%%%%%%%%%%%%%%%%%%%%%%%%%%%%%%%%%%%%%%%%%%%%%%%%%%%%%%%%%%%%%%%%%
%%%%%%%%%%%%%%%%%%%%%%%%%%%%%%%%%%%%%%%%%%%%%%%%%%%%%%%%%%%%%%%%%%%%%%%%%%%

\newpage

\demo{Proof} 
The firt motion decreases the total area by Lemma 3.3. The second motion is
possible if the first motion was sufficiently small, and during this second
motion $V(\Delta BA_1C) = (1/2) a_1c_1 \sin \beta _1$ decreases (and $V(\Delta
ABC)$ remains constant). 
$\blacksquare $
\enddemo

%%%%%%%%%%%%%%%%%%%%%%%%%%%%%%%%%%%%%%%%%%%%%%%%%%%%%%%%%%%%%%%%%%%%%%%%%%
%%%%%%%%%%%%%%%%%%%%%%%%%%%%%%%%%%%%%%%%%%%%%%%%%%%%%%%%%%%%%%%%%%%%%%%%%%%

\proclaim{Lemma 3.5}
Let $p \in (0, \infty )$.
Let the triangles $\Delta ABC$ and $\Delta A_1BC$ 
lie on the opposite sides of their common
side $BC$,
and let $c_1 = b_1 = p$, $p \le b \le p \sqrt{7/2}$, $\gamma \ge  \pi
/2$, $\alpha _1 \ge \pi /2$, $\alpha + \alpha _1 \le \pi $, and the
circumradius of $\Delta ABC$ be $p \sqrt{2}$. Then under these hypotheses
the sum of the areas of our triangles is minimal only 
either for $\alpha _1 = \pi /2$ or for $\alpha + \alpha _1 = \pi $.
\endproclaim

%%%%%%%%%%%%%%%%%%%%%%%%%%%%%%%%%%%%%%%%%%%%%%%%%%%%%%%%%%%%%%%%%%%%%%%%%%
%%%%%%%%%%%%%%%%%%%%%%%%%%%%%%%%%%%%%%%%%%%%%%%%%%%%%%%%%%%%%%%%%%%%%%%%%%%

\demo{Proof} 
Let 
$$
\delta :=\alpha _1 /2 \in (\pi /4, \pi /2) {\text{\,\, and \,\,}} 
\beta _0 = \arcsin \left( 1/\sqrt{8} \right) = 20.7048 \ldots ^{\circ }\,.
\tag 3.5.1
$$
We have $a/2 = p \sqrt{2} \sin \alpha $ and $a_1/2 = p \sin \delta $.
But $a = a_1$, hence 
$$
\sin \delta = \sqrt{2} \cdot \sin \alpha \,.
\tag 3.5.2
$$ 
Also, by $\delta = \alpha _1 /2 \ge \pi /4$ we have
$a = a_1 = 2 p \sin \delta \ge  p \sqrt{2}$ and the circumradius of $\Delta
ABC$ is $p {\sqrt{2}}$, hence 
$$
\alpha \ge \pi/6 \,.
\tag 3.5.3
$$
Next we determine the upper bound for $\alpha $. By increasing $\alpha _1$, 
and thus $\delta $, but retaining the circumcircle of $\Delta ABC$ and $b$
(hence also $\beta $), also
$a=a_1=2 p \sin \delta < 2p$ increases. Then the line $BC$ gets closer to the
circumcentre of $\Delta ABC$, but does not pass over the circumcentre, since
$a < 2p$, while the diameter of the circumcircle is $2p \sqrt{2}$, thus is
larger than $a$. Also the smaller half central angle $\alpha $
subtended by the side $AB$ increases. At the
same time the distance of $A_1$ to the line $BC$, i.e., $p \cos \delta $ 
decreases, so $A_1$ gets ever
closer to the circumcircle. For 
a certain value of $\alpha _1$, the point $A_1$ will lie on
the boundary of the circumcircle, and then we cannot increase $\alpha _1$ this
way any
more (by $\alpha + \alpha _1 \le \pi $). 
Then the quadrangle $BA_1CA$ has a circumcircle, $K$, say. During this
motion the order of the points $A,B,C$ on $K$ does not change: namely, the
smaller central angles subtended both by the sides $BC$ and $CA$ 
remain less than $\pi $. Moreover, the smaller
central angles subtended by the sides $BA_1$ and $A_1C$ in $K$
(which are half the smaller central
angle subtended by $AB$ in $K$) will be $2 \arcsin \left(
(p/2) / (p \sqrt{2} ) \right) = 2 \arcsin (1/\sqrt{8})$. Then $\alpha $ 
is the smaller central angle subtended by the chord $BA_1$ (or $A_1C$) of
length $p$, 
therefore after this increasing $\alpha $ becomes 
$2 \arcsin (1/\sqrt{8}) = 2 \beta _0$.
That is, 
$$
\pi /6 \le \alpha \le 2 \beta _0 \,.
\tag 3.5.4
$$

Similarly, by $b \ge p$ we have that $\beta $ is half the smaller 
central angle 

\newpage

subtended in $K$ by a chord of length in $[p, p{\sqrt{7}}/2]$. 
That is, $\beta $ is at least half the smaller 
central angle subtended in $K$ by a chord of length $p$, i.e.,
$\beta \ge \arcsin \beta _0$. Similarly, $\beta $ is at most half the smaller 
central angle subtended in $K$ by a chord of length $p\sqrt{7/2}$, i.e.,
$\beta \le \arcsin ({\sqrt{7}}/4) = 2 \beta _0$.
That is, 
$$
\beta _0 \le \beta \le 2 \beta _0\,.
\tag 3.5.5
$$
From the upper bounds for $\alpha , \beta $ in \thetag{3.5.4}
and \thetag{3.5.5} we get
$$
\gamma \ge \pi - 4 \beta _0 > \pi /2\,.
\tag 3.5.6
$$
By \thetag{3.5.2} we have 
$$
d \delta / d \alpha =
\sqrt{2} \cdot \cos \alpha / \cos \delta \,.
\tag 3.5.7
$$
The sum $S$ of the areas of the two
triangles $\Delta ABC)$ and $BA_1C$ is 
$$
\cases
S = ab \sin (\alpha + \beta ) /2 + p^2 \left( \sin (2 \delta ) \right) /2 = \\
2 (p \sqrt{2})^2 \sin \alpha \sin \beta \sin (\alpha + \beta ) / 2 + 
p^2 \left( \sin (2 \delta ) \right) /2 = \\
2 p^2 \sin \beta \left( \cos ( \beta ) - \cos (2 \alpha + \beta ) \right) +
p^2 \left( \sin (2 \delta ) \right) /2\,
\endcases
\tag 3.5.8
$$
that depends on $\alpha $ and $\beta $ (recall that $\delta $ is a function of
$\alpha $ by \thetag{3.5.10}). However, we will consider $\beta $ as fixed, and
then $S$ is a function of $\alpha $ only. We have, also using \thetag{3.5.7},
$$
\frac{dS}{d \alpha } = 4p^2 \sin \beta \sin (\beta + 2 \alpha ) +  
p^2 \left( \cos (2 \delta ) \right)
\sqrt{2} \frac{\cos \alpha }{\cos \delta }\,.
\tag 3.5.9
$$
Once more using \thetag{3.5.2}, we see that 
this has the same sign as 
$$
\cases
T := 2\sqrt{2} \sin \beta \sin (\beta + 2 \alpha ) \cos
\delta +\cos (2 \delta ) \cos \alpha = \\
2\sqrt{2} \sin \beta \sin (\beta + 2 \alpha ) \sqrt{\cos ( 2 \alpha )} -
(4 \sin ^2 \alpha -1) \cos \alpha \,.
\endcases
\tag 3.5.10
$$
Observe that under the square root we have by 
$2 \alpha \le 4 \beta _0 <\pi /2$ a positive number, and by \thetag{3.5.4} 
and \thetag{3.5.5} $\beta _0 + \pi /3 \le \beta + 2 \alpha \le 6 \beta _0$, so
$\sin (\beta + 2 \alpha ) > 0$.

We are going to show that $T$ is strictly monotonically decreasing for $\pi /6
\le \alpha \le 2 \beta _0$ (cf. \thetag{3.5.4}). 
Thus either it has always positive or always
negative sign, or below
some value of $\alpha $ it is positive, and above this value of $\alpha $
it is negative. In any of these three cases the same holds for $dS/d \alpha $,
so 
$$
S {\text{\,\,attains\,\,its\,\,minimum\,\,only\,\,for\,\,}} \alpha = \pi /6 
{\text{\,\,or\,\,for\,\,}} \alpha = 2 \beta _0\,.
\tag 3.5.11
$$
Taking in account \thetag{3.5.2} we have that the
equalities $\pi /6 = \alpha $ 

\newpage

and $\alpha = 2 \beta _0$
are equivalent to 
$\alpha _1 =  \pi /2$, and to $\alpha _1 /2 = \arcsin {\sqrt{7/8}}$ when
simultaneously $\sin (\alpha /2) = \sin \beta _0 = 1/{\sqrt{8}}$ and thus
$\alpha _1 + \alpha = \pi $. 

We write 
$$
T_1 := \sin ^2 (\beta + 2 \alpha ) \cos (2 \alpha ) =
\left( 1 - \cos \left( 2 (\beta + 2 \alpha ) \right) \right) \cos (2 \alpha )/2 
\tag 3.5.12
$$ 
and 
$$
T_2 := (4 \sin ^2 \alpha - 1) \cos \alpha \,.
\tag 3.5.13
$$ 
Then 
$$
T = 2\sqrt{2} \sin \beta {\sqrt{T_1}} - T_2 \,.
\tag 3.5.14
$$ 
Since $\beta $ is considered as fixed, and only $\alpha $
varies, 
$$
\cases
{\text{the\,\,strictly\,\,monotonically\,\,decreasing\,\,property\,\,of\,\,}} T
{\text{\,\,for}} \\
\pi /6 \le \alpha \le 2 \beta _0 
{\text{\,\,follows\,\,from\,\,the\,\,strictly\,\,monotonically}} \\
{\text{decreasing\,\,property\,\,of
\,\,}} T_1 {\text{\,\,and\,\,from\,\,the}} \\
{\text{strictly\,\,monotonically\,\,increasing\,\,property\,\,of\,\,}} T_2, \\
{\text{also\,\,for\,\,}} \pi /6 \le \alpha \le 2 \beta _0\,. 
\endcases
\tag 3.5.15 
$$
We are going to show these strict monotonicity properties.

We have (using for the first term the second form of $T_1$ and for the second
term the first form of $T_1$)
$$
dT_1/d \alpha = 2 \sin \left( 2 (\beta + 2 \alpha ) \right) \cos (2 \alpha ) - 
2 \sin ^2 (\beta + 2 \alpha ) \sin (2 \alpha )\,.
\tag 3.5.16
$$
For $\beta + 2 \alpha > \pi /2$ (also noting $\beta + 2 \alpha \le 6 \beta _0 =
124.2288 \ldots ^{\circ }$) this is negative by 
$\pi /3 \le 2 \alpha \le 4 \beta _0
= 82.8192 \ldots ^{\circ } < \pi /2$. 
For $\beta + 2 \alpha \le \pi /2$ we have by $\pi /3 \le 2 \alpha \le 4 
\beta _0 < \pi /2$ and
$80.7048 \ldots ^{\circ } =
\beta _0 + \pi /3 \le \beta + \pi /3 \le \beta + 2 \alpha \le \pi /2$ that
$$
\cases
dT_1/d \alpha \le 2 \cdot 1 \cdot \cos (\pi /3) - 
2 \sin ^2 (\beta + 2 \alpha ) \cdot \sin (\pi /3) \le \\
1 - \sin ^2 (\beta _0 + \pi /3) \cdot {\sqrt{3}} = 
1 - {\sqrt{3}} \cdot 0.9739 \ldots < 0 \,.
\endcases
\tag 3.5.17
$$ 
Hence $T_1$ strictly monotonically decreases for $\pi /6 \le \alpha
\le 2 \beta _0$.

For $T_2$ we have by $\pi /6 \le \alpha \le 2 \beta _0 = 41.4096 \ldots
^{\circ }$ that
$$
\cases
dT_2/d \alpha = 4 \sin (2 \alpha ) \cdot \cos \alpha - 
(4 \sin ^2 \alpha - 1) \sin \alpha \ge \\
4 \sin ( \pi /3) \cdot \cos (2 \beta _0) -
(4 \sin ^2 (2 \beta _0) - 1) \cdot \sin (2 \beta _0) = \\
2 {\sqrt{3}} \cdot (3/4) - (3/4) \cdot ({\sqrt{7}}/4) > 0\,.
\endcases
\tag 3.5.18
$$

\newpage

Hence $T_2$ strictly monotonically increases for $\pi /6 \le \alpha
\le 2 \beta _0$.

By \thetag{3.5.15} these strict 
monotonicity properties of $T_1$ and $T_2$ imply the strictly
monotonically decreasing property of $T$ for $\pi /6 \le \alpha \le 2 \beta
_0$. This shows \thetag{3.5.11}, that is just the statement of this lemma.
$\blacksquare $
\enddemo

%%%%%%%%%%%%%%%%%%%%%%%%%%%%%%%%%%%%%%%%%%%%%%%%%%%%%%%%%%%%%%%%%%%%%%%%%%
%%%%%%%%%%%%%%%%%%%%%%%%%%%%%%%%%%%%%%%%%%%%%%%%%%%%%%%%%%%%%%%%%%%%%%%%%%%

The following Lemmas 3.6 and 3.8 are related to \cite{BKM}, Lemma 3.1 ---
however, in our lemmas also the circumradius of the polygon is variable.

%%%%%%%%%%%%%%%%%%%%%%%%%%%%%%%%%%%%%%%%%%%%%%%%%%%%%%%%%%%%%%%%%%%%%%%%%%
%%%%%%%%%%%%%%%%%%%%%%%%%%%%%%%%%%%%%%%%%%%%%%%%%%%%%%%%%%%%%%%%%%%%%%%%%%%

\proclaim{Lemma 3.6}
Let a convex pentagon $P$ have a circumcircle and let it contain its
circumcentre. Let each side of $P$ have a length at least $p\,\,(>0)$.
Then the area of $P$ is at least $p^2 \cdot (5/4) \cot (\pi /10)  
= p^2 \cdot 1.7204 \ldots $. Equality holds if and only if $P$ is
a regular pentagon of side length $p$.
\endproclaim

%%%%%%%%%%%%%%%%%%%%%%%%%%%%%%%%%%%%%%%%%%%%%%%%%%%%%%%%%%%%%%%%%%%%%%%%%%
%%%%%%%%%%%%%%%%%%%%%%%%%%%%%%%%%%%%%%%%%%%%%%%%%%%%%%%%%%%%%%%%%%%%%%%%%%%

\demo{Proof} 
We may suppose that $P$ is an extremal pentagon (which exists, since
$V(P)/$
\newline
$p^2$ is similarity invariant).

Let us denote the circumradius of our pentagon $P$ by $R$, and the angles
subtended by the sides of $P$ at its circumcentre $O$ by $\varphi _i$, where
$1 \le i \le 5$. Since $P$ contains $O$, for $1 \le i \le 5$ we have 
$\varphi _i \in [0, \pi]$. 
Let us write
$$
\varphi _0 := 2 \arcsin \left( p/(2R) \right) > 0, {\text{ i.e., }} R=p/ \left(
2 \sin (\varphi _0 / 2) \right) \,.
\tag 3.6.1
$$

Since each side of $P$ has a length at least $p$, for $1 \le
i \le 5$ we have $\varphi \ge \varphi _0 $. Thus 
$$
{\text{for }} 1 \le i \le 5 {\text{ we have }} 
\varphi _i \in [\varphi_0, \pi ] {\text{ and }}
\sum _{i=1}^5 \varphi _i = 2 \pi \,.
\tag 3.6.2
$$
We
have 
$$
\cases
V(P)/p^2 = (R^2/2) \sum _{i=1}^5 \sin \varphi _i /
\left( 2R \sin (\varphi _0 /2) \right) ^2 \\
= \sum _{i=1}^5 \sin \varphi _i / \left( 8
\sin ^2 (\varphi _0 /2) \right) \,.
\endcases
\tag 3.6.3
$$ 
Observe that the function $\sin \varphi $ is a strictly concave function on
$[0, \pi ]$. This implies that if $V(P)/p^2$ is
maximal then 
$$
{\text{there is at most one }}\varphi _i  {\text{ in the open interval }}
(\varphi _0, \pi ) \,.
\tag 3.6.4
$$
In fact if there were two
such ones then preserving their sum we could increase their difference a bit, 
and then $V(P)/p^2$ would decrease.

There cannot be two $\varphi _i$'s equal to $\pi $, since then we would have
$\sum _{i=1}^5 \varphi _i > 2 \pi $. So either 
\newline
(1) there is one $i$ with $\varphi = \pi $, or

\newpage

(2) there is no $i$ with $\varphi = \pi $.

{\bf{1.}}
We begin with the proof of case (1). Let, e.g., $\varphi _5 = \pi $. Then
$$
{\text{for }} 1 \le i \le 4 {\text{ we have }}\varphi _i \in [\varphi _0 , \pi )
{\text{ and }} \sum _{i=1}^4 \varphi _i = \pi \,,
\tag 3.6.5
$$
and we have to minimize $ \sum _{i=1}^5 \sin \varphi _i / \left( 8 \sin ^2 
(\varphi_0 /2) \right) 
= \sum _{i=1}^4 \sin \varphi _i/ 
\left( 8 \sin ^2 (\varphi_0 /2) \right) $.

By \thetag{3.6.4} we have, e.g., that $\varphi _1 = \varphi _2 = \varphi _3
= \varphi _0$, and then $\varphi _0 \le
\varphi _4 = \pi - \sum _{i=1}^3 \varphi _i = \pi - 3 \varphi _0$, 
from which there follows 
$$
\varphi _0 \in (0, \pi /4] \,.
\tag 3.6.6
$$

We have by \thetag{3.6.3}
$$
\cases
V(P)/p^2 = \left( 3 \sin \varphi _0 + \sin ( \pi - 3 \varphi _0) \right)
/ \left( 8 \sin ^2 (\varphi_0 /2) \right) \\
= \left[ (3 \sin \varphi _0) / \left( 8 \sin ^2 (\varphi_0 /2) \right) \right]
+ \\
\left[ \sin \varphi _0 \cdot (3 - 4 \sin ^2 \varphi _0) 
/ \left( 8 \sin ^2 (\varphi_0 /2) \right) \right] \,.
\endcases
\tag 3.6.7
$$
We investigate both summands in the last expression in \thetag{3.6.7}.
Its first summand is
$$
(3 \sin \varphi _0) / \left( 8 \sin ^2 (\varphi_0 /2) \right) =
(3/4) \cot (\varphi _0 /2)\,, 
\tag 3.6.8
$$
which is a positive
strictly decreasing function of $\varphi _0 \in (0, \pi /4] $.
Its second summand is
$$
\cases
\sin \varphi _0 \cdot (3 - 4 \sin ^2 \varphi _0) 
/ \left( 8 \sin ^2 (\varphi_0 /2) \right) \\ 
= (1/8) \left( \sin \varphi _0 / \sin ^2 (\varphi_0 /2) \right)  
(3 - 4 \sin ^2 \varphi _0) \\ 
= (1/4) \cot (\varphi _0/2) (3 - 4 \sin ^2 \varphi _0)
\endcases
\tag 3.6.9
$$
Here both $\cot (\varphi _0/2)$ and $(3 - 4 \sin ^2 \varphi _0)$ are
positive strictly 
decreasing functions of $\varphi _0 \in (0, \pi /4]$, hence the
last expression in \thetag{3.6.9} is a positive strictly decreasing function of 
$\varphi _0 \in (0, \pi /4]$. Since we know the same for the last expression in
\thetag{3.6.8}, by \thetag{3.6.6} $V(P)/p^2$ in
\thetag{3.6.7} attains its minimum exactly for $\varphi _0 = \pi
/4$, and this minimum is $1 +{\sqrt{2}} = 2.4142 \ldots $. This is greater
than $(5/4) \cot ( \pi /10) = 1.7204 \ldots $, and this ends the proof of case
(1) of this Lemma.

{\bf{2.}}
We turn to the proof of case (2). In case (2) 

\newpage

$$
{\text{for }} 1 \le i \le 5 {\text{ we have }}
\varphi _i \in [\varphi _0 , \pi ) {\text{ and }}
\sum _{i=1}^5 \varphi _i = 2 \pi \,,
\tag 3.6.10
$$
and we have to minimize $ \sum _{i=1}^5 \sin \varphi _i / \left( 8 \sin ^2
(\varphi_0 /2) \right) $.

By \thetag{3.6.4} we have, e.g., that $\varphi _1 = \varphi _2 = \varphi _3
= \varphi _4 = \varphi _0$, and then by (2)
$\pi > \varphi _5 = 2 \pi - \sum _{i=1}^4 \varphi _i = 2 \pi - 4 \varphi _0$, 
from which there follows $\varphi _0 > \pi /4 $.
Also we have $2 \pi = \sum _{i=1}^5 \varphi _i \ge 5 \varphi _0$, from which
there follows $\varphi _0 \le 2 \pi /5$. That is, we have
$$
\varphi _0 \in (\pi /4, 2 \pi /5] \,.
\tag 3.6.11 
$$
We have (like in \thetag{3.6.8})
$$
\cases
V(P)/p^2 = \left( 4 \sin \varphi _0 + \sin (2 \pi - 4 \varphi _0) \right) 
/ \left( 8 \sin ^2 (\varphi _0 /2) \right)  = \\
\left( 4 \sin \varphi _0 + 4 \sin \varphi _0 \cdot (\cos \varphi _0 - 2 \cos ^3 
\varphi _0) \right)
/\left( 8 \sin ^2 (\varphi _0 /2) \right) =\\
\cot (\varphi _0 /2) \cdot (1+\cos \varphi _0 - 2 \cos ^3 \varphi _0) = \\
\cot (\varphi _0 /2) \cdot (1 - \cos \varphi _0) (1+2 \cos \varphi _0 + 2 \cos
^2 \varphi _0) = \\
\sin \varphi _0 \cdot (1+2 \cos \varphi _0 + 2 \cos ^2 \varphi _0) \,.
\endcases
\tag 3.6.12
$$
We are going to show that for $\varphi _0 \in ( \pi /4, 2 \pi /5]$ the last
expression in \thetag{3.6.12} is decreasing. Its derivative with respect to
$\varphi _0$ is the cosine polynomial
$$
-2 -3 \cos \varphi _0 +4 \cos ^2 \varphi _0 + 6 \cos ^3 \varphi _0 =
(2 \cos ^2 \varphi _0 - 1)(2+ 3 \cos \varphi _0) \,.
\tag 3.6.13
$$
By \thetag{3.6.11} the first factor of the last expression in \thetag{3.6.13} 
is negative, and the second factor is positive, hence 
the last expression in \thetag{3.6.13} is negative. That is, 
for $\varphi _0 \in ( \pi /4, 2 \pi /5]$ the last
expression in \thetag{3.6.2} is strictly
decreasing. Hence its minimum is attained exactly for
$\varphi _0 = 2 \pi /5$, and its value is $(5/4)\cot ( \pi /10)$, as asserted in
the Lemma. 
This ends the proof of case (2) of this Lemma.

The case of equality can be checked from the respective steps of the proof
(actually 
one has to consider only case (2), and there we had a unique minimum). 
This ends the proof of this lemma. 
$\blacksquare $
\enddemo

%%%%%%%%%%%%%%%%%%%%%%%%%%%%%%%%%%%%%%%%%%%%%%%%%%%%%%%%%%%%%%%%%%%%%%%%%%%
%%%%%%%%%%%%%%%%%%%%%%%%%%%%%%%%%%%%%%%%%%%%%%%%%%%%%%%%%%%%%%%%%%%%%%%%%%%

\proclaim{Lemma 3.7}
Let $ABA_1C$ be a convex quadrangle having a circumcircle, of radius at most
$p {\sqrt{2}}$. Let the side $AB$ 
be a diameter of the circumcircle, and let $|BA_1|,
|A_1C| \ge p$ and $|CA| \ge p{\sqrt{2}}$. Then the area of $ABA_1C$ is at
least $p^2 \cdot 1.4048 \ldots $.
\endproclaim

%%%%%%%%%%%%%%%%%%%%%%%%%%%%%%%%%%%%%%%%%%%%%%%%%%%%%%%%%%%%%%%%%%%%%%%%%%
%%%%%%%%%%%%%%%%%%%%%%%%%%%%%%%%%%%%%%%%%%%%%%%%%%%%%%%%%%%%%%%%%%%%%%%%%%%

\newpage

\demo{Proof} 
Let $R$ denote the circumradius of $ABA_1C$. We have $\angle BA_1C > \angle
BA_1A = \pi /2$, hence by $|BA_1|,|A_1C| \ge p$ we have $|BC| > p{\sqrt{2}}$.
Then, also using the hypothesis of the lemma, 
in the right triangle $\Delta ABC$ both legs have lengths at least $p
{\sqrt{2}}$, hence 
$$
2R = |AB| \ge p{\sqrt{2}} \cdot {\sqrt{2}}= 2p \,.
\tag 3.7.1
$$
This implies, by the hypothesis of the Lemma that
$$
p \le R \le p{\sqrt{2}}\,.
\tag 3.7.2
$$
 
If fixing $A,B$, we move $C$ on the circumcircle, so that $|BC|, |CA| \ge
p{\sqrt{2}}$, then the minimum of $V(\Delta ABC)$ occurs if $\min \{ |BC|,
|CA| \} = p{\sqrt{2}}$, and the value of this minimum is 
$$
(1/2) \cdot p{\sqrt{2}} \cdot \sqrt{4R^2 - 2 p^2} \,.
\tag 3.7.3
$$

On the other hand, the area of $\Delta BA_1C$ can be estimated from below by
$R$ and $p$ from Lemma 3.2 ($\Delta ABC $ in Lemma 3.2 corresponding
to $\Delta BA_1C$ here, $a_0$ and $b_0$ in Lemma 3.2
corresponding to $p$ here, and $R_0$
in Lemma 3.2 corresponding to $R$ here). Then the minimal area $V(\Delta
BA_1C)$ occurs for $|BA_1| = |A_1C| = p$. Then $\Delta BA_1C$ is an
isosceles triangle, with $\angle BA_1C = \pi - \varphi $, where $\varphi $ is
the (smaller)
central angle in the circle of radius $R$ corresponding to a chord of
length $p$. 
That is, this minimum is
$$
(1/2) p^2 \sin ( \pi - \varphi ) \,,
\tag 3.7.4
$$
where, also using \thetag{3.7.2}, we have
$$
\sin (\varphi /2) = p / (2R) \in [1/(2{\sqrt{2}}) , 1/2] \,.
\tag 3.7.5
$$ 
This implies that for $R \in [p, p{\sqrt{2}}]$ we have that $\varphi \in
(0,\pi /3]$ is a decreasing function of $R$, and thus also $\sin \varphi
$ is a decreasing function of $R$.

Then 
$$
V(ABA_1C) = V(\Delta ABA_1) + V(\Delta A_1CA) \ge \sqrt{2R^2 - p^2} + (1/2)
p^2 \sin \varphi \,.
\tag 3.7.6
$$
By homogeneity, we may assume for simplicity $p=1$.
Then by \thetag{3.7.4} it remains to prove
$$
\sqrt{2R^2 - 1} + (1/2) \sin \varphi \ge 1.4048 \ldots \,.
\tag 3.7.7
$$

\newpage

Here the first summand is an increasing,
and the second summand is a decreasing function of $R \in [p, p{\sqrt{2}}] = 
[1, {\sqrt{2}}]$.

We cover the interval $[1, {\sqrt{2}}]$ by the intervals $[1 + (i-1) \cdot 0.1,
1+i \cdot 0.1]$, where $1 \le i \le 5$. We show the inequality \thetag{3.7.7}
for $R$ in the larger interval $[1, 1.5]$.

In the $i$'th interval we estimate from below ${\sqrt{2R^2 - 1}}$ by its value
in the left hand endpoint of this interval, and $ (1/2) \sin \varphi $ by its
value in the right hand endpoint of this interval. Thus for their sum
we obtain five lower
estimates in the five subintervals (sums of the lower estimates of the two
summands), 
whose values are $1.4048 \ldots $,
$1.5704\ldots $, $1.7261 \ldots $, $1.8763 \ldots $, $2.0230 \ldots $. 
Then the smallest
of these values, i.e., $1.4048 \ldots $ is a lower bound of the left hand
side of \thetag{3.7.7} in the whole interval
$[1, 1.5]$, thus also in the smaller interval $[1, \sqrt{2}]$. This proves the
lemma.
$\blacksquare $
\enddemo

%%%%%%%%%%%%%%%%%%%%%%%%%%%%%%%%%%%%%%%%%%%%%%%%%%%%%%%%%%%%%%%%%%%%%%%%%%
%%%%%%%%%%%%%%%%%%%%%%%%%%%%%%%%%%%%%%%%%%%%%%%%%%%%%%%%%%%%%%%%%%%%%%%%%

\proclaim{Lemma 3.8}
Let a convex pentagon $P$ 
have a circumcircle and let the circumcentre of $P$ be not an
interior point of $P$. Let each side of $P$ have a length
at least $p\,\,(>0)$, and let the circumradius $R$ of $P$ be at most
$p{\sqrt{2}}$. Then the area of 
$P$ is at least $p^2 \cdot 
2.3977 \ldots $. Equality holds if and only if $R =
p{\sqrt{2}}$, and four sides of $P$ 
have length $p$ (witn the line spanned by the
fifth side separating, not strictly, $P$ and its circumcentre). 
\endproclaim

%%%%%%%%%%%%%%%%%%%%%%%%%%%%%%%%%%%%%%%%%%%%%%%%%%%%%%%%%%%%%%%%%%%%%%%%%%
%%%%%%%%%%%%%%%%%%%%%%%%%%%%%%%%%%%%%%%%%%%%%%%%%%%%%%%%%%%%%%%%%%%%%%%%%%%

\demo{Proof} 
We may suppose that $P=A_1A_2A_3A_4A_5$ 
is an extremal pentagon (which exists, since $V(P)/p^2$ 
is similarity invariant).

Let us denote the circumradius of $P$ by $R$, and the angles subtended by the
sides $A_iA_{i+1}$
of $P$ at its circumcentre $O$ by $\varphi _i$, where $1 \le i \le 5$.
We may suppose that the line spanned by the side $A_5A_1$ separates, not
strictly, $P$ and $O$. This implies 
$$
\sum _{i=1} ^4 \varphi _i \le \pi \,.
\tag 3.8.1
$$

Since by hypothesis $R \le p{\sqrt{2}}$, we write
$$
\cases
\varphi _0 := 2 \arcsin \left( p/(2R) \right) \ge 2 \arcsin \left(
1/(2{\sqrt{2}}) \right) \\
= 41.4096 \ldots ^{\circ }, {\text{ thus }} R=p/ \left(
2 \sin (\varphi _0 / 2) \right) \,.
\endcases
\tag 3.8.2
$$
Since each side has a length at least $p$, for $1 \le
i \le 4$ we have
$$
\varphi _1, \varphi _2, \varphi _3, \varphi _4 \ge \varphi _0 {\text{ for }}
1 \le i \le 4\,,
\tag 3.8.3
$$ 
which implies by \thetag{3.8.1} 
$$
\varphi _0 \le \pi /4 \,.
\tag 3.8.4
$$

\newpage

Thus, also using \thetag{3.8.1},
$$ 
{\text{for }} 1 \le i \le 4 {\text{ we have }} \varphi _i \in [\varphi _0
, \pi ) {\text{ and }} \varphi _5 = 2 \pi - \sum _{1 \le i \le 4} \varphi _i 
\in [ \pi , 2 \pi )
\tag 3.8.5
$$
We have 
$$
V(P)/p^2 = (R^2/2)\left( \sum _{i=1} ^4 \sin \varphi _i - \sin (\sum _{i=1} ^4
\varphi _i) \right) / \left( 2R \sin (\varphi _0/2) \right) ^2 \,.
\tag 3.8.6
$$

Our pentagon $P$ is contained in a half-circle of its circumcircle, and on the
corresponding half boundary of the circumcircle the Euclidean and angle
distances depend strictly monotonically on each other. Therefore,
supposing $\varphi _1 >
\varphi _0$, we may move $A_1$
toward $A_2$ in the boundary of the circumcircle, and then $V(P)=
V(A_2A_3A_4A_5) + V(\Delta A_1A_2A_5)$ strictly decreases (since $|A_1A_5|$
and $\angle A_1A_5A_2$ decrease). Observe that by permuting the angles
$\varphi _i$ (i.e., the sides $A_iA_{i+1}$ of $P$)
for $1 \le i \le 4$ the area $V(P)$ does not change. This implies that
for our extremal $P$ we have 
$$
\varphi _i = \varphi _0  {\text{ for }} 1 \le i \le 4,
\tag 3.8.7
$$
and
$$
V(P)/p^2= \left( 4 \sin \varphi _0 - \sin (4 \varphi _0) \right) / \left(
8 \sin ^2 (\varphi _0 /2) \right) \,.
\tag 3.8.8
$$
From \thetag{3.6.12} we have 
$$
V(P)/p^2 = \sin \varphi _0 \cdot (1+ 2 \cos \varphi _0 + 2 \cos ^2 \varphi _0)
\,.
\tag 3.8.9
$$
The derivative of \thetag{3.8.9} with respect to $\varphi _0$ is
from \thetag{3.6.13} 
$$
(2 \cos ^2 \varphi _0 - 1)(2 + 3 \cos \varphi _0) \,.
\tag 3.8.10
$$ 
By \thetag{3.8.2} and \thetag{3.8.4} we have
$$
\varphi _0 \in [2 \arcsin \left( 1/(2{\sqrt{2}}) \right) , \pi /4]\,,
\tag 3.8.11
$$
hence \thetag{3.8.10} is positive for $\varphi _0 \ne \pi /4$. 
Thus \thetag{3.8.9} is a strictly increasing function
of $\varphi _0$ on the interval in \thetag{3.8.11}, 
i.e., is a strictly decreasing function of $R$ (cf. \thetag{3.8.2}),
hence its minimum is attained for the maximal value $p{\sqrt{2}}$ of $R$.

\newpage

The uniqueness of the extremal pentagon $P$ follows from the proof.
$\blacksquare $
\enddemo

%%%%%%%%%%%%%%%%%%%%%%%%%%%%%%%%%%%%%%%%%%%%%%%%%%%%%%%%%%%%%%%%%%%%%%%%%%
%%%%%%%%%%%%%%%%%%%%%%%%%%%%%%%%%%%%%%%%%%%%%%%%%%%%%%%%%%%%%%%%%%%%%%%%%%

\proclaim{Lemma 3.9}
Let a packing consisting of the triangles $\Delta ABC$, $\Delta BA_1C$, 
$\Delta CB_2A$ satisfy the side hypothesis and the angular hypothesis, and 
let \,$c_1=b_1=a_2=c_2=p$.
Let the circumradius of $\Delta ABC$ be some fixed $R \in [p, p{\sqrt{2}}]$.
Let the (smaller) central angles corresponding to chords of the
circumcircle of 
$\Delta ABC$ of lengths
$a,b$ and $p$ be $\varphi _A$, $\varphi _B$ and $\varphi _0 = 
2 \arcsin \left( p/(2R) \right) $.
Let $\varphi _A, \varphi _B \in [\varphi _{{\text{min}}}, \varphi
_{{\text{max}}}] := [2 \arcsin \left( (p{\sqrt{2}})/(2R) \right), 
4 \arcsin \left( p/(2R) \right) ]$.
Then $V(\Delta ABC) + V(\Delta BA_1C) +V(\Delta $
\newline
$CB_2A)$, as a function of the
variables $\varphi _A, \varphi _B$, is a strictly
concave function both of $\varphi _A$
and of $\varphi _B$. In particular, its minimum is attained at some vertex of
the square $[\varphi _{{\text{min}}}, \varphi _{{\text{max}}}] \times 
[\varphi _{{\text{min}}}, \varphi _{{\text{max}}}]$, and only there.
\endproclaim

%%%%%%%%%%%%%%%%%%%%%%%%%%%%%%%%%%%%%%%%%%%%%%%%%%%%%%%%%%%%%%%%%%%%%%%%%%
%%%%%%%%%%%%%%%%%%%%%%%%%%%%%%%%%%%%%%%%%%%%%%%%%%%%%%%%%%%%%%%%%%%%%%%%%%%

\demo{Proof} 
We have
$$
V(\Delta ABC) 
= (R^2/2)\left( \sin \varphi _A + \sin \varphi _B - \sin (\varphi _A
+ \varphi _B) \right) \,,
\tag 3.9.1
$$
$$
V(\Delta BA_1C) 
= R \sin ( \varphi _A /2) {\sqrt{p^2 - R^2 \sin ^2 (\varphi _A /2)}}
\tag 3.9.2
$$
and
$$
V(\Delta CB_2A) 
= R \sin ( \varphi _B /2) {\sqrt{p^2 - R^2 \sin ^2 (\varphi _B /2)}}\,,
\tag 3.9.3
$$
where the expressions under the square root signs are positive.
Adding these inequalities, and using the identity $\sin ^2 (t/2) = (1 - \cos
t)/2$, and the notation
$$
c:=(2p^2-R^2)/(2R^2) \,\,(\in [0,1/2])
\tag 3.9.4
$$
we get 
$$
\cases
\left( V(\Delta ABC) + V(\Delta BA_1C) + V(\Delta CB_2A) \right) /(R^2/2) =
\sin \varphi _A  + \sin \varphi _B \\
- \sin (\varphi _A + \varphi _B) + 
2 \sin (\varphi _A /2) {\sqrt{c + (\cos \varphi _A)/2}} + \\
2 \sin (\varphi _B /2) {\sqrt{c + (\cos \varphi _B)/2}}\,.
\endcases
\tag 3.9.5
$$
The partial derivative of the right hand side of
\thetag{3.9.5} with respect to $\varphi _A$ is 
$$
\cases
\cos \varphi _A - \cos (\varphi _A + \varphi _B) + 
\cos (\varphi _A /2) {\sqrt{c + (\cos \varphi _A)/2}} - \\
\sin (\varphi _A /2) (\sin \varphi _A) \big/ 
{\sqrt{c + (\cos \varphi _A)/2}} 
= \cos \varphi _A - \cos (\varphi _A + \varphi _B) + \\
\left( \cos (\varphi _A /2) \cdot c + \cos (\varphi _A/2) \cdot
(\cos \varphi _A)/2 - \sin (\varphi _A/2) \cdot (\sin \varphi _A)/2
\right) \\
\big/ {\sqrt{c + (\cos \varphi _A)/2}} = 
\left( \cos \varphi _A - \cos (\varphi _A + \varphi _B) \right) + \\
\bigl( \cos (\varphi _A /2) \cdot c + 
\left( \cos (3\varphi _A/2) \right) /2 \bigr) 
\big/ {\sqrt{c + (\cos \varphi _A)/2}}\,,
\endcases
\tag 3.9.6
$$

\newpage

where the expressions under the square root signs are positive.
(The case of $\varphi _B$ is analogous.)
We have to show that \thetag{3.9.6} is a strictly
decreasing function of $\varphi _A$,
i.e., that its partial
derivative with respect to $\varphi _A$ is negative. The last expression
of \thetag{3.9.6} is a
sum of two summands, the second one being a quotient.
Here $\varphi _A$ varies in the interval 
$$
[\varphi _{{\text{min}}}, \varphi
_{\text{max}}] = [2 \arcsin \left( (p{\sqrt{2}})/(2R) \right), 
4 \arcsin \left( p/(2R) \right) ] \,\,(\subset [\pi /3, 2 \pi /3]) \,.
\tag 3.9.7
$$

The derivative of the first summand of the last expression of
\thetag{3.9.6} with respect to $\varphi
_A$ is 
$$
\sin (\varphi _A + \varphi _B) - \sin \varphi _A \,.
\tag 3.9.8
$$
Observe that \thetag{3.9.7} implies 
$$
\varphi _A , \varphi _B \in [\pi /3, 2 \pi /3] {\text{ and so }} 
\varphi _A + \varphi _B \in [2 \pi /3, 4 \pi /3] \,,
\tag 3.9.9
$$
which in turn implies 
$$
\sin (\varphi _A + \varphi _B) \le {\sqrt{3}}/2 \le \sin \varphi _A \,,
\tag 3.9.10
$$
which in turn implies that \thetag{3.9.8} is nonpositive, and so
$$
{\text{the first summand of the last expression of
\thetag{3.9.6} is non-increasing.}}
\tag 3.9.11
$$

We turn to the derivative of the second 
summand of the last expression of
\thetag{3.9.6} with respect to $\varphi _A$. This summand is a
quotient of the form 
$$
f(\varphi _A)/ {\sqrt{g(\varphi _A)}} \,,
\tag 3.9.12
$$ 
with $g(\varphi _A) > 0$. Its derivative with respect to $\varphi _A$ is 
$\left( f'(\varphi _A) g(\varphi _A) - (1/2) f(\varphi _A) \right.$
\newline
$\left. g'(\varphi _A)
\right) /g(\varphi _A)^{3/2}$. So we have to show  
$$
f'(\varphi _a) g(\varphi _A) - (1/2) f(\varphi _A) g'(\varphi _A) \le 0 \,.
\tag 3.9.13
$$
Writing out \thetag{3.9.13} in detail, we have to show
$$
\cases
\left(- (1/2) \sin (\varphi _A /2) \cdot c - (3/4) \sin (3 \varphi _A /2)
\right) \cdot \left( c + (\cos \varphi _A)/2 \right) + \\
(1/2) \cdot \left( \cos (\varphi _A /2) \cdot c + \left( \cos(3\varphi _A/2)
\right) /2 \right) \cdot \left(  (\sin \varphi _A) /2 \right) \le 0\,.
\endcases
\tag 3.9.14
$$
Now we investigate the signs of the factors of the two summands
of \thetag{3.9.14}. 
By \thetag{3.9.9}

\newpage

$$
\varphi _A \in [\pi /3, 2 \pi /3] \,\,{\text{ and so }}\,\, 
3\varphi _A/2 \in [\pi /2 , \pi ] \,,
\tag 3.9.15 
$$
which implies with \thetag{3.9.4} (and \thetag{3.9.7}) that 
{\it{the first factor of the first summand is nonpositive (is actually
negative for $\varphi _A \in (\varphi _{{\text{min}}} , 
\varphi _{{\text{max}}})$)
and the third factor of the second summand 
of \thetag{3.9.14} is positive}}. 

Next we investigate the second factor of the first summand of \thetag{3.9.14}.
By hypothesis of this lemma $\varphi _A \le 4 \arcsin \left( p/(2R) \right) 
\,\,(\le 4 \cdot \pi /6)$.
Then, since the $\arcsin $ function is monotonically increasing on $[0,
{\sqrt{3}}/2] \,\,(\subset [0,1])$, we have
$$
\cases
\cos \varphi _A = 1 - 2 \sin ^2 (\varphi _A /2) \ge 1 - 2 \sin ^2 \left(
2 \arcsin \left( p/(2R) \right) \right) = \\
1 - 2 \left( 2 \cdot \left( p/(2R) \right)
{\sqrt{1- \left( p/(2R) \right) ^2}} \right) ^2 = 
1- 2 ( p/R ) ^2 +(1/2) ( p/R ) ^4 \,.
\endcases
\tag 3.9.16
$$
Then, recalling the notation \thetag{3.9.4}, 
{\it{the second factor of the first summand of \thetag{3.9.14} is
positive}}, since we have, writing
$$
d:=(p/R)^2\,\,( \in [1/2, 1]) \,\,{\text{ (cf. the hypothesis of the lemma),}}
\tag 3.9.17
$$
and using \thetag{3.9.16}, that
$$
c + (\cos \varphi _A)/2 \ge d - 1/2 + (1 - 2d +d^2/2)/2 = d^2/4 \ge 1/16
> 0 \,.
\tag 3.9.18
$$

Next we investigate 
the second factor of the second summand of \thetag{3.9.14}.
By hypothesis of this lemma and by \thetag{3.9.7} we have
$$
\pi /6 \le 
\arcsin \left( (p{\sqrt{2}})/(2R) \right)
= \varphi _{{\text{min}}} /2
\le \varphi _A/2 \le \varphi _{{\text{max}}} /2 \le \pi /3 \,.
\tag 3.9.19
$$ 
Then, since the cosine function is
monotonically decreasing on $[\pi /6 , \pi /3]\,\,\,\,( \subset [0, $
\newline
$\pi /2])$,
we have 
$$
0 < 1/2 \le \cos (\varphi _A /2) \le {\sqrt{1 - p^2/(2R^2)}}\,.
\tag 3.9.20
$$
Now using the identity $\cos (3t) = 4 \cos ^3 t - 3 \cos t$ we have, 
from \thetag{3.9.17} and \thetag{3.9.20},
$$
\cases
\cos (\varphi _A /2) \cdot c + \left( \cos (3 \varphi _A /2) \right) /2 = \\
\cos (\varphi _A /2) \cdot \left( c + 4 \cos ^2 (\varphi _A /2) -3 \right) /2 \\
\le \cos (\varphi _A /2) \cdot \left( 
c + \left( 4 \left( 1 - p^2/(2R^2) \right) - 3 \right) /2 
\right) = \\ 
\cos (\varphi _A /2) \cdot \left( 
d - 1/2 + \left( 4 (1 - d/2) -3 \right) /2 
\right) = 0 \,.
\endcases
\tag 3.9.21
$$

\newpage

Then, recalling the notation \thetag{3.9.4}, {\it{the second factor of the
second summand of \thetag{3.9.14} is nonpositive}}.

Summing up: for $\varphi \in (\varphi _{{\text{min}}}, \varphi
_{{\text{max}}})$
in \thetag{3.9.14} the first summand is the product of a negative
and a positive factor, and the second summand is $1/2$ times the product of
a nonpositive and a positive factor, hence \thetag{3.9.14} is negative.
This shows that 
$$
\cases
{\text{the second summand of the last expression}} \\
{\text{of \thetag{3.9.6} is strictly decreasing.}}
\endcases
\tag 3.9.22
$$

Together with \thetag{3.9.11} this implies that \thetag{3.9.6} is strictly
decreasing,
i.e., \thetag{3.9.5} is a strictly
concave function of $\varphi _A$. Analogously, 
\thetag{3.9.5} is a strictly
concave function of $\varphi _B$ as well. This shows the
statement of the lemma. (The statement about attaining of the minimum is 
immediate from these strict concavity properties.)
$\blacksquare $
\enddemo

%%%%%%%%%%%%%%%%%%%%%%%%%%%%%%%%%%%%%%%%%%%%%%%%%%%%%%%%%%%%%%%%%%%%%%%%%%
%%%%%%%%%%%%%%%%%%%%%%%%%%%%%%%%%%%%%%%%%%%%%%%%%%%%%%%%%%%%%%%%%%%%%%%%%%

\proclaim{Lemma 3.10}
Let the hypotheses of Lemma 3.9 hold. Additionally, let $\varphi _A = \varphi
_{{\text{min}}}$ and $\varphi _B = \varphi _{{\text{max}}}$. If now we allow
$R$ to vary in the interval $[p, p{\sqrt{2}}]$, then $V(\Delta ABC) + V(\Delta
BA_1C) + V(\Delta CB_2A)$, as a function of $R$, is always at least $
p^2 \cdot 1.8720 \ldots $.
\endproclaim

%%%%%%%%%%%%%%%%%%%%%%%%%%%%%%%%%%%%%%%%%%%%%%%%%%%%%%%%%%%%%%%%%%%%%%%%%%
%%%%%%%%%%%%%%%%%%%%%%%%%%%%%%%%%%%%%%%%%%%%%%%%%%%%%%%%%%%%%%%%%%%%%%%%%%

\demo{Proof}
By \thetag{3.9.2} and $\varphi _A = \varphi _{{\text{min}}} := 2 \arcsin \left(
(p{\sqrt{2}})/(2R) \right) $ we have
$$
\cases
V(\Delta BA_1C) = R \sin (\varphi _A/2) {\sqrt{p^2 - R^2 \sin ^2 (\varphi
_A/2)}} \\
= R \cdot \left( (p{\sqrt{2}})/(2R) \right) {\sqrt{p^2 - R^2 \left(
p^2/(2R^2) \right) }} \\
= (p/{\sqrt{2}}) {\sqrt{ p^2 - p^2/2}} = p^2/2 \,.
\endcases
\tag 3.10.1
$$
On the other hand, by $\varphi _B = \varphi _{{\text{max}}} :=
4 \arcsin \left( p/(2R) \right) $ we have that $B_2$ lies on the circumcircle
of $\Delta ABC$. In fact, let us consider the convex
deltoid $OCB_2A$, where $O$ is the
circumcentre of $\Delta ABC$. Its diagonal $OB_2$, which is also its axis of 
symmetry, is cut by its other diagonal into two parts. The part having as one
endpoint $O$ has length 
$$
R \cos (\varphi _{{\text{max}}} /2) = 
R  \cos \left( 2 \arcsin \left( p/(2R) \right) \right)
= R - p^2/(2R) \,,
\tag 3.10.2
$$
and the part having as one endpoint $B_2$ has length 
$$
{\sqrt{p^2 - R^2 \sin ^2 (\varphi _{{\text{max}}} / 2) }} = 
{\sqrt{p^2 - R^2 \sin ^2 \left( 2 \arcsin \left( p/(2R) \right) \right) }} = 
p^2/(2R) \,,
\tag 3.10.3
$$ 
hence
$$
|OB_2|=R \,,
\tag 3.10.4 
$$

\newpage

as asserted.

Writing $O$ for the circumcentre of $\Delta ABC$, we have, also using that
$B_2$ lies on the circumcircle of $\Delta ABC$,
$$
\cases
\left( V(\Delta ABC) + V(\Delta BA_1C) + V(\Delta CB_2A) \right) /p^2 = 
\left( V (\Delta OBC) + \right. \\
\left. V(\Delta OCB_2) + V(\Delta OB_2A) + 
V(\Delta OAB) + V(\Delta BA_1C) \right) / p^2 = \\
\left( (R/p)^2 /2 \right) \left( \sin \varphi _{{\text{min}}} + 2 \sin (\varphi
_{{\text{max}}}/2) - \sin \left( \varphi _{{\text{min}}} + 2 (\varphi
_{{\text{max}}}/2) \right) \right) \\
+ 1/2 =  
\left( (R/p)^2 /2 \right) \left[ \sin \left( 2 \arcsin \left(
(p{\sqrt{2}})/(2R) \right) \right) +
2 \sin \left( 2 \arcsin \right. \right. \\
\left. \left. \left( p/(2R) \right) \right)   
- \sin \left( 2 \arcsin \left( (p{\sqrt{2}})/(2R) \right) + 
2 \cdot 2 \arcsin \left( p/(2R) \right) \right) \right] \\
+ 1/2  \,.
\endcases
\tag 3.10.5
$$
Now recall $R \in [p, p{\sqrt{2}}]$. In \thetag{3.10.5}, last expresssion,
the first factor $\left( (R/p)^2/2 \right) $ of the first summand 
is a positive increasing function of $R$ on $[p, p{\sqrt{2}}]$, and
the second factor of the first summand 
is positively proportional to $V(ABCB_2)$, hence is positive for $R \in 
[p, p{\sqrt{2}}]$ as well. We are going to show that 
$$
\cases
{\text{the second factor of the first summand of}} \\
{\text{the last expression in }} \thetag{3.10.5} {\text{ is a positive}} \\ 
{\text{decreasing function of }} R {\text{ on }} [p, p{\sqrt{2}}] \,.
\endcases
\tag 3.10.6
$$

We have 
$(p{\sqrt{2}})/(2R) \in [1/2, 1/{\sqrt{2}}]$, and on $[1/2, 1/{\sqrt{2}}]$
the arcsine function is increasing, has values there in $[\pi /6, \pi /4]$, 
whose double lie in $[ \pi /3, \pi /2]$. Moreover, 
the sine function is increasing
on this last interval. Hence $\sin \left( 2 \arcsin \left(
(p{\sqrt{2}})/(2R) \right) \right) $ is a decreasing function of $R 
\in  [p, p{\sqrt{2}}]$.

Similarly, we have $p/(2R) \in  [1/{\sqrt{8}}, 1/2]$,
and on $[1/{\sqrt{8}}, 1/2]$ the arcsine function is increasing, 
has values there in $[\arcsin (1/{\sqrt{8}}), \pi /6]= [20.7048 \ldots ^{\circ
}, \pi
/6]$, whose double lie in $[2 \arcsin (1/{\sqrt{8}}), \pi /3]$. Moreover, 
the sine
function is increasing on this last interval. Hence $\sin \left(
2 \arcsin \left( p/(2R) \right) \right) $ is a decreasing function of $R 
\in  [p, p{\sqrt{2}}]$.

These imply that also $2 \arcsin \left( (p{\sqrt{2}})/(2R) \right) +
4 \arcsin \left(p/(2R) \right) $ is a decreasing function of $R 
\in  [p, p{\sqrt{2}}]$. Its values for $R \in [p, p{\sqrt{2}}]$ form the
interval $[\pi /3 + 4 \arcsin (1/{\sqrt{8}}), \pi /2 + 2 \pi /3] =
[142.8192 \ldots ^{\circ } , 210 ^{\circ }]$, on which
the minus sine function is increasing. Hence $- \sin \left(
2 \arcsin \left( (p{\sqrt{2}})/(2R) + 4 \arcsin \left(p/(2R) \right) \right) 
\right) $ is a decreasing function of $R \in  [p, p{\sqrt{2}}]$ as well.

Recapitulating: in the last expression of \thetag{3.10.5} 
the second factor of the first summand, being positive, and being equal to the
sum of three decreasing functions of $R$ on $[p, p{\sqrt{2}}]$, 
is itself a positive decreasing function of $R$ on $[p, p{\sqrt{2}}]$. That
is, \thetag{3.10.6} is proved.

\newpage

Now we turn to the lower estimate of the last expression of \thetag{3.10.5}.
Its second summand $1/2$ is constant. Its first summand is a product of two
positive functions, the first factor being an increasing, the second factor
being a decreasing function of $R \in [p,p{\sqrt{2}}]$, by \thetag{3.10.6}. 

Now we subdivide the
interval $[p,p{\sqrt{2}}]$ into five subintervals: $[p, p \cdot 1.1]$, 
$[p \cdot 1.1, p \cdot 1.2]$, $[p \cdot 1.2, p \cdot 1.3]$, $[p \cdot 1.3, 
p \cdot 1.4]$, $[p \cdot 1.4, p {\sqrt{2}}]$. On each of these
subintervals we estimate from
below the above mentioned first factor by its value at the
left hand endpoint of the subinterval, and 
we estimate from below the above mentioned second factor by its value at the
right hand endpoint of the subinterval. This way we obtain the following
numerical values as lower estimates for the first summand of \thetag{3.10.5} in
the above five subintervals: $1.3730 \ldots $, $1.3949 \ldots $, 
$1.3912 \ldots $, $1.3720 \ldots $ and $1.5528 \ldots $ (observe that the last
interval is much shorter than the previous ones, therefore do
we have there a much
higher value). The minimum 
of these five numbers, i.e., $1.3720 \ldots $ is a
lower estimate of the first summand of the last expression in \thetag{3.10.5}
on the whole interval $[p,p{\sqrt{2}}]$. The second summand 
of the last expression in \thetag{3.10.5} being the constant $1/2$, we
obtain the statement of the lemma.
$\blacksquare $
\enddemo

%%%%%%%%%%%%%%%%%%%%%%%%%%%%%%%%%%%%%%%%%%%%%%%%%%%%%%%%%%%%%%%%%%%%%%%%%%
%%%%%%%%%%%%%%%%%%%%%%%%%%%%%%%%%%%%%%%%%%%%%%%%%%%%%%%%%%%%%%%%%%%%%%%%%%

\head{\S 4 The main lemmas}\endhead

In this paragraph we will prove statements in ${\Bbb R}^2$.

We say that a packing of triangles, any two of which join by entire edges, 
or by common vertices, or being disjoint, satisfies the
{\it{angular hypothesis}}, if for any two triangles joining with a common
edge the sum of the angles opposite to the common side is at most $\pi $. 
An example is the Delone triangulation corresponding to an $(r,R)$-system
in ${\Bbb R}^2$, where $0 < r < R < \infty $ (with Delone polygons
with more than three sides triangulated in an arbitrary way). We say that this
packing satisfies the {\it{side hypothesis}} if all the sides of these
triangles are at least $p$, where $p \in (0, \infty )$ is fixed. In our
lemmas there will occur also 
right and obtuse triangles. Then of course the side
hypothesis is sufficient to be supposed for these triangles only 
for their sides adjacent to the right or obtuse angle.

In the proofs of our Lemmas we will use the notations $a_i,b_i,c_i$ and
$\alpha _i , \beta _i , \gamma _i$ for $1 \le i \le 3$ as introduced before
Lemma 3.3. 

%%%%%%%%%%%%%%%%%%%%%%%%%%%%%%%%%%%%%%%%%%%%%%%%%%%%%%%%%%%%%%%%%%%%%%%%%%
%%%%%%%%%%%%%%%%%%%%%%%%%%%%%%%%%%%%%%%%%%%%%%%%%%%%%%%%%%%%%%%%%%%%%%%%%%%

\proclaim{Lemma 4.1}
Let the packing consisting of the triangles $\Delta ABC,
\Delta BA_1C$ satisfy the angular and side
hypotheses and let 
$\alpha , \beta , \gamma \le \pi /2 $, $\alpha _1 \ge \pi /2$.
Then the sum of the areas of the two triangles is at least $p^2$.
Equality holds if and only if $b=c=c_1=b_1=p$ and $\alpha _1 = \pi /2$.
\endproclaim

%%%%%%%%%%%%%%%%%%%%%%%%%%%%%%%%%%%%%%%%%%%%%%%%%%%%%%%%%%%%%%%%%%%%%%%%%%
%%%%%%%%%%%%%%%%%%%%%%%%%%%%%%%%%%%%%%%%%%%%%%%%%%%%%%%%%%%%%%%%%%%%%%%%%%%

\demo{Proof} 
Consider an extremal configuration (this exists). 

We make a case distinction. 

\newpage

(1)
$\alpha + \alpha _1 < \pi $ and $\beta ,\gamma < \pi /2$.
\newline
(2)
$\alpha + \alpha _1 < \pi, {\text{ and, e.g., }} \beta = \pi /2 $ (the case
when here $\gamma = \pi /2$ is analogous).
\newline
(3)
$\alpha + \alpha _1 = \pi $, i.e., the quadrangle $ABA_1C$ has a circumcircle.

{\bf{1.}}
In case (1) we apply Lemma 3.1 to $\Delta BCA_1$ and $\Delta BCA$
(rather than $\Delta ABC$ 
in Lemma 3.1), respectively. Then by $\beta _1, \gamma _1
<\pi /2$ and $\alpha _1 < \pi
- \alpha $ we obtain $b_1=c_1=p$, and by $\alpha < \pi - \alpha _1$ and 
$\beta ,\gamma < \pi /2$ we obtain
$b=c=p$. Hence $ABC$ and $A_1BC$ are congruent triangles, thus
$\pi /2  \ge  \alpha = \alpha _1 \ge \pi /2$, thus $\alpha = \alpha _1 = \pi
/2$, a contradiction to our hypothesis $\alpha + \alpha _1 < \pi $.
 
{\bf{2.}}
In case (2) fixing $\Delta CA_1B$ 
and the side length $c=|AB|$, decrease $\beta
$ a bit. Then the hypotheses of the lemma
remain valid, and the area of $\Delta ABC$
decreases, a contradiction.

{\bf{3.}}
In case (3), unless $\alpha _1 = \pi /2$, we can apply
Lemma 3.3 (with the same notations there as here), 
obtaining $\beta = \pi /2$ or $\gamma = \pi /2$. 
Thus we have to investrigate the cases
\newline
(a)
$\alpha _1 = \pi /2$, and
\newline
(b)  
$\beta = \pi /2$ (the case $\gamma = \pi /2$ is analogous).

In case (a), observe that for $\alpha _1 = \pi /2$ we have by (3)
also $\alpha  = \pi /2$, and the
total area of the triangles $\Delta ABC$ and $\Delta A_1BC$ is
$$
(bc+b_1c_1)/2 \ge p^2\,,
\tag 4.1.1
$$
proving the Lemma in case (a).

In case (b), let, e.g., $\beta = \pi /2$ (the case of $\gamma = \pi /2$ is
analogous).
By (3) all four
vertices of our two triangles lie on a circle, and by $\beta = \pi /2$ the
side $AC$ is a diameter of this circle. Let $O$ be the centre of this circle.
Then the quadrangle $ABA_1C$ is inscribed to this circle, and its area (i.e.,
the sum of the areas of our two triangles in the lemma) is
$$
(b/2)^2 \cdot (\sin \angle AOB + \sin \angle BOA_1 + \sin \angle A_1OC)/2\,.
\tag 4.1.2
$$
The longest side of this quadrangle is $b$, the other three sides are
$c, c_1, b_1 \ge p$. The central angles corresponding to the sides $c,
c_1,b_1$ will be denoted by $\varphi _1, \varphi _2, \varphi _3$, respectively.
Then $\varphi _i \ge \varphi _0 := 2 \arcsin \left( p/(2R) \right) $ for
$1 \le i \le 3$ and $\varphi _1 + \varphi _2 + \varphi _3 = \pi $, which implies
$\varphi _i \in [\varphi _0, \pi )$ for $1 \le i \le 3$. We will 
use concavity of the function $\sin x$ in $[0, \pi ]$. Thus at most one
$\varphi _i$ can be in the open interval $(\varphi _0, \pi )$, since else we
decrease one and increase another one by a bit, preserving their sum and
increasing their difference, and then
$\sum _{1 \le i \le 3} \sin \varphi _i$ decreases (while $b$ is fixed), 
and then \thetag{4.1.2}
decreases also, a contradiction. So two of the $\varphi _i$'s equal $\varphi
_0$, and then the third one equals $\pi - 2 \varphi _0 \,\,(\ge \varphi _0)$.
The last inequality implies 
$6 \arcsin \left( p/(2R) \right) = 3 \varphi _0 \le \pi $
thus $2 \arcsin \left( p/(2R) \right) = \varphi _0 \le \pi /3$, which in
turn implies $b = 2R \ge 2p$.

\newpage

Recapitulating, the area of our quadrangle is minimal, e.g., for a symmetric
trapezoid of longer base of length $b$, which is also the diameter of the
circumscribed circle of our quadrangle, and equal sides of lengths $p$, 
and other base $c_1 \ge p$. Observe that here we have a free parameter $b$,
whose value determines our quadrangle $ABA_1C$ uniquely, and we have to find
the minimal area of our quadrangle when $b$ varies. 

We calculate the area of our symmetric trapezoid as the arithmetic mean of the
two bases times the height corresponding to the bases. The equal sides of our
symmetric trapezoid have lenghts $p$, and enclose with the longer base an
angle $\angle CAB = \angle ACA_1 = \arccos (p/b) \ge \pi /3$. Then the longer
base is $b \ge 2p$, the shorter base is 
$b - 2p \cos \angle CAB \ge b - 2 p \cos (\pi /3) = b-p
\ge p$, and the height corresponding to the bases is 
$p \sin \angle CAB \ge p \sin (\pi /3)
= p \sqrt{3}/2$. Putting all these estimates together, we gain that the area
of our quadrangle is at least $(2p + p)/2 \cdot p \sqrt{3}/2 =   
p^2 \cdot (3\sqrt{3}/4) > p^2$, proving the lemma in this case as well. (This
lower estimate is sharp for a ``half of a regular hexagon or side length $p$''.)

The case of equality follows from the proof. (All cases except the one
investigated in {\thetag{{4.1.1}} were contradictory, or gave better estimates.)
$\blacksquare $
\enddemo

%%%%%%%%%%%%%%%%%%%%%%%%%%%%%%%%%%%%%%%%%%%%%%%%%%%%%%%%%%%%%%%%%%%%%%%%%%
%%%%%%%%%%%%%%%%%%%%%%%%%%%%%%%%%%%%%%%%%%%%%%%%%%%%%%%%%%%%%%%%%%%%%%%%%%%

\proclaim{Lemma 4.2}
Let the packing consisting of the triangles $\Delta ABC,
\Delta BA_1C$, $\Delta CB_2A$ satisfy the angular and side hypotheses,
let the circumradii of $\Delta BA_1C$ and $\Delta CB_2A$ be at most
$p{\sqrt{2}}$, and let 
$\alpha , \beta , \gamma \le \pi /2 $, $\alpha _1 , \beta _2 \ge \pi /2$.
Then the sum of the areas of the three
triangles is at least $p^2(\sqrt{7}/4+1)= p^2 \cdot 1.6614 \ldots $, with
equality if and only if 
$c_1=b_1=a_2=c_2=c=p$ and $\alpha _1 = \beta _2 = \pi /2$. 
If additionally $c \ge p{\sqrt{2}}$, then the sum of the areas of
the three triangles is at least $p^2 \cdot (5/4) \cot (\pi /10)  
= p^2 \cdot 1.7204 \ldots $.
The second lower 
estimate cannot be strengthened by replacing its right hand side
by any number greater than 
$p^2(1+{\sqrt{3}}/2)=p^2 \cdot 1.8660 \ldots $.
\endproclaim

%%%%%%%%%%%%%%%%%%%%%%%%%%%%%%%%%%%%%%%%%%%%%%%%%%%%%%%%%%%%%%%%%%%%%%%%%%
%%%%%%%%%%%%%%%%%%%%%%%%%%%%%%%%%%%%%%%%%%%%%%%%%%%%%%%%%%%%%%%%%%%%%%%%%%%

\demo{Proof} 
We begin with showing an example for the second inequality
having total area 
$p^2(1+{\sqrt{3}}/2)$. We let 
$c_1=b_1=a_2=c_2=p$ and $c = p{\sqrt{2}}$ 
and $\alpha _1 = \beta _2 = \pi /2$. 

Now we turn to prove the two lower estimates.
Consider an extremal configuration (this exists). 

We have two angular hypotheses, namely $\alpha + \alpha _1 \le \pi $
and $\beta + \beta _2 \le \pi $. We distinguish three cases. Either 
\newline
(1)
$\alpha + \alpha _1 = \beta + \beta _2 = \pi $, or, e.g.,
\newline
(2)
$\alpha + \alpha _1 = \pi > \beta + \beta _2 $, or
\newline
(3)
$\alpha + \alpha _1 , \beta + \beta _2 < \pi $.

{\bf{1.}}
We begin with case (1). That is, the pentagon $ABA_1CB_2$ has a circumcircle.
Then by Lemma 3.6 we have $V(\Delta ABC)+V(\Delta A_1BC)+V(\Delta AB_2C)=
V(ABA_1CB_2) \ge p^2 \cdot (5/4) \cot ( \pi /10) = p^2 \cdot 1.7204 \ldots >
p^2({\sqrt{7}}/4 +1) = p^2 \cdot 1.6614 \ldots $. This proves the first
inequality of our Lemma, and also the second inequality (for $c \ge
p{\sqrt{2}}$) in case {\bf{1}}.
 
{\bf{2.}}
We continue with case (2). We apply Lemma 3.3 to the triangles 

\newpage

$\Delta ABC$
and $\Delta BA_1C$. Then $\alpha $ and $\alpha _1$ decrease while $\beta
+ \beta _2 < \pi $, so the
angular hypotheses remain preserved (observe $\beta + \beta _2 < \pi $). 
The side hypothesis is preserved since 
by $\alpha _1 \ge \pi /2$
$$
a=a_1 \ge {\sqrt{b_1^2+c_1^2}} \ge p {\sqrt{2}} > p \,.
\tag 4.2.1
$$
For $\alpha _1 > \pi /2$ the circumradius hypothesis remains preserved, since 
$\Delta CB_2A$ remains preserved, and the circumradius of $\Delta BA_1C$ is
$a/(2 \sin \alpha _1)$, where $a$ decreases and $\alpha _1\,\,( > \pi
/2)$ decreases.
Since $\beta _2$ remains
preserved, this motion can be prevented by
\newline
(a)
$\alpha _1 = \pi /2 $ and then by (2) also $\alpha = \pi /2$, or
\newline
(b) 
$\beta = \pi /2$, or
\newline
(c) $\gamma = \pi /2$.

In case (a) we use the analogue of \thetag{4.2.1} for $b$, i.e., $b \ge
p{\sqrt{2}}$,
which yields 
$$
a = {\sqrt{b^2+c^2}} \ge {\sqrt{(p{\sqrt{2}})^2+p^2}}
=p{\sqrt{3}} \,.
\tag 4.2.2
$$ 
Then 
$$
V(\Delta ABC) = (1/2) \cdot bc \ge (1/2) \cdot
p{\sqrt{2}} \cdot p = p^2{\sqrt{2}}/2
\tag 4.2.3
$$ 
Further, moving $A_1$ on the
circumcircle of $\Delta ABC$, the minimal area $V(\Delta BA_1C)$
occurs when $\min \{ b_1,c_1 \} =
p$, and the minimal area is by \thetag{4.2.2}
$$
V(\Delta BA_1C) = (1/2) \cdot p{\sqrt{a^2 - p^2}}
\ge (1/2) \cdot p{\sqrt{3p^2 - p^2}} = p^2 \cdot {\sqrt{2}}/2
\tag 4.2.4
$$ 
For $\Delta AB_2C$ we use the circumradius hypothesis, which gives by Lemma
3.2, (3)
$$
V(\Delta AB_2C) \ge p^2{\sqrt{7}}/8 \,.
\tag 4.2.5
$$
Adding \thetag{4.2.3}, \thetag{4.2.4} and \thetag{4.2.5} we have
$$
\cases
V(\Delta ABC) + V(\Delta BA_1C) + V(\Delta AB_2C) \ge p^2
({\sqrt{2}}+{\sqrt{7}}/8) \\
= p^2 \cdot 1.7449 \ldots > p^2 \cdot (5/4) \cot (\pi /10) = p^2 \cdot
1.7204 \ldots \\
> p^2({\sqrt{7}}/4 + 1) = p^2 \cdot 1.6614 \ldots \,,
\endcases
\tag 4.2.6
$$
which proves both lower estimates of the lemma in case (a). 

\newpage

In case (b), i.e., when $\beta = \pi /2$, we have by (2) $\beta _2 < \pi /2$,
while by hypothesis of the lemma $\beta _2 \ge \pi /2$, a contradiction.

In case (c), i.e., when $\gamma = \pi /2$, the quadrangle $ABA_1C$ has a
circumcircle, and its circumradius is at most $p{\sqrt{2}}$ (by the
hypothesis about the circumradii in the lemma). By $\gamma = \pi /2$ the side
$AB$ is a diameter of the circumcircle, and by the side hypothesis we have
$|BA_1|,|A_1C| \ge p$ and also $|AC| \ge {\sqrt{a_2^2 + c_2^2}} \ge
p{\sqrt{2}}$. Thus we may apply Lemma 3.7, yielding 
$$
V(ABA_1C) \ge p^2 \cdot 1.4048 \ldots \,.
\tag 4.2.7
$$
Then by the circumradius hypothesis and Lemma 3.2 we have
$$
V(\Delta CB_2A) \ge 0.3307 \ldots \,.
\tag 4.2.8
$$
Adding \thetag{4.2.7} and \thetag{4.2.8} we get
$$
\cases
V(\Delta ABC) + V(\Delta BA_1C) + V(\Delta CB_2A) \ge p^2 \cdot
1.7355 \ldots \\
> p^2 \cdot (5/4) \cot (\pi /10) = p^2 \cdot 1.7204 \ldots \\
> p^2({\sqrt{7}}/4 + 1) = p^2 \cdot 1.6614 \ldots \,,
\endcases
\tag 4.2.9
$$
proving both lower estimates of the lemma in case (c).

{\bf{3.}}
We turn to case (3).
We distinguish the cases
\newline
(a) 
$\alpha , \beta , \gamma < \pi /2$, and
\newline
(b)
$\beta = \pi /2$ (the case $\alpha = \pi /2$ is analogous)
\newline
(c)
$\gamma = \pi /2$.

In case (a)
we apply Lemma
3.3 to the triangles $\Delta ABC$ and $\Delta BA_1C$ (with the same notations
here as in Lemma 3.3). The case of the triangles $\Delta ABC$ and $\Delta CB_2A$
can be settled in an identical way.

By $\alpha _1 \ge \pi /2$ we have $a=a_1 \ge p{\sqrt{2}} >p$, so the side
hypothesis is satisfied if we apply Lemma 3.3 (the other side lengths 
are unchanged). The angular hypothesis $\alpha + \alpha _1 \le \pi $
plays no role here. For $\alpha _1 > \pi /2$ the circumradius hypotheses are
preserved by the same reasoning as in {\bf{2}}.
So by Lemma 3.3 either $\alpha _1 = \pi /2$ or $\beta + \beta _2 
= \pi /2$, but the second case is impossible by case (3), so 
$$
\alpha _1 = \pi/2 \,.
\tag 4.2.10
$$
Analogously, by Lemma 3.3 applied to $\Delta ABC$ and $\Delta CB_2A$,
one obtains 

\newpage

$$
\beta _2 = \pi /2 \,.
\tag 4.2.11
$$

Thus $\alpha _1 = \beta _2 = \pi /2$.
Then $b_1, c_1 \ge p$ gives $a=a_1 \ge p{\sqrt{2}}$, 
and analogously, $b=b_1 \ge p{\sqrt{2}}$. Moreover, we have $c \ge p$.
Then \thetag{3.0} yields 
$$
V(\Delta ABC) \ge p^2 \cdot {\sqrt{7}}/4 \,.
\tag 4.2.12
$$ 
Further,
$$
V(\Delta BA_1C) = b_1c_1 \sin \alpha _1 /2 \ge p^2/2 {\text{ and similarly }}
V(\Delta CB_2A) \ge p^2/2 \,,
\tag 4.2.13
$$ 
Adding \thetag{4.2.12} and the two inequalities in \thetag{4.2.13} we obtain
the first inequality of the lemma in case (a). 

If moreover we
have $c \ge p{\sqrt{2}}$, then \thetag{3.0} yields even 
$$
V(\Delta ABC) \ge p^2 \cdot {\sqrt{3}}/2\,,
\tag 4.2.14
$$ 
and then adding \thetag{4.2.14} and the two inequalities in \thetag{4.2.13}
the total area of the
three triangles is at least $p^2 (1+{\sqrt{3}}/2) = p^2 \cdot
1.8660 \ldots > p^2 \cdot (5/4) \cot ( \pi /10) = p^2 \cdot 1.7204 \ldots $.
This proves the second inequality of the lemma in case (a).

In case (b), i.e., when $\beta = \pi /2$, we have by (3) $\beta _2 < \pi /2$,
while by hypothesis of the lemma $\beta _2 \ge \pi /2$, a contradiction.

In case (c) we have 
$\gamma = \pi /2$. We decrease $\gamma $ a bit,
rotating $\Delta BA_1C$ and $\Delta CB_2A$ towards each other.
Then $V(\Delta ABC)$ decreases, $V(\Delta BA_1C)$ and $V(\Delta CB_2A)$ remain
constant, so their sum decreases. 
The circumradius hypotheses remain preserved.
The side hypothesis is to be checked only
for $c$, but that was originally by $a,b \ge p$ and $\gamma = \pi /2$ at
least $p {\sqrt{2}}$, so the side hypothesis remains valid. Similarly,
originally we had $\alpha , \beta < \pi /2$, so they remain acute angles. So
all the hypotheses of this lemma (and also those of
case (3)) remain preserved, and the total area of our
three triangles decreased. This is a contradiction.

The case of equality in the first inequality of the lemma follows from the
proof. (In {\bf{1}} and {\bf{2}} we had strict inequalities or contradiction, 
{\bf3}}, cases
(b) and (c) were contradictory, and {\bf{3}}, case (a) is easily discussed.)
$\blacksquare $
\enddemo

%%%%%%%%%%%%%%%%%%%%%%%%%%%%%%%%%%%%%%%%%%%%%%%%%%%%%%%%%%%%%%%%%%%%%%%%%%
%%%%%%%%%%%%%%%%%%%%%%%%%%%%%%%%%%%%%%%%%%%%%%%%%%%%%%%%%%%%%%%%%%%%%%%%%%%

\proclaim{Lemma 4.3}
Let the packing consisting of the triangles $\Delta ABC$,
$\Delta BA_1C$, $\Delta CB_2A$ and $\Delta AC_3B$ 
satisfy the angular and side hypotheses, let the circumradii of 
$\Delta BA_1C$, $\Delta CB_2A$ and $\Delta
AC_3B$ be at most $p \cdot {\sqrt{2}}$, and let
$\alpha , \beta , \gamma \le \pi /2 $ and $\alpha _1 , \beta _2 , \gamma _3 
\ge \pi /2$.
Then the sum of the
areas of the four triangles is at least $p^2 \cdot \left[ (5/4) \cot ( \pi /10
) + {\sqrt{7}}/ \right. $
\newline
$\left. 8 \right] = p^2 \cdot 2.0511 \ldots $.
This estimate cannot be strengthened by replacing the right hand side
of this inequality by any number greater than  
$p^2(\sqrt{3}/2 + 3/2) = p^2 \cdot 2.366 \ldots $. 
\endproclaim

%%%%%%%%%%%%%%%%%%%%%%%%%%%%%%%%%%%%%%%%%%%%%%%%%%%%%%%%%%%%%%%%%%%%%%%%%%
%%%%%%%%%%%%%%%%%%%%%%%%%%%%%%%%%%%%%%%%%%%%%%%%%%%%%%%%%%%%%%%%%%%%%%%%%%%

\newpage

\demo{Proof} 
We begin with showing an example having a total area $p^2({\sqrt{3}}/2 +
3/2)$. We let $c_1=b_1=a_2=c_2=b_3=a_3=p$ and
$\alpha _1 = \beta _2 = \gamma _3 = \pi /2$.

Now we turn to the lower estimate.
By $a_3,b_3 \ge p$ and $\gamma _3 \ge \pi /2$ we have $c \ge p{\sqrt{2}}$.
Therefore we can apply Lemma 4.2 to the triangles $\Delta ABC$, $\Delta BA_1C$,
$\Delta CB_2A$, with the same notations here as in Lemma 4.2. Therefore the
second inequality of Lemma 4.2 applies, and gives 
$$
\cases
V(\Delta ABC) + V(\Delta BA_1C) + V(\Delta CB_2A) \\
\ge p^2 \cdot (5/4) \cot ( \pi /10) = p^2 \cdot 1.7204 \ldots \,.
\endcases
\tag 4.3.1
$$
Lemma 3.2, (3) implies
$$
V(\Delta AC_3B) \ge p^2 \cdot {\sqrt{7}}/8 = p^2 \cdot 0.3307 \ldots \,.
\tag 4.3.2
$$
Adding these we obtain
$$
\cases
V(\Delta ABC) + V(\Delta BA_1C) + V(\Delta CB_2A) + V(\Delta AC_3B) \\
\ge p^2 \cdot \left[ (5/4) \cot ( \pi /10 ) + {\sqrt{7}}/8 \right] =
p^2 \cdot 2.0511 \ldots \,,
\endcases
\tag 4.3.3
$$
proving the lemma.
$\blacksquare $
\enddemo

%%%%%%%%%%%%%%%%%%%%%%%%%%%%%%%%%%%%%%%%%%%%%%%%%%%%%%%%%%%%%%%%%%%%%%%%%%
%%%%%%%%%%%%%%%%%%%%%%%%%%%%%%%%%%%%%%%%%%%%%%%%%%%%%%%%%%%%%%%%%%%%%%%%%%%

\proclaim{Lemma 4.4}
Let the packing consisting of the triangles $\Delta ABC,
\Delta BA_1C$ satisfy the angular and side hypotheses and let 
$\gamma , \alpha _1 \ge \pi /2$. Further let the
circumradius of $\Delta ABC$ be at most $p\sqrt{2}$ (or $p \sqrt{5/2}$).
Then $V(\Delta ABC)$ 
is at least $p^2 \cdot (\sqrt{7}+\sqrt{3})/8 = p^2 \cdot 0.5472 \ldots $ (or
$p^2/2$). Equality holds if and only if $b=b_1=c_1=p$, $a=p \sqrt{2}$ and
the circumradius of $\Delta ABC$ is $p\sqrt{2}$ (or $p\sqrt{5/2}$).
\endproclaim

%%%%%%%%%%%%%%%%%%%%%%%%%%%%%%%%%%%%%%%%%%%%%%%%%%%%%%%%%%%%%%%%%%%%%%%%%%
%%%%%%%%%%%%%%%%%%%%%%%%%%%%%%%%%%%%%%%%%%%%%%%%%%%%%%%%%%%%%%%%%%%%%%%%%%%

\demo{Proof} 
We have $b_1,c_1 \ge p$ and $\alpha _1 \ge \pi /2$, which imply $a \ge
p \sqrt{2}$. Moreover, we have $b \ge p$.
Now we can apply Lemma 3.2, (1) (or (2)) 
to $\Delta ABC$ (with the same notation of vertices 
as in Lemma 3.2), which yields both statements of the Lemma.
The case of equality follows from the proof of this lemma and from Lemma 3.2.
$\blacksquare $
\enddemo

%%%%%%%%%%%%%%%%%%%%%%%%%%%%%%%%%%%%%%%%%%%%%%%%%%%%%%%%%%%%%%%%%%%%%%%%%%%

\proclaim{Lemma 4.5}
Let the packing consisting of the triangles $\Delta ABC,
\Delta A_1BC$ satisfy the angular and side hypotheses and let 
$\gamma , \alpha _1 \ge \pi /2$. Further let the
circumradius of $\Delta ABC$ be at most $p\sqrt{2}$.
Then the sum of the areas of these two
triangles is at least $p^2 \cdot (\sqrt{7} + \sqrt{3} +4)/8 = p^2 \cdot 1.0472
\ldots $. Equality holds if and only if $c_1=b_1=b=p$, $\alpha _1 = \pi /2$ and
the circumradius of $\Delta ABC$ is $p\sqrt{2}$.
\endproclaim

%%%%%%%%%%%%%%%%%%%%%%%%%%%%%%%%%%%%%%%%%%%%%%%%%%%%%%%%%%%%%%%%%%%%%%%%%%
%%%%%%%%%%%%%%%%%%%%%%%%%%%%%%%%%%%%%%%%%%%%%%%%%%%%%%%%%%%%%%%%%%%%%%%%%%%

\demo{Proof} 
Consider an extremal configuration (this exists).

We apply Lemma 3.2 to $\Delta ABC$ (with the same notation of vertices 
there as here),
with $a_0:=a$, $b_0:=p$ and $R_0:=p{\sqrt{2}}$. Observe that the angular
hypothesis remains preserved by the application of Lemma 3.2, since in the
proof of Lemma 3.2 both acute angles of $\Delta ABC$ there do not increase.
Thus we get that 

\newpage

$$
b=p {\text{ and the circumradius of }} \Delta ABC {\text{ is }} p {\sqrt{2}}\,.
\tag 4.5.1
$$

{\bf{1.}}
Now we make a case distinction.
\newline
(1) 
$\alpha + \alpha _1 = \pi $, i.e., that the quadrangle $ABA_1C$ has a
circumcircle, and
\newline
(2)
$\alpha + \alpha _1 <  \pi $.

First we deal with case (1).
Then we apply Lemma 3.1 to $\Delta BA_1C$ (rather than $\Delta ABC$ in Lemma
3.1), with $\gamma _0 $ there replaced by $\alpha _1$ here.
Thus we have either $b_1=p$ or $c_1=p$. Although these are two different cases,
the areas of the quadrangle $ABA_1C$ coincide in these two cases. Thus for
estimating the area of our quadrangle below, we may suppose that 
$$
b_1=p\,.
\tag 4.5.2
$$

The area of our quadrangle also equals 
$V(\Delta ABA_1) + V(\Delta A_1CA)$. Observe
that by $\gamma \ge \pi /2$ our quadrangle lies in a closed half-circle $H$
of its
circumcircle, partly bounded by a diameter of the circumcircle, parallel to
the line $AB$. The orthogonal bisector line of side $AA_1$ meets the shorter arc
${\widehat{AA_1}}$ 
of the circumcircle in the midpoint $M$ of this arc, which lies in the
relative interior of $H$ (with 
respect to the entire circumcircle). The same orthogonal
bisector meets the longer arc $AA_1$ of the circumcircle in the opposite point
of $M$ of the circumcircle, thus in the relative interior of the complementary
halfcircle to $H$. Therefore moving
$B$ on the circumcircle of our quadrangle, towards $A_1$, $V(\Delta ABA_1)$
strictly decreases, while $V(\Delta A_1CA)$ remains constant, therefore the
area of our quadrangle $ABA_1C$ also strictly decreases.
Therefore in the extremal configuration we have also 
$$
c_1=p\,.
\tag 4.5.3
$$ 

If we had, rather than {\thetag{4.5.2}},
$$
c_1=p\,,
\tag 4.5.4
$$
then by equality of the areas in these two cases (which are in fact
axially symmetric images of each other), 
in the extremal configuration we have also 
$$
b_1=p\,.
\tag 4.5.5
$$

We turn to case (2).
Then we again apply Lemma 3.1 to $\Delta BA_1C$ 
(rather than $\Delta ABC$ in Lemma 3.1), with $\gamma _0 $ there replaced by
$\pi - \alpha $ here.
By $\alpha + \alpha _1 < \pi $ we get that 

\newpage

$$
b_1=c_1=p\,.
\tag 4.5.6
$$

{\it{Therefore \thetag{4.5.6} is valid in both cases}} (1) {\it{and}} (2).

{\bf{2.}}
Now we apply Lemma 3.5 (with the
same notations there as here). All its hypotheses are satisfied by
\thetag{4.5.1} (observe $p = b \le p{\sqrt{7/2}}$) and \thetag{4.5.6}.
By the conclusion of Lemma 3.5 we may suppose at the lower estimate of
$V(\Delta ABC) + V(BA_1C)$ that either
\newline
(1) 
$\alpha _1 = \pi/2$, or
\newline
(2) 
$\alpha + \alpha _1 = \pi $.

In case (1) we have $V(\Delta BA_1C) = b_1c_1/2 = p^2/2$.
Moreover, by Lemma 4.4 
we have $V(\Delta ABC) \ge p^2 \cdot ({\sqrt{7}} +
{\sqrt{3}})/8$, hence their sum is at least $p^2 \cdot ({\sqrt{7}} +
{\sqrt{3}} + 4)/8 = p^2 \cdot 1.0472 \ldots $, as asserted in the Lemma,
finishing the proof of case (1).

In case (2) the quadrangle $ABA_1C$ has a
circumcircle, which has by \thetag{4.5.1} radius $p
{\sqrt{2}}$. By \thetag{4.5.1} and \thetag{4.5.6} we have $c_1=b_1=b=p$. Thus
the quadrangle $ABA_1C$ is completely determined, is a symmetric trapezoid
inscribed to a circle of radius $p{\sqrt{2}}$,
and a straightforward calculation gives that its area is $p^2 \cdot
7\sqrt{7}/16 = p^2 \cdot 1.1575 \ldots $ which is a larger value than asserted
in the lemma, finishing the proof of case (2).

The case of equality follows from the proof.
$\blacksquare $
\enddemo

%%%%%%%%%%%%%%%%%%%%%%%%%%%%%%%%%%%%%%%%%%%%%%%%%%%%%%%%%%%%%%%%%%%%%%%%%%
%%%%%%%%%%%%%%%%%%%%%%%%%%%%%%%%%%%%%%%%%%%%%%%%%%%%%%%%%%%%%%%%%%%%%%%%%%%

\proclaim{Lemma 4.6}
Let the packing consisting of the triangles $\Delta ABC,
\Delta BA_1C$, $\Delta CB_2A$ satisfy the angular and side hypotheses and let 
$\gamma , \alpha _1 , \beta _2 \ge \pi /2$. Further let the
circumradius of $\Delta ABC$ be at most $p\sqrt{2}$.
Then the sum of the areas of these three
triangles is at least $p^2 \cdot (1 + \sqrt{3}/2) = p^2 \cdot 1.8660 \ldots $.
Equality holds if and only if $c_1=b_1=a_2=c_2=p$ and $\alpha _1 = \beta _2
= \pi /2$ and the circumradius of $\Delta ABC$ is $p{\sqrt{2}}$.
\endproclaim

%%%%%%%%%%%%%%%%%%%%%%%%%%%%%%%%%%%%%%%%%%%%%%%%%%%%%%%%%%%%%%%%%%%%%%%%%%
%%%%%%%%%%%%%%%%%%%%%%%%%%%%%%%%%%%%%%%%%%%%%%%%%%%%%%%%%%%%%%%%%%%%%%%%%%%

\demo{Proof} 
Consider an extremal configuration (this exists).

We have two angular hypotheses, namely $\alpha + \alpha _1 \le \pi $ and
$\beta + \beta _2 \le \pi $.
Analogously as in the proof of Lemma 4.2, we make the following case
distinctions. Either
\newline
(1) 
$\alpha + \alpha _1 = \beta + \beta _2 = \pi $, or
\newline
(2)
$\alpha + \alpha _1 , \beta + \beta _2 < \pi $ or, e.g.,
\newline
(3)
$\alpha + \alpha _1 = \pi > \beta + \beta _2$.

{\bf{1.}}
We begin with case (1). That is, the pentagon $ABA_1CB_2$ has a circumcircle,
and by $\gamma \ge \pi /2$ its circumcentre is not an interior point of it.
Applying Lemma 3.8 we get $V(ABA_1CB_2) \ge p^2 \cdot 2.3977 \ldots > 
p^2 \cdot (1 + \sqrt{3}/2) = p^2 \cdot 1.8660 \ldots $.

{\bf{2.}}
We turn to case (2).
We apply Lemma 3.1 to $\Delta BA_1C$ and $\Delta CB_2A$ 

\newpage

(in place of $\Delta
ABC$ in Lemma 3.1, with $ \pi - \alpha $ and $\pi - \beta $ in place of 
$\gamma $ in Lemma 3.1), 
obtaining
$$
c_1=b_1=a_2=c_2=p
\tag 4.6.1
$$
We apply Lemma 3.2 to $\Delta ABC$ (with the same notations here and there,
and $a_0:=a$ and $b_0:=b$) 
obtaining for the circumradius $R$ of $\Delta ABC$
$$
R=p{\sqrt{2}} 
\tag 4.6.2
$$
(observe that now by (2) we need not care the angular hypotheses, and the side
hypothesis remains preserved by $c \ge {\sqrt{a^2 + b^2}} \ge p{\sqrt{2}} >
p$).

{\bf{3.}}
Leaving the further investigation of case (2) later, for a while
we turn to case (3).
Analogously as in the proof of \thetag{3.8.7}, we have also here 
$$
b_1=c_1=p\,.
\tag 4.6.3 
$$
On the other hand, by Lemma 3.1 applied to $\Delta CB_2A$ (in place of $\Delta
ABC$ in Lemma 3.1, with $\pi - \beta $ in place of $\gamma $ in Lemma 3.1), 
we have 
$$
a_2=c_2=p\,.
\tag 4.6.4
$$
{\it{Hence \thetag{4.6.1} holds also now.}}

By 
$$
a \ge {\sqrt{c_1^2 + b_1 ^2}} \ge p{\sqrt{2}} {\text{ and similarly }} b \ge 
p{\sqrt{2}}
\tag 4.6.5
$$
we have 
$$
2R \ge c \ge {\sqrt{a^2+b^2}} \ge 2p \,.
\tag 4.6.6
$$
thus, also using the hypothesis of the lemma,
$$
R \in [p, p{\sqrt{2}}] /,.
\tag 4.6.7
$$
In {\bf{2}} we had {\thetag{4.6.2}}, in {\bf{3}} we have
obtained {\thetag{4.6.7}}.
This means that {\thetag{4.6.7}} holds in both of these cases.

{\bf{4.}}
Both in {\bf{2}} and {\bf{3}} we denote the central angles belonging to the
sides $a$ and $b$ in the circumcircle of $\Delta ABC$ by $\varphi _A$ and
$\varphi _B$. We denote the (smaller) central angle corresponding to a chord of
length $p$ of this circumcircle by $\varphi _0$; then
$$
\sin (\varphi _0 /2) = p/(2R) \,.
\tag 4.6.8
$$

\newpage

Then the minimum of $\alpha _1$ or of $\beta _2$ occurs when $a = p{\sqrt{2}}$
or $b = p{\sqrt{2}}$ (cf. {\thetag{4.6.5}}).
By $2p \le 2R$ we have $\varphi _A \le \pi $. Then $\alpha _1$, $a_1=a$ and 
$\varphi _A=2 \alpha $ increase together, so their maxima occur when $\alpha
+ \alpha _1$ increases to $\pi $, i.e., when $A_1$ gets to the boundary of the
circumcircle of $\Delta ABC$. Analogously, the maximum of $\beta _2$ occurs
when $B_2$ gets to the boundary of the circumcircle of $\Delta ABC$. 
That is, 
$$
\varphi _A , \varphi _B \in [2 \arcsin \left( (p{\sqrt{2}})/(2R) \right) ,
4 \arcsin \left( p /(2R) \right) ]
\tag 4.6.9
$$

Thus we have obtained that both in {\bf{2}} and in {\bf{3}}
all hypotheses 
of Lemma 3.9 are satisfied. Then the conclusion of Lemma 3.9
holds as well, i.e., the minimum of the total area of our three triangles is
attained (only) for 
$$
\varphi_A, \varphi _B \in \{ \varphi _{{\text{min}}} ,
\varphi _{{\text{max}}} \} = \{ 2 \arcsin \left( (p{\sqrt{2}})/(2R) \right) ,
4 \arcsin \left( p /(2R) \right) \} 
\tag 4.6.10
$$

{\bf{5.}}
The case 
$$
\varphi _A = \varphi _B = \varphi _{{\text{max}}}
\tag 4.6.11
$$
is covered by {\bf{1}}.

In the case
$$
\varphi _A = \varphi _B = \varphi _{{\text{min}}}
\tag 4.6.12
$$
we have that both $\Delta BA_1C$ and $\Delta CB_2A$ are isosceles right
triangles of legs of length $p$, hence have total area $p^2$. On the other
hand, by Lemma 3.2 (with the same notations there and here)
$V(\Delta ABC)$ is minimal when $a$ and $b$ are minimal, i.e.,
$a=b=p{\sqrt{2}}$, and the circumradius $R$ is maximal, i.e.,
$R=p{\sqrt{2}}$, when $V(\Delta ABC) = {\sqrt{3}}/2$, and thus the total area
of the three triangles is 
$$
p^2(1+{\sqrt{3}}/2) \,,
\tag 4.6.13
$$
proving the inequality of the lemma, with case of equality only as given in
the lemma, as follows from this proof.

There remains the case when one of $\varphi _A$ and $\varphi _B$ 
is $\varphi _{{\text{min}}}$,
and the other one is $\varphi _{{\text{max}}}$. Then we apply Lemma 3.10,
obtaining
$$
\cases
V(\Delta ABC) + V(\Delta BA_1C) + V(\Delta CB_2A) \ge \\
p^2 \cdot 1.8720 \ldots 
> p^2(1+{\sqrt{3}}/2) = p^2 \cdot 1.8660 \ldots \,.
\endcases
\tag 4.6.14
$$
$\blacksquare $
\enddemo

%%%%%%%%%%%%%%%%%%%%%%%%%%%%%%%%%%%%%%%%%%%%%%%%%%%%%%%%%%%%%%%%%%%%%%%%%%
%%%%%%%%%%%%%%%%%%%%%%%%%%%%%%%%%%%%%%%%%%%%%%%%%%%%%%%%%%%%%%%%%%%%%%%%%%%

\newpage

\head{\S 5 Proofs of the results from \S 2}\endhead

%%%%%%%%%%%%%%%%%%%%%%%%%%%%%%%%%%%%%%%%%%%%%%%%%%%%%%%%%%%%%%%%%%%%%%%%%%
%%%%%%%%%%%%%%%%%%%%%%%%%%%%%%%%%%%%%%%%%%%%%%%%%%%%%%%%%%%%%%%%%%%%%%%%%%%

Now we prove the statements from \S 2, practically in the inverse order, as
they were stated in \S 2.

%%%%%%%%%%%%%%%%%%%%%%%%%%%%%%%%%%%%%%%%%%%%%%%%%%%%%%%%%%%%%%%%%%%%%%%%%%
%%%%%%%%%%%%%%%%%%%%%%%%%%%%%%%%%%%%%%%%%%%%%%%%%%%%%%%%%%%%%%%%%%%%%%%%%%%

\demo{Proof of Theorem 2.9}
We consider an $(r,R)$-system 
$P \subset {\Bbb{R}}^2$ satisfying the hypotheses of the theorem. There will
be no problems with the non-obtuse triangles in the Delone
triangulation, then case (1) of Theorem 2.9
could work. However, as we will see later,
they will not be always settled via case (1) of Theorem 2.9.

To handle the obtuse triangles in the Delone triangulation, we
define an oriented graph, whose vertices are the triangles
of the Delone triangulation of ${\Bbb{R}}^2$. The oriented edges pass
from an obtuse triangle to the triangle joining to it along the longest
side of the obtuse triangle. Thus between any two triangles at most one
oriented edge passes (since
in a Delone triangulation the sum of the angles opposite to the common
edge of two triangles is at most $\pi $). The number of edges starting from a
triangle is at most $1$ and the number of edges ending in a triangle is at
most $3$.

Therefore, if the corresponding non-oriented graph contained a cycle, our
oriented graph would contain an oriented cycle. In fact, if we have a
non-oriented cycle $T_1 \ldots T_m$ (with cyclic notation), then to an
oriented edge $\overrightarrow{T_iT_{i+1}}$ 
there cannot join an oriented edge $\overrightarrow{T_iT_{i-1}}$,
but only an oriented edge $\overrightarrow{T_{i-1}T_i}$. 
Repeat these considerations for the
oriented edge $\overrightarrow{T_{i-1}T_i}$, etc., and we get our claim.

However, an oriented cycle in our graph is impossible. Namely, if an obtuse
triangle $T_1$ joins by its longest side to a side of another obtuse triangle
$T_2$, then the common side is opposite to an acute angle in $T_2$. Then the
square of the longest side of $T_2$ is at least the square of the longest side
of $T_1$ (which is at least $(2r)^2$) plus $(2r)^2$. Thus a cycle in our graph
is impossible, since passing on the oriented 
edges the length of the longest side of the triangle
strictly increases, and returning to the starting triangle we obtain a
contradiction.

Since each triangle in the Delone triangulation has a
circumradius at most $R$, hence longest side at most 
$2R$, similarly we get that an oriented path has a length at most $7$. In
fact, if there were an oriented path $\overrightarrow{T_1T_2 \ldots T_{9}}$,
then by induction
the longest side of $T_i$ would have square at least $i(2r)^2$. In
particular for $i=9$ the square of the longest side were at least
$9(2r)^2$, but on the other hand, as we have just seen,
it is at most $(2R)^2$. Thus $9(2r)^2 \le
(2R)^2$, i.e., $3 \le R/r \le 2{\sqrt{2}}$, a contradiction to the hypotheses
of our theorem.

By the above facts we see that the connected components of the corresponding
non-oriented graph contain a bounded number of triangles. In fact, each
oriented edge path has an endpoint, and that is 
a non-obtuse triangle $T_0$. Then $T_0$ is the
endpoint of at most three oriented edges, from obtuse triangles $T_{1i}$.
Repeating, each $T_{1i}$ is the
endpoint of at most two oriented edges, from 

\newpage

obtuse triangles $T_{2ij}$.
Continuing similarly, since the lengths of oriented paths are at most $7$, a
connected component of the corresponding unoriented graph can have a size at
most $1+3 \cdot (1+2+ \ldots +2^6)$. 

A non-trivial
connected component (i.e., which contains at least two triangles)
of the corresponding non-oriented graph (later we will
say just ``non-trivial component'') is thus a tree,
and contains just one non-obtuse triangle, namely the one from which no edge of
the graph starts. On the other hand, there are triangles in the non-trivial
component, at which no edge ends, but there starts just one edge. 
These will be divided in classes, two such triangles being
in the same class if the oriented edge starting from them end at the same 
triangle.

Next we divide all triangles of a non-trivial
component in classes. An above class of
triangles, in which no edge ends, together with the triangle which is 
the
common endpoint of the edges starting from them will form a new class. All
other triangles in the same non-trivial
component, not contained in any of the new classes,
i.e., which are the end-points of oriented paths of length (at least) 
$2$, with the exception of the triangle from which no edge starts 
(unless if it is in a new class),
will form one-element new classes. (In fact, we will use only that they
are endpoints of an oriented edge, cf. Lemma 4.4, which will be applied to
these triangles.)

Each Delone 
triangle, either the triangle in a non-trivial component from which no edge
starts (provided that it is not in a new class of triangles), 
or a triangle not contained in any non-trivial
component of our graph, which triangle is therefore in both cases
non-obtuse, will form a one element class. Thus now all the Delone
triangles are divided to classes.

Observe that the diameter of the
union of any class is at most twice the maximal possible diameter of a Delone
triangle, i.e., is at most $4R$. By boundedness of the diameters, one can
estimate from below 
the average area of a Delone triangle in each class separately, and the same
lower bound will be valid on the whole plane.

The hypotheses of our theorem imply the hypotheses of each of Lemmas 4.1 to
4.6, with $p:=2r$. By these lemmas, the average area of Delone 
triangles in one class is at least $p^2/2 = 2r^2$. In fact,
one element classes in a non-trivial component, except
the triangle from which no edge starts (unless it is in a new class), 
are settled by Lemma 4.4, giving cases (1) and (2) of the theorem. 
Two element classes are settled in Lemmas 4.1 and
4.5, giving cases (3) and (4) of the theorem.
Three element classes are settled in Lemmas 4.2 and 4.6, giving cases (3) and
(4) of the theorem.
Four element classes are settled in Lemma 4.3, giving case (3) of the
theorem. (Except Lemma 4.1 we have even
that the average area is strictly greater than $2r^2$.)
Last, if a class is a single
non-obtuse triangle, either the triangle in a non-trivial component from which
no edge starts (if it is not in a new class), 
or a triangle not contained in any non-trivial component of our graph,
then the area of this triangle is at least $V_0$. Therefore the 
average area of all 

\newpage

Delone triangles is at least $\min \{ V_0, 2r^2 \} $. 

This bound is sharp. In fact, for any $\varepsilon > 0$ there is a non-obtuse
Delone triangle $T$ with area less than $V_0 + \varepsilon $. Now consider the
point lattice generated by the vertices of $T$. Then in the Delone triangulation
of this point lattice
all Delone triangles are congruent to $T$, hence all of them have area 
less than $V_0 + \varepsilon $, so the average area is also
less than $V_0 + \varepsilon $. The average area $2r^2$, even all Delone
triangles having area $2r^2$ can be attained for a
square lattice of side length $2r$, where the respective covering radius $R$ is
$r{\sqrt{2}}\,\,(<r \cdot 2{\sqrt{2}})$.
$\blacksquare $
\enddemo

%%%%%%%%%%%%%%%%%%%%%%%%%%%%%%%%%%%%%%%%%%%%%%%%%%%%%%%%%%%%%%%%%%%%%%%%%%
%%%%%%%%%%%%%%%%%%%%%%%%%%%%%%%%%%%%%%%%%%%%%%%%%%%%%%%%%%%%%%%%%%%%%%%%%%%

\demo{Proof of Theorem 2.8}
{\bf{1.}}
By hypothesis we have
$$
1/{\sqrt{2}} \le m(L) \le M(L)=1\,.
\tag 2.8.1
$$

Clearly we may suppose that our packing of translates of $L$ in ${\Bbb{R}}^n$
is saturated, i.e., we cannot add any new translate of $L$ to the packing
without violating the packing property. This implies that any empty circle
associated to the point system on the $x_1x_2$-coordinate plane, consisting of
the intersection points of the axes of rotation of the translates of $L$ in our
packing with the $x_1x_2$-coordinate plane,
has a radius strictly smaller than $2$. In fact, if one such radius would be
at least $2$, then we could add to our packing of translates of $L$ even a
translate of the cylinder $L_0$ with base a circle of radius $1$ and centre
the origin in the $x_1x_2$-coordinate plane, and with axis the $x_3 \ldots
x_n$-coordinate plane, without violating the packing property. Then this
translate of $L_0$ would contain some translate of $L$, and addition of this
translate of $L$ to our packing of translates of $L$ would not violate the
packing property, contradicting saturatedness of our packing of translates of
$L$.

Later we will use only that each above empty circle has a radius at most $2$.
 
So we may restrict our attention to saturated packings in this sense.
Now consider the intersection points of the axes of rotation of the
translates of $L$ consisting 
our packing with the $x_1x_2$-coordinate plane. We consider the Delone
triangulation of the $x_1x_2$-coordinate plane 
belonging to this system of points (if there are
Delone-polygons with larger numbers of sides, they are triangulated in an
arbitrary way). By saturatedness any triangle in this triangulation has a
circumradius at most (actually, smaller than) $2$. 
The number density of this system of points is
$1/(2\overline{V})$, where $\overline{V}$ is the average area of the
Delone-triangles. (For this observe that the total angle sum of all Delone
triangles in some large circle about the origin is about $2 \pi $ times
the total number of vertices.) So we have to estimate ${\overline{V}}$ from
below.

{\bf{2.}}
Let $\Delta ABC$ be a non-obtuse triangle of our Delone
triangulation. We are going
to show that 
$$
V(\Delta ABC) \ge V_0(L) {\text{ (cf. \thetag{2.8}).}} 
\tag 2.8.2
$$

\newpage

The points $A,B,C$ are the intersection points 
of the axes of rotation of suitably
translated copies of $L$ with the $x_1x_2$-coordinate plane.
Let these translated copies be $L + (x_{11},x_{12}, \ldots ,x_{1d})$, 
$L + (x_{21},x_{22}, \ldots ,x_{2n})$ and 
$L + (x_{31},x_{32}, \ldots ,x_{3n})$. 
Then by the definition of the function $g$ we have
$$
\cases
|AB|= {\sqrt{(x_{21} - x_{11})^2 + (x_{22} - x_{12})^2}} \\
\ge 2g(x_{23} - x_{13}, \ldots, x_{2n} - x_{1n})\,.
\endcases
\tag 2.8.3 
$$ 
The analogous inequalities hold also for $|BC|$ and $|CA|$.

That is, the side lengths of the
triangle $\Delta ABC$ are at least certain values of the function $g$.
By \thetag{2.8.1} and \thetag{2.7} 
the values of the function $g$ lie in the interval $[ {1/\sqrt{2}}, 1]$. 
We define the triangle $T$ as follows.
$$
\cases
{\text{Denote }} T {\text{ a triangle with side lengths }}
2g(x_{23} - x_{13}, \ldots, x_{2n} - x_{1n}), \\
2g(x_{33} - x_{23}, \ldots, x_{3n} - x_{2n}) {\text{ and }}
2g(x_{13} - x_{33}, \ldots, x_{1n} - x_{3n}) \,.
\endcases
\tag 2.8.4
$$
Observe that a triangle with any given side lengths,
with maximal side$/$minimal side $\le {\sqrt{2}}$ exists (by the triangle
inequality) and is non-obtuse (by the cosine law), 
thus $T $ is a non-obtuse triangle.
Then using \thetag{3.0} we have
$$
V(\Delta ABC) \ge V(T)\,.
\tag 2.8.5
$$ 

By rotational symmetry of $L$ with respect to the
$x_3 \ldots x_n$-coordinate plane, we can change the first two coordinates of
the translation vectors $v_1:=(x_{11},x_{12}, \ldots ,$
\newline
$x_{1n})$, 
$v_2:=(x_{21},x_{22}, \ldots ,x_{2n})$ and 
$v_3:=(x_{31},x_{32}, \ldots ,x_{3n})$, thus
obtaining new translation vectors 
$v_1':=(x_{11}',x_{12}', x_{13}, \ldots ,x_{1d})$, 
$v_2':=(x_{21}',x_{22}', x_{23}, \ldots ,x_{2n})$ and 
$v_3':=(x_{31}',x_{32}', x_{33}, \ldots ,x_{3n})$ such that the respective new
translated copies of $L$ are pairwise touching
(in the sense precised in (2.5)). The
projections of $v_1', v_2', v_3'$ to the $x_1x_2$-coordinate plane are denoted
by $A',B',C'$.
Let $T':=\Delta v_1'v_2'v_3'$; its projection to the $x_1x_2$-coordinate
plane is $\Delta A'B'C'$, which is congruent to $T$. We may and will
suppose that actually 
$$
\Delta A'B'C' = T \,.
\tag 2.8.6
$$

By the definition of $V_0(L)$ (otherwise said,
varying the coordinates $x_{13}, \ldots ,x_{1n}$, 
$x_{23}, \ldots ,x_{2n}$ and $x_{33}, \ldots ,x_{3n}$ of the translation
vectors in an arbitrary way, and then taking infimum over these coordinates),
we obtain that 
$$
V(T) \ge V_0(L) \,.
\tag 2.8.7
$$
Then \thetag{2.8.5} and \thetag{2.8.7} together give \thetag{2.8.2}.

\newpage

{\bf{3.}}
Clearly
$$
\cases
{\text{the isosceles right triangle with legs of length }} 2 m(L) {\text{
has}} \\
{\text{at least such an area, as the isosceles triangle with sides }} \\
2m(L), 2m(L), 2M(L) \,.
\endcases
\tag 2.8.8
$$
In formula,
$$
2m(L)^2 \ge M(L)\sqrt{\left( 2m(L) \right) ^2 - M(L)^2} 
\tag 2.8.9
$$
(this also follows from the arithmetic-geometric mean inequality).

Now we show that,
$$
\cases
{\text{for any }} \varepsilon >0{\text{, for a suitable two-dimensional
lattice packing of}} \\
{\text{translates of }} L 
{\text{ generated by three mutually touching translates of }} \\
L
{\text{ (in the sense of \thetag{2.5}), a non-obtuse Delone triangle on the}} \\
x_1x_2{\text{-co\-or\-di\-na\-te plane can have an area at most }} \\ 
M(L) \sqrt{\left( 2m(L) \right) ^2 - M(L)^2} + \varepsilon \,.
\endcases
\tag 2.8.10
$$ 

Let us consider the following
translates of $L$: $L+(M(L),0, \ldots ,0)$, $L-(M(L),0, $
\newline
$\ldots ,0)$ (these are
touching translates in the sense of \thetag{2.5}) 
and $L+(0,x_2,x_3, \ldots , x_n)$, where $x_2$ is chosen
so that $L$ and any of $L \pm (0,x_2,x_3, \ldots , x_n)$ should be touching
(in the sense of \thetag{2.5}),
for any given $(x_3, \ldots , x_n) \in {\Bbb{R}}^{n-2}$. (This is possible by
the rotational symmetry of $L$ about the $x_3 \ldots  x_n$-coordinate plane.)
Then the intersection points 
of the axes of rotation (translates of the $x_3 \ldots
x_n$-coordinate plane) of these three translates of $L$
with the $x_1x_2$-coordinate plane are the vertices of an isosceles triangle
$T(\delta )$,
with base of length $2M(L)$ and equal sides of lengths $2g(x_3, \ldots , x_n)$.
By a suitable choice of $x_3, \ldots , x_n$ we can attain that 
$g(x_3, \ldots , x_n) - m(L) \in [0, \delta )$, 
which for suitably small $\delta > 0$ implies
$$
V \left( T(\delta ) \right) - 
M(L) \sqrt{\left( 2m(L) \right) ^2 - M(L)^2} \in [0,  \varepsilon ) \,,
\tag 2.8.11
$$
proving the inequality in \thetag{2.8.10}, {\it{provided $T(\delta )$ is a 
Delone triangle, for a suitable packing of translates of $L$}}. 
 
Still we have to show that $T(\delta )$ is a non-obtuse Delone triangle, 
for a suitable packing of translates of $L$. 
The base of $T(\delta )$ has length $2M(L)$, and its
equal sides have length in $[2m(L), 2 \left( m(L) + \delta \right) ]$. By
$M(L)/m(L) \le {\sqrt{2}}$ this triangle is non-obtuse. There remains to show
that $T(\delta )$ is a Delone triangle for a suitable packing of translates of
$L$.

\newpage

{\bf{4.}}
We consider the (inhomogeneous) two-dimensional
lattice of translates of $L$ spanned by the translates of $L$ by the three new
translation vectors $v_1', v_2', v_3'$ from {\bf{2}}. 
Any two of these three translates of
$L$ are touching. Then, in the spanned body lattice, any two translates of $L$,
which are simultaneous translates of any two of these three mutually
touching translates, are also touching. 
If $L$ and $L+v'$ are two translates of $L$ from this body lattice which are
not simultaneous translates of any two of our three mutually touching
translates $L+v_1', L+v_2', L+v_3'$, then the following holds. If we decompose
our two-dimensional 
lattice of translates of $L$ to one-dimensional lattices of
translates of $L$ 
spanned by suitable two of our three mutually touching translates of $L$, then 
$L$ and $L+v'$ are not in the same above 
one-dimensional lattice of translates of
$L$, but also not in neighbourly above 
one-dimensional lattices of translates of $L$.

The distance of the orthogonal projections to the $x_1x_2$-coordinate plane of
the straight lines spanned by homologous points in translates
of $L$ for two neighbourly above one-dimensional lattices of translates of $L$,
is the height of
a triangle with vertices the points of intersection of the axes of rotation
of three mutually touching copies of $L$,
with sides in $[2m(L), 2M(L)] = [2m(L), 2] \,\,(\subset 
[{\sqrt{2}}, 2])$ with the $x_1x_2$-ccordinate plane. 
The minimal height belongs to the
maximal side (for which the adjacent angles are acute), 
and for given length of the maximal side, it 
is attained if the other two sides are minimal possible, i.e., have length 
$2m(L)$. 
Hence the height is at least $\sqrt{\left( 2m(L) \right) ^2 - M(L)^2}
=  \sqrt{\left( 2m(L) \right) ^2 - 1} \ge \sqrt{2-1} = 1$. 

Therefore the distance of the orthogonal projections to the
$x_1x_2$-coordinate plane of the
straight lines spanned by homologous points in translates
of $L$ for two neighbourly above one-dimensional lattices of translates of $L$
is at least $1$. 
This implies that the distance of the orthogonal projections to the
$x_1x_2$-coordinate plane of the
straight lines spanned by homologous points in translates
of $L$ for two at least second neighbour above 
one-dimensional lattices of translates of $L$ is at least $2\,\,\,\,$($= \sup
2g$, cf. \thetag{2.7}).
Therefore the translates 
$L$ and $L+v'$ are disjoint. This shows that our two-dimensional lattice
arrangement of translates of $L$ is in fact a two-dimensional lattice packing,
spanned by some three mutually touching translates of $L$ (in the sense
of \thetag{2.5}).

Then the points of intersection of the axes of rotation
of all translates of $L$ in the
two-dimensional lattice packing of translates of $L$ is the point lattice in
the $x_1x_2$-coordinate plane generated by the projections of the three new
translation vectors $v_1', v_2', v_3'$ from {\bf{2}} to the
$x_1x_2$-coordinate plane, i.e., the points $A',B',C'$, which are the vertices
of $T(\delta )$ by \thetag{2.8.6} (observe that the role of 

\newpage

$T$ is now taken
over by $T(\delta )$).

The corresponding two-dimensional point lattice projects orthogonally to the
$x_1x_2$-plane injectively onto a two-dimensional point lattice, generated by
$A',B',C'$ in the $x_1x_2$-coordinate plane, since $m(L) > 0$. 
Moreover, the Delone triangulation
of the $x_1x_2$-coordinate plane corresponding to this point lattice 
has as Delone triangles $T(\delta )=A'B'C'$ 
and all its lattice translates (by this
point lattice) and $A' + B' - T(\delta )$ 
and all its lattice translates (by this point lattice). Therefore the 
number density of this point lattice in the $x_1x_2$-coordinate plane
is 
$$
\cases
1/\left( 2 V\left( T(\delta ) \right) \right) \in 
\left( \left( 
2 M(L){\sqrt{\left( 2 m(L) \right) ^2 - M(L)^2}} + 2 \varepsilon \right) ^{-1}, 
\right. 
\\
\left. 
\left( 2 
M(L){\sqrt{\left( 2 m(L) \right) ^2 - M(L)^2}} \right) ^{-1} \right] = \\
\left( 
\left( 2 {\sqrt{4 \left( m(L) \right) ^2 - 1}} + 2 \varepsilon \right) ^{-1},
\left( 2 {\sqrt{4 \left( m(L) \right) ^2 - 1}} \right) ^{-1} \right]
\endcases
\tag 2.8.12
$$

All this shows the second statement of the theorem, taking in
consideration \thetag{2.9}.

{\bf{5.}}
\thetag{2.8.8}, \thetag{2.8.9} and \thetag{2.8.10} imply 
that 
$$
2m(L)^2 \ge M(L){\sqrt{\left( 2m(L) \right) ^2 - M(L)^2}} =
\inf _{\delta >0} V \left( T(\delta ) \right) \ge V_0(L) \,,
\tag 2.8.13
$$
which in turn implies
$$
\min \{ V_0(L), 2m(L)^2 \} = V_0(L) \,.
\tag 2.8.14
$$
Hence, for any packing of translates of $L$,
the number density of the intersections
of the axes of rotation of the translates of $L$ in the packing
is by Theorem 2.9 at most 
$$
1/ \min \{ 2V_0(L), 4m(L)^2 \} = 1/ \left( 2V_0(L) \right) \,,
\tag 2.8.14
$$ 
which implies by hypothesis (2.9) the theorem.
$\blacksquare $
\enddemo

%%%%%%%%%%%%%%%%%%%%%%%%%%%%%%%%%%%%%%%%%%%%%%%%%%%%%%%%%%%%%%%%%%%%%%%%%%%
%%%%%%%%%%%%%%%%%%%%%%%%%%%%%%%%%%%%%%%%%%%%%%%%%%%%%%%%%%%%%%%%%%%%%%%%%%%

JON

\demo{Proof of Proposition 2.6} 
This follows from Proposition 2.2. We claim that the packing
density of translates of $L$ is at most the density of a body lattice of
translates of $K$. Here the corresponding point lattice is the sum 

\newpage

of the
$(n-2)$-dimensional point lattice in the $x_3 \ldots x_n$-coordinate plane from
this Proposition, and a two dimensional point lattice that intersects the
$x_3 \ldots x_n$-coordinate plane only in $\{ 0 \} $ (by $\inf g >0$). Still
recall that there exists a densest lattice packing of translates of $K$.
$\blacksquare $
\enddemo

%%%%%%%%%%%%%%%%%%%%%%%%%%%%%%%%%%%%%%%%%%%%%%%%%%%%%%%%%%%%%%%%%%%%%%%%%%
%%%%%%%%%%%%%%%%%%%%%%%%%%%%%%%%%%%%%%%%%%%%%%%%%%%%%%%%%%%%%%%%%%%%%%%%%%%

\demo{Proof of Proposition 2.5} 
Let us consider the Dirichlet-Voronoj (DV) decomposition 
of the $x_3 \ldots x_n$-coordinate
plane associated to the centres of our balls occurring in the packing with
union $L$.
Then the DV-cell of 
any of these unit balls is contained in the concentric ball
of radius ${\sqrt{2}}$. 

Clearly we have $\sup g = 1$. Next we show that $\inf g \ge 1/\sqrt{2}$.
Let $L$ and $L+(x_1, \ldots ,x_n)$ 
be two touching copies of translates of $L$. Let $B_0$
be one of the balls constituting the packing $L+(x_1, \ldots ,x_n)$, 
with centre $(x_{10}, \ldots ,x_{n0})$. 
We project this centre orthogonally to the $x_3 \ldots x_n$-coordinate plane, 
obtaining
the point $(0,0,x_{30}, \ldots ,x_{n0})$. 
This last point is in the DV-cell of the centre $(0,0,x_{31},$
\newline 
$\ldots ,x_{n1})$
of one of the balls constituting the packing $L$. Hence by hypothesis of the
Proposition the distance of $(0,0,x_{30}, 
\ldots ,x_{n0})$ and $(0,0,x_{31}, \ldots ,x_{n1})$ is at most
$\sqrt{2}$. By Pythagoras theorem, the distance of the 
$x_3 \ldots x_n$-coordinate plane and $(x_{10}, \ldots ,x_{n0})$ is at least
$\sqrt{2}$. Hence 
$\inf g \ge \sqrt{2}$, and therefore $1=M(L) \ge m(L) \ge 1/{\sqrt{2}}
$ (cf. \thetag{2.7}). 
Now applying Theorem 2.8 this proposition is proved.
$\blacksquare $
\enddemo

%%%%%%%%%%%%%%%%%%%%%%%%%%%%%%%%%%%%%%%%%%%%%%%%%%%%%%%%%%%%%%%%%%%%%%%%%%
%%%%%%%%%%%%%%%%%%%%%%%%%%%%%%%%%%%%%%%%%%%%%%%%%%%%%%%%%%%%%%%%%%%%%%%%%%

\demo{Proof of Proposition 2.4}
{\bf{1.}} By hypothesis \thetag{2.2} the 
density of the packing consisting of some translated copies of
$L$ is proportional to the number density of the intersection of their
rotation axes with the $xy$-plane. We have $\sup g = 1$. We are going to show
that $\inf g \ge 1/\sqrt{2}$. Let us consider two different 
translates of $L$ from our packing of translates of $L$, 
say, $L+(x_1,y_1,z_1)$ and $L+(x_2,y_2,z_2)$. 
Let us consider a segment
$[z-R,z+R] \subset {\Bbb R}$. Let us project the centres of the balls
consisting $L+(x_1,y_1,z_1)$ and $L+(x_2,y_2,z_2)$ to the $z$-axis. Then
the number of the projection points in $[z-R,z+R]$ is $2R/(2d) + o(R)$, where 
$o(R)$ is uniform for all $z$'s and all 
translates $L+(x_1,y_1,z_1)$ and $L+(x_2,y_2,z_2)$ consisting our packing of
translates of $L$.
Then the number of projection points (possibly with multiplicities)
of the centres of balls consisting 
$L+(x_1,y_1,z_1)$ and $L+(x_2,
y_2,z_2)$ is $2R/d + o(R)$. Hence some neighbouring projections
have a distance is at most $d + o(1) \le \sqrt{2} + o(1)$, 
with $o(1)$ uniform (as above with $o(R)$). 
Then both projections cannot come from
$L+(x_1,y_1,z_1)$, or from $L+(x_2,
y_2,z_2)$, since they are unions of packings of unit balls, and
then this distance should be at least $2$. Hence these neighbouring
projections come one from $L+(x_1,y_1,z_1)$, and the
other one from $L+(x_2,y_2,z_2)$. By Pythagoras theorem,
the distance of the axes of rotation of 
$L+(x_1,y_1,z_1)$ and 
$L+(x_2,y_2,z_2)$ is a least $\sqrt{2} + o(1)$. Letting
$R \to \infty $, we get that the distance of these axes of rotation is at
least $\sqrt{2}$. 

Thus the hypotheses of Theorem 2.8 are satisfied, hence also its
conclusion is satisfied, i.e., the statement of Proposition 2.4 is proved.
$\blacksquare $
\enddemo

%%%%%%%%%%%%%%%%%%%%%%%%%%%%%%%%%%%%%%%%%%%%%%%%%%%%%%%%%%%%%%%%%%%%%%%%%%%
%%%%%%%%%%%%%%%%%%%%%%%%%%%%%%%%%%%%%%%%%%%%%%%%%%%%%%%%%%%%%%%%%%%%%%%%%%

\newpage

\demo{Proof of Theorem 2.1} 
By $d \le \sqrt{2}$ we have that the concentric $1$-dimensional
balls of radius $\sqrt{2}$ cover ${\Bbb R}$, hence we may apply Proposition 2.5.
Thus the density of our packing of translates of $L$ is at most the maximal
density of the two-dimensional lattice packings of $L$,
where the corresponding point lattice projects orthogonally
to the $xy$-plane injectively, onto a two-dimensional
point lattice in the $xy$-plane. (This maximum exists.)
We have to prove that these two-dimensional lattice packings of
translates of $L$ have a density at most $\pi /(3d\sqrt{3-d^2})$, and the
unique densest such lattice packing is the one given in this Theorem.

Evidently in a
densest lattice packing of translates of $L$ there are two translates of $L$
touching each
other. These together generate a two-dimensional lattice $ B^3 + \Lambda 
:=\{ B^3 + \lambda \mid \lambda \in \Lambda \} $
of balls (recall that $B^3$ is the closed unit ball centred at the origin). 
We may assume that the $2$-plane $\Pi $ spanned by the centres of these balls
is the $xy$-plane.

Join the consecutive centres of the balls in translates of $L$ 
belonging to $B^3 + \Lambda $
by segments, and also join the centres of each pair of balls
belonging to two neighbouring translates of $L$ in $B^3 + \Lambda $ 
and touching each other. Thus the $xy$-plane 
is subdivided to a lattice of parallelograms of sides $2d$
and $2$. 

Still draw the diagonals of these parallelograms joining two obtuse 
angles, thus cutting these parallelograms into two triangles. 
If the parallelogram is a rectangle, we draw one of its diagonals
(chosen parallel to each other). Thus we
obtain a tiling ${\Cal{T}}$ of the $xy$-plane to triangles which are
translates and centrally symmetric images of some fixed triangle from
${\Cal{T}}$.

Let $\Delta ABC$ be one of the triangles from ${\Cal{T}}$, 
with $A$ and $B$ being the
centres of balls belonging to the same translate of $L$, and with $|BC|=2$. 
Then, denoting the angles of 
$\Delta ABC$ at $A,B,C$ by $\alpha , \beta , \gamma $, we have 
$$
\alpha ,\beta \le \pi /2 \,.
\tag 2.1.1
$$
We have also 
$$
\gamma \le \pi /2
\tag 2.1.2
$$ 
since else 
${\sqrt{2}} < |AB|/\min \{ |BC|, |CA| \} = |AB|/|BC| = d$, contradicting the
hypothesis of the theorem. 

Our lattice of unit balls decomposes to horizontal layers which are translates
of $\{ B^3 + \lambda 
\mid \lambda \in \Lambda \} $. We determine how close a neighbouring horizontal
layer can be to $B^3+\Lambda
$ (i.e., how close their mid-planes can be). 
Let $E$ be the centre of a unit ball in a neighbouring horizontal
layer. Then its
distance to the centre of each ball in $\Lambda $ is at least $2$. Let $E'$ be
the orthogonal projection of $E$ to the $xy$-plane. 
We want to determine the minimum of $|EE'|$. 

\newpage

We may suppose that $E'$ lies in the above triangle $\Delta ABC \in {\Cal{T}}$.
Then 
$$
|EA|,|EB|,|EC| \ge 2 \,.
\tag 2.1.3
$$
Observe that we have right triangles $\Delta AEE'$, $\Delta BEE'$, $\Delta
CEE'$, whose hypotenuses $EA$, $EB$, $EC$ have lengths
at least $2$ (by \thetag{2.1.3}). Then, by Pythagoras' theorem,
$$
\cases
|EE'| \ge \max \{ {\sqrt{2^2-|E'A|^2}}, {\sqrt{2^2-|E'B|^2}}, 
{\sqrt{2^2-|E'C|^2}} \} \\ 
= {\sqrt{4- (\min \{ |E'A|,|E'B|,|E'C| \} )^2 }} \,.
\endcases
\tag 2.1.4
$$
Since by \thetag{2.1.1} and \thetag{2.1.2} $\Delta ABC$ is a non-obtuse
triangle, with the circumradius $R$ of $\Delta ABC$ we have
$$
\cases
\min \{ |E'A|,|E'B|,|E'C| \} \le R \,,
{\text{ with equality if and only if}}\\
 E' {\text{ is the circumcentre }} c {\text{ of }} \Delta ABC \,.
\endcases
\tag 2.1.5
$$
Thus \thetag{2.1.4} implies
$$
\cases
|EE'| \ge {\sqrt{4-R^2}} \,,
{\text{ with equality if and only if}} \\
E'=c {\text{ and }} |EA|=|EB|=|EC| =2 \,.
\endcases
\tag 2.1.6
$$
In what follows we assume that 
$$
\cases
E {\text{ and }} E' {\text{ are the points for
which \thetag{2.1.5} and \thetag{2.1.6}}} \\
{\text{become equalities, i.e., }} E'=c {\text{ and }} |EA|=|EB|=|EC| =2 \,.
\endcases
\tag 2.1.7
$$

Still we have to show that 
$$
A,B,C {\text{ and the above chosen }} E {\text{ generate a lattice 
{\rm{packing}}}}.
\tag 2.1.8
$$
In a horizontal layer the
corresponding open unit balls are disjoint. Two open unit balls in neighbourly
horizontal
layers are disjoint by the following reason. The minimal distance of $E$ and
any lattice point in $\Lambda $ is attained if and only if the minimal
distance of $E'$ and any lattice point in $\Lambda $ is attained.
Observe that the DV-cells of $\Lambda $ in the $xy$-plane satisfy the
following: $\Delta ABC$ is covered by the DV-cells of $A,B,C$, whose common
vertex is $c$. Thus $c$ is closest
among any points of $\Lambda $ to $A,B,C$, and this minimal
distance is $R$. In particular, all other lattice points of $\Lambda $ have a
distance at least $R$ from $E'=c$. Then also $E$ has a distance at least 
${\sqrt{R^2+|EE'|^2}}=2$ from any lattice points of $\Lambda $.

It remains to show that the open unit 
balls in at least second neighbour horizontal layers are disjoint. For this it
is sufficient to show that the height $|EE'|$ 
of the tetrahedron $ABCE$ corresponding to the vertex $E$ 
is at least $1$. 

\newpage

Since $|EE'|={\sqrt{2^2-R^2}}$, therefore 

$$
{\text{we have to prove that }} R \le {\sqrt{3}} \,. 
\tag 2.1.9
$$

Observe that $c$ lies on the perpendicular bisectors of the sides $AB$
and $BC$ of $\Delta ABC$. 
If $|AC|$ increases from its minimum $2$ till its maximum
$2{\sqrt{d^2+1}}$ then the
angle $\beta = \angle ABC$ strictly increases, and by elementary geometric
considerations $c$ moves on the perpendicular
bisector of side $AB$ farther from side $AB$. 
Then $c$ will be the farthest from side 
$AB$ when $|AC|$ attains its maximum $2{\sqrt{d^2+1}}$, which happens if
$\Lambda $ is
a rectangular lattice on the $xy$-plane, when by $d \in [1, {\sqrt{2}}]$
$$
R={\sqrt{d^2+1}} \le {\sqrt{3}} \,, 
\tag 2.1.10
$$
proving our claim \thetag{2.1.9} and thus also \thetag{2.1.8} .

Now we investigate the volume of the basic parallelepiped of our lattice.
We consider the tetrahedron $ABCE$. Its face $\Delta BCE$ is a regular
triangle of side $2$, i.e., $|BC|=|CE|=|EB|=2$,
and $|AE|=2$, $|AB|=2d$ and lastly $2 \le |AC| \le
2{\sqrt{d^2+1}}$. Thus the faces $\Delta BCE$ and $\Delta ABE$ are given up to
congruence. They
join at their common edge $BE$. If their angle is $\varphi $, then the
volume $V(ABCE)$
of our tetrahedron is proportional to $\sin \varphi $. Now we calculate
$|AC|$. Let us denote the projections of $A$ and $C$ to the line $BE$ by $A'$
and $C'$. Then $|AC|^2= |A'C'|^2 + |A'A|^2 + |C'C|^2 
-2 |A'A| \cdot |C'C| \cdot \cos \varphi $, 
thus $|AC|$ is a strictly increasing function of $\varphi $. We have $2 \le 
|AC| \le 2 {\sqrt{d^2+1}}$. We denote the values of $\varphi $ belonging to 
$|AC| = 2$ and $|AC| = 2 {\sqrt{d^2+1}}$ by $\varphi _{\min }$ and $\varphi
_{\max }$. Then we have $\varphi _{\min } \le \varphi _{\max }$. 
Since the sine function is strictly concave on $[0, \pi ]$, on the interval
$[\varphi _{\min }, \varphi _{\max }]$ the function $\sin \varphi $, and thus
also $V(ABCE)$ attains its minimum either at $\varphi = \varphi _{\min }$, or
at $\varphi = \varphi _{\max }$, but not at any $\varphi \in 
(\varphi _{\min }, \varphi _{\max })$. 

For $\varphi = \varphi _{\max }$ we have $|AC| = 2 {\sqrt{d^2+1}}$, hence our
parallelogram lattice on the $xy$-plane is a lattice of $2 \times 2d$
rectangles. Then our lattice, $\Lambda _1 $, say, 
is generated by the vertices of a rectangular
pyramid $ABCDE$ with basis edges of lengths $2$ and $2d$, and lateral edges
$2$, which is one of the forms of giving the lattice in the Theorem.
The density of the corresponding unit ball packing is
$$
(4 \pi /3)/[2 \cdot 2d \cdot {\sqrt{2^2-(1+d^2)}} \,,
\tag 2.1.11
$$
as asserted in the Theorem.

For $\varphi = \varphi _{\min }$ the triangulation of the $xy$-plane consists
of isosceles triangles of base $2d$ and lateral sides $2$. Then as above, the
next horizontal layer has a centre $E$ of a unit ball
such that $|AE|=|BE|=|CE|=2$, where
$\Delta ABC$ is the triangle of the triangulation ${\Cal{T}}$
of the $xy$-plane, containing the 

\newpage

projection $E'$ of $E$ to the $xy$-plane.
The remaining three edges of the tetrahedron $ABCE$ have lengths $|AC|=|BC|=2$
and $|AB|=2d$. It remains
to show that this lattice, $\Lambda _2 $, say, 
is, up to congruence, the same lattice as $\Lambda _1$.

For convenience, we identify the vertices with the respective vectors. Then
$\Delta ACB$ and
$\Delta ACE$ and $\Delta BCE$ can be completed to parallelograms with vertices 
$A,C,B,A':=A+E-C$ and $A,C,E,B':=B+E-C$ and $B,C,E,C':=A+B-C$. 
The new vertices $A',B',C'$ belong to the respective point lattice. Then
$A'-A=B'-B=E-C$, hence the quadrangle $ABB'A'$ is a parallelogram, with edge
lengths $|AA'|=|BB'|=|CE|=2$ and $|AB|=|B'A'|=2d$. However, this
parallelogram is in fact a rectangle, since $\langle B-A, A'-A \rangle = 
\langle B-A,E-C \rangle =0$, since $|AC|=|BC|=|AE|=|BE|\,\,\,\,(=2)$,
and thus $[E,C]$ lies on the orthogonal bisector plane of $[A,B]$.
We assert that the distance
of $C'$ from each vertex of the parallelogram $ABB'A'$ is $2$. In fact, we
have $C'-A=B-C$ while $|BC|=2$, and similarly $C'-B$ has length $2$. Further,
$C'-A'=B-E$ while $|BE|=2$, and similarly $C'-B'$ has length $2$. Hence the
lattice $\Lambda _1$ is a sublattice of this new lattice $\Lambda _2$. 
Now we assert that
$$
\cases
{\text{the volume of the basic tetrahedron }} ABCE {\text{ of }} \Lambda _1 \\
{\text{and of the lattice tetrahedron }} \\
ACA'C'=AC(A+E-C)(A+B-C) {\text{ of }} \Lambda _2 \\
{\text{are equal, from which the equality of the lattices}} \\
\Lambda _1 {\text{ and }} \Lambda _2 {\text{( up to congruence) follows.}}
\endcases
\tag 2.1.12
$$

By a translation we achieve $E=0$ 
Then we have tetrahedra $0ABC$
and $AC(A-C)(A+B-C)$, or by a translation, $0(C-A)(-C)(B-C)$. The second
tetrahedron is the linear image of the first one by the matrix
$$
\left(
\matrix
-1 & 0 & 0 \\
0 & 0 & 1 \\
1 & -1 & -1
\endmatrix
\right)
$$
which has determinant $-1$. Hence we have showed \thetag{2.1.12}, which ends
the proof of the theorem.
$\blacksquare $
\enddemo

%%%%%%%%%%%%%%%%%%%%%%%%%%%%%%%%%%%%%%%%%%%%%%%%%%%%%%%%%%%%%%%%%%%%%%%%%%
%%%%%%%%%%%%%%%%%%%%%%%%%%%%%%%%%%%%%%%%%%%%%%%%%%%%%%%%%%%%%%%%%%%%%%%%%%%

\demo{Proof of Theorem 2.3}
For the two-dimensional square lattice $\Lambda '$ $( \subset
x_3x_4$-co\-or\-di\-na\-te plane) with edge length
$2$ the concentric balls of radius
${\sqrt{2}}$ cover the $x_3x_4$-coordinate plane. 
For the two-dimensional regular triangular lattice 
$\Lambda ''$ $( \subset x_3x_4$-coordinate plane) with edge length
$2$ the concentric balls of radius
$2/\sqrt{3}$ cover the $x_3x_4$-coordinate plane (and $2/\sqrt{3} <
\sqrt{2}$). Then we may apply Proposition 2.5,
since a lattice packing of unit balls satisfies its hypothesis \thetag{2.2}.
Thus any packing of translates of these sets $L$ has a density 
at most the 

\newpage

supremum of the densities of the lattice packings of $B^4$, with
corresponding point lattices $\Lambda + \Lambda '$, or $\Lambda + \Lambda ''$,
respectively, where $\Lambda $ is like in Proposition 2.5. Now recall that
for lattice packings of $B^4$ in ${\Bbb{R}}^4$
the maximal density is $\pi ^2/16$.
$\blacksquare $
\enddemo

%%%%%%%%%%%%%%%%%%%%%%%%%%%%%%%%%%%%%%%%%%%%%%%%%%%%%%%%%%%%%%%%%%%%%%%%%%
%%%%%%%%%%%%%%%%%%%%%%%%%%%%%%%%%%%%%%%%%%%%%%%%%%%%%%%%%%%%%%%%%%%%%%%%%%

\head{\S 6 Some remarks about the Main Lemmas in \S 4}\endhead

Observe that in Lemma 4.1 
we did not use the hypothesis about the upper bound
of the circumradii of the triangles. Probably also in the first inequality of
Lemma 4.2 and in Lemma 4.3 we may omit the
circumradius hypotheses and we have even then $p^2({\sqrt{7}}/4 + 1)$ and
$p^2({\sqrt{3}}/2 + 3/2)$ as lower bounds.
In Lemmas 4.4-4.6, i.e., when $\Delta ABC$ was right or obtuse
angled, we used the upper bound of the circumradius only for $\Delta
ABC$. However, the hypotheses in Lemmas 4.4-4.6 
imply that the circumradii of the
triangles $\Delta A_1BC$, $\Delta AB_2C$, $\Delta ABC_3$ 
are bounded above by the 
circumradius of $\Delta ABC$, so actually the upper bound of the circumradii
is valid for all triangles in Lemmas 4.4-4.6.

Moreover, in Lemmas 4.2-4.6 the average value of the areas of the considered
triangles
is strictly greater than $p^2/2$ (in Lemma 4.4 in the first case). 
This means that we could relax the hypothesis
that the circumradii are at most $p{\sqrt{2}}$ a bit, and still statements
corresponding to Lemmas 4.2-4.6
would give that the average values of the areas of the triangles in these
lemmas are at most
$p^2/2$, which would still make it possible to prove our Theorem 2.9 (as we
saw in \S 5, its proof in \S 5 only used Lemmas 4.1-4.6).

Of course, this would not imply a sharpening of our packing theorems in \S 2.
We show this on the example of 
Theorem 2.1. If in Theorem 2.1 we would allow $d=\sqrt{2}
+ \varepsilon > \sqrt{2}$, then we consider the second way of describing the
densest lattice packing of translates of $L$. It is generated by the vertices
of a rectangular right pyramid with base edges of lengths $2$ and $2d$, and
lateral edges of length $2$. The vertices of the base form a $2$ by $2d$
rectangular lattice in a plane, and closely packed translates of this 
rectangular lattice form the densest packing of translates of $L$. 
Observe that the height of our
pyramid is $\sqrt{3-d^2}$. This is at least $1$ for $1 \le d \le
{\sqrt{2}}$, but for $d>{\sqrt{2}}$ it is less than $1$, which means that
the closely packed translates of the rectangular lattice being
just above and just below the plane of the original rectangular lattice
overlap. Locally we may attain the density corresponding to two closely packed
neighbourly translates of our original rectangular lattice, but this cannot be
continued over the whole space.
Thus we cannot obtain a sharp upper bound for the
packing densities in Theorem 2.1 for $d > {\sqrt{2}}$. This implies that also
Propositions 2.4-2.6 and Theorem 2.8 cannot be extended to any $d >
{\sqrt{2}}$, and any covering radius greater than ${\sqrt{2}}$, and any $m(L) <
1/{\sqrt{2}}$, and any $m(L) < 1/{\sqrt{2}}$, respectively.

Now we return to Lemmas 4.4-4.6. Clearly
Lemma 4.4 is 

\newpage

false for any $R > p
{\sqrt{5/2}} = p \cdot 1.5811 \dots $. Probably
in Lemma 4.5 we could admit such circumradii $R$, 
that the convex quadrangle $ABA_1C$ satisfying that the circumradius of
$\Delta ABC$ is $R$ and $c_1=b_1=b=p$ and $\alpha _1 = \pi /2$ should
have an area at least $p^2$, which gives 
that $R \le p {\sqrt{5/2}} $, like for Lemma 4.4. 
Similarly, probably also in Lemma 4.6 
we could admit circumradii $R \le p  {\sqrt{5/2}} $. As numerical
evidence, we give for $R = p {\sqrt{5/2}} $ and $c_1=b_1=a_2=c_2=p$ and
either for (1) $\alpha _1 = \beta _2 = \pi /2$, or for (2) $\alpha _1 = \pi /2$
and $\beta + \beta _2 = \pi $, or for (3) $\alpha + \alpha _1 = \beta + \beta
_2 = \pi $ the numerical values of the total area of the three triangles. These
are (1) $p^2 \cdot 1.8 $, or (2) $p^2 \cdot 1.7 $ (exact values!), or (3) 
$p^2 \cdot 1.8624 \ldots ???$, all of which are greater than $p^2 \cdot 1.5$.
(Observe that for $R=p{\sqrt{2}}$ the corresponding values decreased in this
order, so the
eventual proofs for $R = p {\sqrt{5/2}}$ should be more complicated than those
in our paper.)
However, because of the rather technical
proofs in \S 3, where we essentially used for our numerical calculations
the inequality that the circumradius of $\Delta ABC$
is at most $p {\sqrt{2}}$, it is not clear
whether the proofs of Lemmas 4.3, 
4.5 and 4.6, with average area of the triangles at least
$p^2/2$, can be done for all circumradii at most $p{\sqrt{5/2}}$.

%%%%%%%%%%%%%%%%%%%%%%%%%%%%%%%%%%%%%%%%%%%%%%%%%%%%%%%%%%%%%%%%%%%%%%%%%%
%%%%%%%%%%%%%%%%%%%%%%%%%%%%%%%%%%%%%%%%%%%%%%%%%%%%%%%%%%%%%%%%%%%%%%%%%%%

\vskip.1cm

{\bf{Acknowledgement.}} We express our gratitude to J. Moln\'ar for calling
our attention to the fact that also his $L^*$-decomposition (cf. \cite{M77},
\cite{M78}) makes possible a proof of our theorems.

%%%%%%%%%%%%%%%%%%%%%%%%%%%%%%%%%%%%%%%%%%%%%%%%%%%%%%%%%%%%%%%%%%%%%%%%%%
%%%%%%%%%%%%%%%%%%%%%%%%%%%%%%%%%%%%%%%%%%%%%%%%%%%%%%%%%%%%%%%%%%%%%%%%%%%

\end%document 

%%%%%%%%%%%%%%%%%%%%%%%%%%%%%%%%%%%%%%%%%%%%%%%%%%%%%%%%%%%%%%%%%%%%%%%%%%%%%
%%%%%%%%%%%%%%%%%%%%%%%%%%%%%%%%%%%%%%%%%%%%%%%%%%%%%%%%%%%%%%%%%%%%%%%%%%%%%
%%%%%%%%%%%%%%%%%%%%%%%%%%%%%%%%%%%%%%%%%%%%%%%%%%%%%%%%%%%%%%%%%%%%%%%%%%%%

\demo{Proof of Theorem 2.8}
We may suppose that $$
S = \sup _i(\sup g_{ii})=1 \,.
\tag 2.1.1
$$ 
Then by hypothesis of the theorem we have for $I =
\inf_{i \ne j} (\inf g_{ij})$ 
that 
$$
I \ge 1/{\sqrt{2}} \,.
\tag 2.1.2
$$
We begin the proof with
saturating the packing, in the following way. We pick some infinitely many
translates of $L_i$'s in our packing, whose rotation axes intersected with the
$xy$-plane give points whose distances from the origin tend to $\infty $
sufficiently fast, so that their union has a
density (in the sense described in the Theorem) $0$. Let this subset of the
packing be ${\Cal{L}}_0$. This cannot be done only
if the packing is finite, but then its density is $0$. Let $k=1,2,\ldots $.
For each $k$, consecutively, we saturate our packing in a circle of radius $k$
about the origin, in the $x_1x_2$-plane, in the following way. We translate
some translates of $L_i$'s, belonging to ${\Cal{L}}_0$, to saturate our
packing in the circle of radius $R$ about the origin in the $x_1x_2$-coordinate
plane, but we only use elements of ${\Cal{L}}_0$ whose rotation axes intersect
the $x_1x_2$-coordinate plane at a distance at least $k+100$ from the
origin. This ensures that an already saturated packing in a circle of radius
$k$ will not be changed in the later steps, either by adding or by removing
some translate of some $L_i$. Thus our process ``converges'', and finally we
get a saturated packing of translates of the $L_i$'s. Since ${\Cal{L}}_0$ has
a density $0$, by this way of saturation the density does not decrease.
 
So we may restrict our attention to saturated packings in this sense.
Now consider the intersections of the axes of rotational symmetry of the
translates of the $L_i$'s with the $x_1x_2$-plane. We consider the Delone
triangulation of this plane belonging to this system of points (if there are
Delone-polygons with larger numbers of sides, they are triangulated in an
arbitrary way). By saturatedness any triangle in this triangulation has a
circumradius at most (actually, smaller than) $2$. 
The number density of this system of points is
$1/(2\overline{A})$, where $\overline{A}$ is the average area of the
Delone-triangles. (For this observe that the total angle sum of all Delone
triangles in some large circle about the origin is about $2 \pi $ times
the total number of vertices.)

Let $\Delta ABC$ be a non-obtuse triangle of our Delone
triangulation. We are going
to show that its area is at least 
$$
\inf \{ f \left( i(1),i(2),i(3) \right) \mid i(1),i(2),i(3) \in \{ 1,2,
\ldots \} {\text{ are different}} \} \,.
$$

Then $A,B,C$ are the intersections of the axes of rotation of suitably
translated copies of some $L_{i(1)}, L_{i(2)}$ and 
$L_{i(3)}$ with the $x_1x_2$-plane.
Let these translated copies be $L_{i(1)} + (x_{11},x_{12}, \ldots ,x_{1d})$, 
$L_{i(2)} + (x_{21},x_{22}, \ldots ,x_{2n})$ and 
$L_{i(3)} + (x_{31},x_{32}, \ldots ,x_{3n})$. 
Then by the definition of the function $g_{i(1)i(2)}$ we have
$$
\cases
|AB|= {\sqrt{(x_{21} - x_{11})^2 + (x_{22} - x_{12})^2}} \\
\ge 2g_{i(1)i(2)}(x_{23} - x_{13}, \ldots, x_{2n} - x_{1n})\,.
\endcases
\tag 2.1.3 
$$ 
The analogous inequalities hold also for $|BC|$ and $|CA|$ with the
functions $g_{i(2)i(3)}$ and $g_{i(3)i(1)}$. That is, the sides of the
triangle $\Delta ABC$ are at least certain values of the functions 
$g_{i(1)i(2)}$, $g_{i(2)i(3)}$ and $g_{i(3)i(1)}$. 
By \thetag{2.1.2} the sides of 
$\Delta ABC$ are at least ${\sqrt{2}}$. Further, evidently we have for $i \ne
j$ that $2 \sup g_{ij} \le \sup g_{ii} + \sup g_{jj}$. Thus the 
values of the functions 
$g_{i(1)i(2)}$, $g_{i(2)i(3)}$ and $g_{i(3)i(1)}$ lie in the interval 
$[ {1/\sqrt{2}}, 1]$. Then using \thetag{3.0} the area of $\Delta ABC$ is at
least the area of the triangle $T$ with sides 
$2g_{i(1)i(2)}(x_{23} - x_{13}, \ldots, x_{2n} - x_{1n})$, 
$2g_{i(2)i(3)}(x_{33} - x_{23}, \ldots, x_{3n} - x_{2n})$ and
$2g_{i(3)i(1)}(x_{13} - x_{33}, \ldots, x_{1n} - x_{3n})$,
i.e.,
$$
V(\Delta ABC) \ge V(T)\,.
\tag 2.1.4
$$ 
Observe that a triangle
with maximal side$/$minimal side $\le {\sqrt{2}}$ exists (by the triangle
inequality) and is non-obtuse (by the cosine law), 
thus $T$ is a non-obtuse triangle (as promised in the beginning of \S 2). 
By rotational symmetry of the $L_i$'s with respect to the
$x_3 \ldots x_n$-coordinate plane, we can change the first two coordinates of
the translation vectors $v_1:=(x_{11},x_{12}, \ldots ,x_{1d})$, 
$v_2:=(x_{21},x_{22}, \ldots ,x_{2n})$ and 
$v_3:=(x_{31},x_{32}, \ldots ,x_{3n})$, thus
obtaining new translation vectors 
$v_1':=(x_{11}',x_{12}', x_{13}, \ldots ,$
\newline
$x_{1d})$, 
$v_2':=(x_{21}',x_{22}', x_{23}, \ldots ,x_{2n}')$ and 
$v_3':=(x_{31}',x_{32}', x_{33}, \ldots ,x_{3n})$ such that the respective new
translated copies of $L_{i(1)},L_{i(2)}$ and $L_{i(3)}$ are pairwise touching
(in the sense precised in the beginning of \S 2). The
projections of $v_1', v_2', v_3'$ to the $x_1x_2$-coordinate plane are denoted
by $A',B',C'$.
Let $T':=\Delta v_1'v_2'v_3'$; its projection to the $x_1x_2$-coordinate
plane is $T=\Delta A'B'C'$. Then $A',B',C'$ 
are the intersections of the axes of rotations of $L_{i(1)}+v_1'$, $L_{i(2)}+
v_2'$ and $L_{i(3)}+v_3'$ with the $x_1x_2$-coordinate plane. By the
definition of the function $f\left( i(1),i(2),i(3) \right) $ (otherwise said,
varying the coordinates $x_{13}, \ldots ,x_{1n}$, 
$x_{23}, \ldots ,x_{2n}$ and $x_{33}, \ldots ,x_{3n}$ of the translation
vectors in an arbitrary way, and then taking infimum over these coordinates),
we obtain
$$
V(T) \ge f\left( i(1),i(2),i(3) \right) \,.
\tag 2.1.5
$$
Then \thetag{2.1.2} and \thetag{2.1.3} together give the first statement of the
theorem.

If all $L_i$'s are translates of the same set $L$, then
we consider the (inhomogeneous) two-dimensional
body lattice spanned by the translates of $L$ by the three new
translation vectors $v_1', v_2', v_3'$. Any two of these three translates of
$L$ are touching. Then, in the spanned body lattice, any two translates of $L$
which are simultaneous translates of any two of these three mutually
touching translates are also touching. 
If $L$ and $L+v'$ are two translates of $L$ from this body lattice which are
not simultaneous translates of any two of our three mutually touching
translates $L+v_1', L+v_2', L+v_3'$, then the following holds. If we decompose
our body lattice to strings (one-dimensional body lattices) 
spanned by suitable two of our three mutually touching translates of $L$, then 
$L$ and $L+v'$ are not in the same string, but also not in neighbourly strings.
The distance of the axes of rotation of two neighbourly strings is the height of
a triangle corresponding to three mutually touching copies of $L$,
with sides in $[2I, 2S]
\,\,\,\, = [2I, 2] \,\,(\subset 
[{\sqrt{2}}, 2])$. The minimal height belongs to the
maximal side (for which the adjacent angles are acute), 
and for given length of the maximal side, it 
is attained if the other two sides are minimal possible, i.e., have length 
$2I$. 
Hence the height is at least $\sqrt{(2I)^2 - S^2}
=  \sqrt{(2I)^2 - 1} \ge \sqrt{2-1} = 1$. Therefore the distance of the
axes of rotation of two neighbourly strings is at least $1$, and the distance
between the axes of rotation of two different strings which are not neighbourly,
is at least $2 \,\,\,\,(\ge \sup 2g_{11})$ 
(for $L_1:=L$). Therefore the translates 
$L$ and $L+v'$ are disjoint. This shows that our two-dimensional lattice
arrangement of translates of $L$ is in fact a two-dimensional lattice packing,
spanned by some three mutually touching translates of $L$.

Once more, 
varying the coordinates $x_{13}, \ldots ,x_{1n}$, 
$x_{23}, \ldots ,x_{2n}$ and $x_{33}, \ldots ,x_{3n}$ of the translation
vectors in an arbitrary way, and then taking infimum over these coordinates,
we get that \thetag{2.1.1} and \thetag{2.1.2} hold, with $L_{i(1)}=L_{i(2)}=
L_{i(3)}=L$, which shows the second statement of the theorem.
$\blacksquare $
\enddemo

%%%%%%%%%%%%%%%%%%%%%%%%%%%%%%%%%%%%%%%%%%%%%%%%%%%%%%%%%%%%%%%%%%%%%%%%%%%%

Now we distinguish cases according to the values of $\alpha , \beta , \gamma $.
We may have 
\newline
(a)
$\gamma = \pi /2$, or, e.g., 
\newline
(b)
$\alpha = \pi /2$, or
\newline
(c)
$\alpha , \beta , \gamma < \pi /2$.

In case (a) we apply Lemma 3.1 to $\Delta BA_1C$ and $\Delta CB_2A$ (rather
than to $\Delta ABC$ in Lemma 3.1), with lower bound $p$ for $c_1, b_1, a_2,
c_2$, and upper bound $\pi - \alpha $ on 
$\alpha _1$ and $ \pi - \beta $ on $\beta_2$. Thus we obtain 
$$
\min \{ b_1, c_1 \} = \min \{ a_2, c_2 \} = p\,.
$$
We may suppose that $c_1=c_2=p$. Then the area of the pentagon $BA_1CB_2A$
is the sum of the areas of the symmetric trapezoid $BA_1B_2A$ and of $\Delta
A_1CB_2$. We may move $C$ on the circumcircle till e.g., $a_2=p$, while
$b_1 \ge p$. Then the
pentagon is already uniquely determined, if the circumradius $R$ is given. The
total area of our pentagon is $3 \cdot p {\sqrt{R^2-p^2/4}}/2 + R^2 \sin
\left( \pi /2 - 3 \arcsin (p/(2R)) \right) $.
This area is unchanged if we permute the sides, with preservation of their
central angles. So we may rather consider the case when $c_2=a_2=b_1=p \le c_1$.
Then 
Let $\varphi \le \pi /4$ 
be 
the central angle $\varphi \,\,( \le \pi )$ of a chord of length $p$ in the 
circumcircle of our pentagon satisfies $\varphi \le \pi /4$. Let the 
circumradius and the centre of our pentagon be $R$ and $O$. 
Then the area of our pentagon is 
$(R^2/2) \cdot \left( 3 \sin \varphi + \sin ( \pi - 3 \varphi )\right) = (R^2/2)
\cdot \left( 3 \sin \varphi + \sin (3 \varphi ) \right) $, where $\varphi /2 =
p/(2R)$. Here $(R^2/2) \cdot \sin \varphi = R^2 \cdot \sin
(\varphi /2) \cos (\varphi /2) = R^2 \cdot \left( p / (2R) \right) \sqrt{ 1 - 
\left( p / (2R) \right) ^2} = (p/2) \sqrt{R^2 - (p/2)^2}$ monotonically
increases with $R$. Further, $R^2 \cdot \sin (3 \varphi ) =  $

%$V(\Delta BAB_2)+V(\Delta BB_2C)+V(\Delta BCA_1)$.

%%%%%%%%%%%%%%%%%%%%%%%%%%%%%%%%%%%%%%%%%%%%%%%%%%%%%%%%%%%%%%%%%%%%%%%%%%%%%

\demo{Proof of Theorem 2.8}
We may suppose that $$
S = \sup _i(\sup g_{ii})=1 \,.
\tag 2.1.1
$$ 
Then by hypothesis of the theorem we have for $I =
\inf_{i \ne j} (\inf g_{ij})$ 
that 
$$
I \ge 1/{\sqrt{2}} \,.
\tag 2.1.2
$$
We begin the proof with
saturating the packing, in the following way. We pick some infinitely many
translates of $L_i$'s in our packing, whose rotation axes intersected with the
$xy$-plane give points whose distances from the origin tend to $\infty $
sufficiently fast, so that their union has a
density (in the sense described in the Theorem) $0$. Let this subset of the
packing be ${\Cal{L}}_0$. This cannot be done only
if the packing is finite, but then its density is $0$. Let $k=1,2,\ldots $.
For each $k$, consecutively, we saturate our packing in a circle of radius $k$
about the origin, in the $x_1x_2$-plane, in the following way. We translate
some translates of $L_i$'s, belonging to ${\Cal{L}}_0$, to saturate our
packing in the circle of radius $R$ about the origin in the $x_1x_2$-coordinate
plane, but we only use elements of ${\Cal{L}}_0$ whose rotation axes intersect
the $x_1x_2$-coordinate plane at a distance at least $k+100$ from the
origin. This ensures that an already saturated packing in a circle of radius
$k$ will not be changed in the later steps, either by adding or by removing
some translate of some $L_i$. Thus our process ``converges'', and finally we
get a saturated packing of translates of the $L_i$'s. Since ${\Cal{L}}_0$ has
a density $0$, by this way of saturation the density does not decrease.
 
So we may restrict our attention to saturated packings in this sense.
Now consider the intersections of the axes of rotational symmetry of the
translates of the $L_i$'s with the $x_1x_2$-plane. We consider the Delone
triangulation of this plane belonging to this system of points (if there are
Delone-polygons with larger numbers of sides, they are triangulated in an
arbitrary way). By saturatedness any triangle in this triangulation has a
circumradius at most (actually, smaller than) $2$. 
The number density of this system of points is
$1/(2\overline{A})$, where $\overline{A}$ is the average area of the
Delone-triangles. (For this observe that the total angle sum of all Delone
triangles in some large circle about the origin is about $2 \pi $ times
the total number of vertices.)

Let $\Delta ABC$ be a non-obtuse triangle of our Delone
triangulation. We are going
to show that its area is at least 
$$
\inf \{ f \left( i(1),i(2),i(3) \right) \mid i(1),i(2),i(3) \in \{ 1,2,
\ldots \} {\text{ are different}} \} \,.
$$

Then $A,B,C$ are the intersections of the axes of rotation of suitably
translated copies of some $L_{i(1)}, L_{i(2)}$ and 
$L_{i(3)}$ with the $x_1x_2$-plane.
Let these translated copies be $L_{i(1)} + (x_{11},x_{12}, \ldots ,x_{1d})$, 
$L_{i(2)} + (x_{21},x_{22}, \ldots ,x_{2n})$ and 
$L_{i(3)} + (x_{31},x_{32}, \ldots ,x_{3n})$. 
Then by the definition of the function $g_{i(1)i(2)}$ we have
$$
\cases
|AB|= {\sqrt{(x_{21} - x_{11})^2 + (x_{22} - x_{12})^2}} \\
\ge 2g_{i(1)i(2)}(x_{23} - x_{13}, \ldots, x_{2n} - x_{1n})\,.
\endcases
\tag 2.1.3 
$$ 
The analogous inequalities hold also for $|BC|$ and $|CA|$ with the
functions $g_{i(2)i(3)}$ and $g_{i(3)i(1)}$. That is, the sides of the
triangle $\Delta ABC$ are at least certain values of the functions 
$g_{i(1)i(2)}$, $g_{i(2)i(3)}$ and $g_{i(3)i(1)}$. 
By \thetag{2.1.2} the sides of 
$\Delta ABC$ are at least ${\sqrt{2}}$. Further, evidently we have for $i \ne
j$ that $2 \sup g_{ij} \le \sup g_{ii} + \sup g_{jj}$. Thus the 
values of the functions 
$g_{i(1)i(2)}$, $g_{i(2)i(3)}$ and $g_{i(3)i(1)}$ lie in the interval 
$[ {1/\sqrt{2}}, 1]$. Then using \thetag{3.0} the area of $\Delta ABC$ is at
least the area of the triangle $T$ with sides 
$2g_{i(1)i(2)}(x_{23} - x_{13}, \ldots, x_{2n} - x_{1n})$, 
$2g_{i(2)i(3)}(x_{33} - x_{23}, \ldots, x_{3n} - x_{2n})$ and
$2g_{i(3)i(1)}(x_{13} - x_{33}, \ldots, x_{1n} - x_{3n})$,
i.e.,
$$
V(\Delta ABC) \ge V(T)\,.
\tag 2.1.4
$$ 
Observe that a triangle
with maximal side$/$minimal side $\le {\sqrt{2}}$ exists (by the triangle
inequality) and is non-obtuse (by the cosine law), 
thus $T$ is a non-obtuse triangle (as promised in the beginning of \S 2). 
By rotational symmetry of the $L_i$'s w.r.t. the
$x_3 \ldots x_n$-coordinate plane, we can change the first two coordinates of
the translation vectors $v_1:=(x_{11},x_{12}, \ldots ,x_{1d})$, 
$v_2:=(x_{21},x_{22}, \ldots ,x_{2n})$ and 
$v_3:=(x_{31},x_{32}, \ldots ,x_{3n})$, thus
obtaining new translation vectors 
$v_1':=(x_{11}',x_{12}', x_{13}, \ldots ,$
\newline
$x_{1d})$, 
$v_2':=(x_{21}',x_{22}', x_{23}, \ldots ,x_{2n}')$ and 
$v_3':=(x_{31}',x_{32}', x_{33}, \ldots ,x_{3n})$ such that the respective new
translated copies of $L_{i(1)},L_{i(2)}$ and $L_{i(3)}$ are pairwise touching
(in the sense precised in the beginning of \S 2). The
projections of $v_1', v_2', v_3'$ to the $x_1x_2$-coordinate plane are denoted
by $A',B',C'$.
Let $T':=\Delta v_1'v_2'v_3'$; its projection to the $x_1x_2$-coordinate
plane is $T=\Delta A'B'C'$. Then $A',B',C'$ 
are the intersections of the axes of rotations of $L_{i(1)}+v_1'$, $L_{i(2)}+
v_2'$ and $L_{i(3)}+v_3'$ with the $x_1x_2$-coordinate plane. By the
definition of the function $f\left( i(1),i(2),i(3) \right) $ (otherwise said,
varying the coordinates $x_{13}, \ldots ,x_{1n}$, 
$x_{23}, \ldots ,x_{2n}$ and $x_{33}, \ldots ,x_{3n}$ of the translation
vectors in an arbitrary way, and then taking infimum over these coordinates),
we obtain
$$
V(T) \ge f\left( i(1),i(2),i(3) \right) \,.
\tag 2.1.5
$$
Then \thetag{2.1.2} and \thetag{2.1.3} together give the first statement of the
theorem.

If all $L_i$'s are translates of the same set $L$, then
we consider the (inhomogeneous) two-dimensional
body lattice spanned by the translates of $L$ by the three new
translation vectors $v_1', v_2', v_3'$. Any two of these three translates of
$L$ are touching. Then, in the spanned body lattice, any two translates of $L$
which are simultaneous translates of any two of these three mutually
touching translates are also touching. 
If $L$ and $L+v'$ are two translates of $L$ from this body lattice which are
not simultaneous translates of any two of our three mutually touching
translates $L+v_1', L+v_2', L+v_3'$, then the following holds. If we decompose
our body lattice to strings (one-dimensional body lattices) 
spanned by suitable two of our three mutually touching translates of $L$, then 
$L$ and $L+v'$ are not in the same string, but also not in neighbourly strings.
The distance of the axes of rotation of two neighbourly strings is the height of
a triangle corresponding to three mutually touching copies of $L$,
with sides in $[2I, 2S]
\,\,\,\, = [2I, 2] \,\,(\subset 
[{\sqrt{2}}, 2])$. The minimal height belongs to the
maximal side (for which the adjacent angles are acute), 
and for given length of the maximal side, it 
is attained if the other two sides are minimal possible, i.e., have length 
$2I$. 
Hence the height is at least $\sqrt{(2I)^2 - S^2}
=  \sqrt{(2I)^2 - 1} \ge \sqrt{2-1} = 1$. Therefore the distance of the
axes of rotation of two neighbourly strings is at least $1$, and the distance
between the axes of rotation of two different strings which are not neighbourly,
is at least $2 \,\,\,\,(\ge \sup 2g_{11})$ 
(for $L_1:=L$). Therefore the translates 
$L$ and $L+v'$ are disjoint. This shows that our two-dimensional lattice
arrangement of translates of $L$ is in fact a two-dimensional lattice packing,
spanned by some three mutually touching translates of $L$.

Once more, 
varying the coordinates $x_{13}, \ldots ,x_{1n}$, 
$x_{23}, \ldots ,x_{2n}$ and $x_{33}, \ldots ,x_{3n}$ of the translation
vectors in an arbitrary way, and then taking infimum over these coordinates,
we get that \thetag{2.1.1} and \thetag{2.1.2} hold, with $L_{i(1)}=L_{i(2)}=
L_{i(3)}=L$, which shows the second statement of the theorem.
$\blacksquare $
\enddemo

%%%%%%%%%%%%%%%%%%%%%%%%%%%%%%%%%%%%%%%%%%%%%%%%%%%%%%%%%%%%%%%%%%%%%%%%%%%%%%%%

$p^2 \cdot (5/4) \left( \cot ( \pi /10) + \right.$
\newline
$\left. {\sqrt{7}}/8 \right) = p^2 \cdot
(1.7204 \ldots + 0.3307 \ldots ) = 2.0511 \ldots $.

%%%%%%%%%%%%%%%%%%%%%%%%%%%%%%%%%%%%%%%%%%%%%%%%%%%%%%%%%%%%%%%%%%%%%%%%%%%%%%

Then fixing $C,A,B$ and
moving $A_1$, $V(\Delta BA_1C)$ 
becomes minimal if $\min \{
b_1,c_1 \} = p$ ---  e.g., for $b_1=p$. 
Then fixing $A_1, C, A$ and moving
$B$, $V(\Delta ABA_1)$ becomes minimal, like above,
e.g., for $c=p$. Since
now $V(\Delta A_1CA)$ remains constant, the total area of the
triangles  $\Delta ABA_1$ and $\Delta A_1CA$, 
that is the area of our quadrangle,
becomes minimal, e.g., for $c=p$. 

%%%%%%%%%%%%%%%%%%%%%%%%%%%%%%%%%%%%%%%%%%%%%%%%%%%%%%%%%%%%%%%%%%%%%%%%%%%%%

By $c=b_1=p$ and $b > c_1 \ge p$ we have 
$$
\cases
\pi = \angle AOB + \angle BOA_1+\angle A_1OC = 
2 \arcsin ( c/b ) + 
2 \arcsin ( c_1/b) \\ 
+ 2 \arcsin ( b_1/b ) 
= 4 \arcsin ( p/b ) + 2 \arcsin ( c_1/b ) \ge 
6 \arcsin ( p/b ) \,,
\endcases
$$
therefore $b \ge 2p$. 

%%%%%%%%%%%%%%%%%%%%%%%%%%%%%%%%%%%%%%%%%%%%%%%%%%%%%%%%%%%%%%%%%%%%%%%%%%%%%

{\bf{3.}}
We turn to case (3).
First we suppose that $\alpha , \beta , \gamma < \pi /2$. Then we apply Lemma
3.3 to the triangles $\Delta ABC$ and $\Delta BA_1C$ (with the same notations
here as in Lemma 3.3). The case of the triangles $\Delta ABC$ and $\Delta CB_2A$
can be settled in an identical way.

By $\alpha _1 \ge \pi /2$ we have $a=a_1 \ge p{\sqrt{2}} >p$, so the side
hypothesis is satisfied if we apply Lemma 3.3 (the other side lengths 
are unchanged). The angular hypothesis is to be investigated for $\alpha
+ \alpha _1$, but that decreases if we apply Lemma 3.3, and for $\beta + \beta
_2$, but that is samller than $\pi $ by the hypothesis of (3). 
So by Lemma 3.3 either $\alpha _1 = \pi /2$ or $\beta + \beta _2 
= \pi /2$. 
Analogously one obtains either $\beta _2 = \pi /2$ or $\alpha + \alpha _1
= \pi $. Thus we have four cases, but two of them have a symmetrical role, and
can be settled as one case.

We begin with the case $\alpha _1 = \beta _2 = \pi /2$.
Then $b_1, c_1 \ge p$ gives $a=a_1 \ge p{\sqrt{2}}$, 
and analogously, $b=b_1 \ge p{\sqrt{2}}$. Moreover, we have $c \ge p$.
Then \thetag{3.0} yields $V(\Delta ABC) \ge p^2 \cdot {\sqrt{7}}/4$. Further,
$V(\Delta BA_1C) = b_1c_1 \sin \alpha _1 /2 \ge p^2/2$ and similarly $V(\Delta
CB_2A) \ge p^2/2$, which give the inequality in the Lemma. If moreover, we
have $c \ge p{\sqrt{2}}$, then \thetag{3.0} yields even 
$V(\Delta ABC) \ge p^2 \cdot {\sqrt{3}}/2$, and then the total area of the
three triangles is at least $p^2 \cdot (1+{\sqrt{3}}/2) = p^2 \cdot
1.8660 \ldots > p^2 \cdot (5/4) \cot ( \pi /10) = p^2 \cdot 1.7204 \ldots $.

Therefore we may suppose that $\Delta ABC$ is a right triangle. It suffices to
investigate the cases $\beta = \pi /2$ (the case $\alpha = \pi /2$ 
is analogous) and $\gamma = \pi /2$.

We begin with the case $\beta = \pi /2$. By the angular hypothesis then also
$\beta _2 = \pi /2$. However then $\beta + \beta _2 = \pi $, contrary to case
(3).

We apply Lemma 3.1 to $\Delta BA_1C$ (in place of $\Delta
ABC$ from Lemma 3.1), with $\pi - \alpha $ in place of $\gamma _0$ from Lemma
3.1. We gain that either $\alpha + \alpha _1 < \pi $, or $\alpha + \alpha _1
= \pi /2$ and one of $b_1$ and $c_1$ equals $p$. 

In case $\alpha + \alpha _1 < \pi $ we apply Lemma 3.4 to $\Delta ABC$ and
$\Delta CB_2A$ (in place of $\Delta BCA$ and $\Delta AC_1B$ in Lemma 3.4). The
side hypothesis is preserved, since $b \ge p \cdot {\sqrt{2}} >p$. The angular
hypothesis remains preserved, since $\beta + \beta _2$ decreases. 
 
We continue with the case $\gamma = \pi /2$. We decrease $\gamma $ a bit,
rotating $\Delta BA_1C$ and $\Delta CB_2A$ towards each other.
Then $V(\Delta ABC)$ decreases, $V(\Delta BA_1C)$ and $V\Delta CB_2A)$ remain
constant, so their sum decreases. The side hypothesis is to be checked only
for $c$, but that was originally by $a,b \ge p$ and $\gamma = \pi /2$ at
least $p {\sqrt{2}}$, so the side hypothesis remains valid. Similarly,
originally we had $\alpha , \beta < \pi /2$, so they remain acute angles. So
all the hypotheses of this Lemma (and also those of
case (3)) remain preserved, and the total area of our
three triangles decreased. This is a contradiction.
$\blacksquare $
\enddemo

%%%%%%%%%%%%%%%%%%%%%%%%%%%%%%%%%%%%%%%%%%%%%%%%%%%%%%%%%%%%%%%%%%%%%%%%%%%%%%

\demo{Proof} 
Consider an extremal configuration (this exists).

We will copy the proof of Lemma 4.5, with the respective modifications.

We apply Lemma 3.2 to $\Delta ABC$ (with the same notations as there), with
$a_0=a$, $b_0:=b$ and $r_0:=p {\sqrt{2}}$. Observe that by application of
Lemma 3.2 the side hypotheses 
remain preserved, and also the angular hypotheses are preserved since by
strictly increasing $\gamma \,\,( \ge \pi /2)$, both angles $\alpha $ and
$\beta $ of $\Delta ABC$ strictly decrease. The supposed inequalities $\alpha
_1 , \beta _2, \gamma \ge  \pi /2$ also remain preserved. Therefore Lemma 3.2
implies that 
$$
{\text{the circumradius of }} \Delta ABC {\text{ is }}  p {\sqrt{2}}\,.
\tag 4.6.1
$$

We are going to show that 
$$
b_1=c_1=a_2=c_2=p\,.
$$
Clearly it suffices to prove 
$$
b_1=c_1=p
$$
(the proof of the remaining two equalities is analogous).
We restrict our attention to the quadrangle $ABA_1C$.

First suppose
\newline
(1) $\alpha + \alpha _1 = \pi $, 
\newline
i.e. that the quadrangle $ABA_1C$ has a
circumcircle. Like in the proof of lemma 4.5, we may suppose that $b_1=p$.
The area of the quadrangle $ABA_1C$ also equals $V(\Delta ABA_1)+V(\Delta
A_1CA)$. Like in the proof of Lemma 4.5, moving $B$ on the circumcircle of the
quadrangle $ABA_1C$ towards $A_1$, we have that $V(\Delta ABA_1)$ strictly
decreases, while $V(\Delta A_1CA)$ and $V(CB_2A)$ remain constant. Therefore
also $V(\Delta ABC) + V(\Delta BA_1C) + V(\Delta CB_2A)$ strictly decreases.
Therefore in the extremal configuration we have also $c_1=p$.
 
Second suppose
\newline 
(2) $\alpha + \alpha _1 < \pi $.
\newline
Like in the proof of Lemma 4.5, by an application of Lemma 3.1 we get
$b_1=c_1=p$.

Therefore in any of cases (1) and (2) we have
$$
b_1=c_1=p, {\text{ and analogously, also }} a_2=c_2=p\,.
\tag 4.6.2
$$

Now we apply Lemma 3.5 to the quadrangles $ABA_1C$ and $ABCB_2$, in order 
to show that 
$$
{\text{either }} \alpha _1 = \pi /2 {\text{ or }} \alpha + \alpha _1 = \pi \,,
\tag 4.6.3
$$
and
$$
{\text{either }} \beta _2 = \pi /2 {\text{ or }} \beta + \beta _2 = \pi \,.
\tag 4.6.4
$$
It suffices
to consider the quadrangle $ABA_1C$, and to show \thetag{4.6.3}. 
We use here the same notations as in
Lemma 3.5. All the hypotheses of Lemma 3.5 are satisfied by \thetag{4.6.1}
and \thetag{4.6.2}, except for $b=|AC| \le p \sqrt{7/2}$. Now we show this
last inequality. 

By the angular hypothesis and $\beta _2 \ge \pi /2$ we have that the
circumradius of $\Delta CB_2A$ is at most the circumradius of $\Delta ABC$,
i.e., $p {\sqrt{2}}$ (cf. \thetag{4.6.1}).
By \thetag{4.6.2} we have $a_2=c_2=p$. We apply Lemma 3.2 to $\Delta
CB_2A$, rather than to $\Delta ABC$ in Lemma 3.2, 
with $p$ in place of $a_0$ and $b_0$ in Lemma 3.2 and with $p{\sqrt{2}}$ in
place of $r_0$ in Lemma 3.2. Thus we will get that the circumradius of $\Delta
CB_2A$ is $p{\sqrt{2}}$. (Observe that by the proof of Lemma 3.2 
$\beta _2 \,\,(\ge \pi /2)$ and the circumradius of $\Delta CB_2A$ increase
simultaneously.)
To establish this, however we have to show that by
changing $\Delta CB_2A$ suitably we can
preserve the side hypothesis and the angular
hypothesis, and also the inequalities $\gamma , \alpha _1, \beta _2 \ge \pi
/2$.

Preservation of the side hypothesis is clear: the only changed sides are
$a_2,c_2,b_2$. However, $a_2$ and $c_2$ are fixed, and equal $p$. Then, since
$\beta _2 \ge  \pi /2$, therefore $b_2 \ge {\sqrt{a_2^2 + c_2^2}} \ge
p{\sqrt{2}} > a_2=c_2=p$.

Preservation of the angular hypothesis is necessary to verify only for $\beta
+ \beta _2$. Originally we had $\beta + \beta _2 \le \pi $. Then we increased
the common side $b_2=b$ of $\Delta ABC$ and $\Delta CB_2A$. This means
increasing of $\beta + \beta _2$. If during the changing of $\Delta CB_2A$ we
arrive to $\beta + \beta _2 = \pi $, then \thetag{4.6.4} is proved, and no
more changing of $\Delta CB_2A$ is necessary. 
 
Lastly, the inequalities $\gamma , \alpha _1, \beta _2 \ge \pi /2$ are
preserved, since $\gamma $ and $\alpha _1$ are not changed, and $\beta _2$ is
increased.

This proves \thetag{4.6.4}, and analogously, \thetag{4.6.3}. That is, either
$\alpha _1 = \beta _2 = \pi /2$, or $\alpha + \alpha _1 = \beta + \beta _2
= \pi $, or, e.g., $\alpha _1 = \pi /2$ and $\beta + \beta _2 = \pi $. In all
of these cases the figure is uniquely determined, and simple calculations show
that the total area of our three triangles equals either 
$p^2 \cdot ({\sqrt{3}}/2 +
1)= p^2 \cdot 1.8660 \ldots $, or $p^2 \cdot (5{\sqrt{7}} + 7 {\sqrt{3}}
+8)/16 = p^2 \cdot 2.0845 \ldots $, or $p^2 \cdot 
(29 {\sqrt{7}})/32 = p^2 \cdot 
2.3977 \ldots $. This finishes the proof of our lemma.
$\blacksquare $
\enddemo

%%%%%%%%%%%%%%%%%%%%%%%%%%%%%%%%%%%%%%%%%%%%%%%%%%%%%%%%%%%%%%%%%%%%%%%%%%%%%%%

First we show that the circumradius $R$
of $\Delta ABC$ equals $p{\sqrt{2}}$. In fact, if $R < p{\sqrt{2}}$, then we
increase $\gamma $, thus rotating apart the triangles $\Delta BA_1C$ and
$\Delta CB_2A$ about their common vertex $C$, the side hypotheses are
preserved, since the only changed side is $c$ which increases, and the angualr
hypotheses are not affected by this rotating apart.

%%%%%%%%%%%%%%%%%%%%%%%%%%%%%%%%%%%%%%%%%%%%%%%%%%%%%%%%%%%%%%%%%%%%%%%%%%%%%%%

\proclaim{Lemma 3.9}
Let a packing consisting of the triangles $\Delta ABC$, $\Delta BA_1C$, 
$\Delta CB_2A$ satisfy the side hypothesis and
let $c_1=b_1=a_2=c_2=p$
and let $\alpha + \alpha _1, \beta + \beta _2 < \pi$. 
Let the circumradius of $\Delta ABC$ be
$p{\sqrt{2}}$. Let the central angles corresponding to chords of the
circumcircle of $\Delta ABC$ of lengths
$a,b,p$ be $\varphi _A$, $\varphi _B$ and $\varphi _0 = 
\arcsin ({\sqrt{7}}/4) = 41.4096\ldots ^{\circ }$. 
Let $\pi /3 \le \varphi _A, \varphi _B \le 2 \varphi
_0$. 
Then $V(\Delta ABC) + V(\Delta BA_1C) +V(\Delta CB_2A)
\ge p^2(1+{\sqrt{3}}/2) = 1.8660 \ldots $. Equality holds if and only if 
$\varphi _A = \varphi _B = \pi /3 $.
\endproclaim

%%%%%%%%%%%%%%%%%%%%%%%%%%%%%%%%%%%%%%%%%%%%%%%%%%%%%%%%%%%%%%%%%%%%%%%%%%
%%%%%%%%%%%%%%%%%%%%%%%%%%%%%%%%%%%%%%%%%%%%%%%%%%%%%%%%%%%%%%%%%%%%%%%%%%%

\demo{Proof} 
Let $R$ be the circumradius of $\Delta ABC$. Then we have
$$
V(\Delta ABC) 
= (R^2/2)\left( \sin \varphi _A + \sin \varphi _B - \sin (\varphi _A
+ \varphi _B) \right) \,,
\tag 3.9.1
$$
$$
V(\Delta BA_1C) 
= R \sin ( \varphi _A /2) {\sqrt{p^2 - R^2 \sin ^2 (\varphi _A /2)}}
\tag 3.9.2
$$
and
$$
V(\Delta CB_2A) 
= R \sin ( \varphi _B /2) {\sqrt{p^2 - R^2 \sin ^2 (\varphi _B /2)}}\,.
\tag 3.9.3
$$
Adding these inequalities and taking in consideration $R = p{\sqrt{2}}$ and 
$\left( 1 - 2 \sin ^2 (\varphi _A / \right.$
\newline
$\left. 2) \right) ^{1/2} = (\cos \varphi _A)^{1/2}$ 
(and its
analogue for $\varphi _B$, where under the square root sign there are positive
numbers by the hypotheses of this lemma)
we get
$$
\cases
\left( V(\Delta ABC) + V(\Delta BA_1C) + V(\Delta CB_2A) \right) /p^2 =
\sin \varphi _A  + \sin \varphi _B \\
- \sin (\varphi _A + \varphi _B) + 
{\sqrt{2}} \sin (\varphi _A /2) {\sqrt{\cos \varphi _A}} +
{\sqrt{2}} \sin (\varphi _B /2) {\sqrt{\cos \varphi _B}}\,.
\endcases
\tag 3.9.4
$$
Then \thetag{3.9.4} is a function of two variables $\varphi _A$ and $\varphi
_B$.
The derivative of \thetag{3.9.4} with respect to $\varphi _A$ is
$$
\cases
\cos \varphi _A  - \cos (\varphi _A + \varphi _B) + 
{\sqrt{2}}(1/2) \cos (\varphi _A / 2) {\sqrt{\cos \varphi _A}} \\
+ {\sqrt{2}} \sin (\varphi _A
/2) \left( 1/(2{\sqrt{ \cos \varphi _A}}) \right) (-\sin \varphi _A) = \\
\left( \cos \varphi _A - \cos ( \varphi _A + \varphi _B) \right) + \\
(1/{\sqrt{2}}) \left( \cos (\varphi _A /2) \cos \varphi _A - \sin (\varphi
A/2) \sin \varphi _A \right) / {\sqrt{\cos \varphi _A}} \\
= \left( \cos \varphi _A - \cos ( \varphi _A + \varphi _B) \right) +
(1/{\sqrt{2}}) \cos (3 \varphi _A /2) / {\sqrt{\cos \varphi _A}}
\,.
\endcases
\tag 3.9.5
$$

We assert strict monotone 
decreasing of \thetag{3.9.5} as a function of $\varphi _A$,
that is, strict convexity of \thetag{3.9.4} as a function of $\varphi _A$. 

The last expression in \thetag{3.9.5} is a sum of two summands. 
Differentiating the first summand we obtain
$$
\sin ( \varphi _A + \varphi _B) - \sin \varphi _A \,. 
\tag 3.9.6
$$
This is nonpositive (in fact negative in $( \pi /3, 2 \varphi _0)$)
by 
$$
\varphi _A +
\varphi _B \in [2 \pi /3, \pi ) {\text{ and }} \varphi _A \in [\pi
/3, \pi /2)\,,
\tag 3.9.7
$$ 
implying strict monotone decreasing of the first summand in \thetag{3.9.5}.
The second summand is a fraction, with numerator nonpositive, and denumerator
positive. Its absolute value 
$$
\left( -(1/{\sqrt{2}}) \cos (3 \varphi _A /2) \right) / {\sqrt{\cos \varphi _A}}
\tag 3.9.8
$$
has nonnegative numerator and positive denumerator, the first one
strictly 
increasing, the second one strictly
decreasing, hence \thetag{3.9.8} is strictly increasing, 
implying that the second summand of \thetag{3.9.5} is strictly
decreasing as well.

This shows that also \thetag{3.9.5} is a strictly
decreasing function of $\varphi _A$,
i.e., \thetag{3.9.4} is a strictly convex function of $\varphi
_A$. Therefore \thetag{3.9.4} as a function of $\varphi _A$
attains its minimum only at an endpoint of the
interval $[\pi /3, 2 \varphi _0]$. 
Similarly, \thetag{3.9.4} as a function of $\varphi _B$
attains its minimum only at an endpoint of the
interval $[\pi /3, 2 \varphi _0]$. 
Therefore the only possible minimum points of \thetag{3.9.4} for $(\varphi _A,
\varphi _B) \in [\pi /3, 2 \varphi _0] \times [\pi /3, 2 \varphi _0]$
are $(\varphi _A, \varphi _B) = 
(\pi /3 , \pi /3), \,\,( \pi
/3, 2 \varphi _0),\,\, (2 \varphi _0, \pi /3),\,\, (2 \varphi _0, 2 \varphi _0
)$. The values of \thetag{3.9.4} in these points are $p^2(1+ {\sqrt{3}}/2
p^2 \cdot 1.8660 \ldots $,
$p^2 \cdot 2.0845 \ldots $, $p^2 \cdot 2.0845 \ldots $, 
$p^2 \cdot 2.3977\ldots $, respectively, hence the first one is the minimum.

Unicity of the minimum follows from the proof.
$\blacksquare $
\enddemo

%%%%%%%%%%%%%%%%%%%%%%%%%%%%%%%%%%%%%%%%%%%%%%%%%%%%%%%%%%%%%%%%%%%%%%%%%%

\demo{Proof} 
Consider an extremal configuration (this exists).

We have two angular hypotheses, namely $\alpha + \alpha _1 \le \pi $ and
$\beta + \beta _2 \le \pi $.
Analogously as in the proof of Lemma 4.2, we make the following case
distinctions. Either
\newline
(1) 
$\alpha + \alpha _1 = \beta + \beta _2 = \pi $, or, e.g.,
\newline
(2)
$\alpha + \alpha _1 , \beta + \beta _2 < \pi $.
\newline
(3)
$\alpha + \alpha _1 = \pi > \beta + \beta _2$, or

{\bf{1.}}
We begin with case (1). That is, the pentagon $ABA_1CB_2$ has a circumcircle.
Applying Lemma 3.8 we get $V(ABA_1CB_2) \ge p^2 \cdot 2.3977 \ldots > 
p^2 \cdot (1 + \sqrt{3}/2) = p^2 \cdot 1.8660 \ldots $.

{\bf{2.}}
We turn to case (2).
We apply Lemma 3.1 to $\Delta BA_1C$ and $\Delta CB_2A$ (in place of $\Delta
ABC$ in Lemma 3.1, with $ \pi - \alpha $ and $\pi - \beta $ in place of 
$\gamma $ in Lemma 3.1), 
obtaining
$$
c_1=b_1=a_2=c_2=p
\tag 4.6.1
$$
We apply Lemma 3.2 to $\Delta ABC$ (with the same notations here and there) 
obtaining for the circumradius $R$ of $\Delta ABC$
$$
R=p{\sqrt{2}} 
\tag 4.6.2
$$
(observe that now we need not care the angular hypotheses, and the side
hypothesis remains preserved by $c \ge {\sqrt{a^2 + b^2}} \ge p{\sqrt{2}} >
p$).

We denote the central angles belonging to the sides $a$ and $b$ in the
circumcircle of $\Delta ABC$ by $\varphi _A$ and $\varphi _B$. We denote the
central angle corresponding to a chord of length $p$ of this circumcircle
by $\varphi _0$; then we have by \thetag{4.6.2}
$$
\cases
\sin ( \varphi _0 /2) = p/(2R) = 1/{\sqrt{8}}\,, {\text{ i.e., }} R  =
p/\left( 2 \sin (\varphi _0 /2) \right) \,, \\
{\text{or }} \varphi _0 = \arcsin ({\sqrt{7}}/4) = 41.4096 \ldots ^{\circ }
\endcases
\tag 4.6.3
$$ 

Then we have
$$
V(ABC) = (R^2/2)\left( \sin \varphi _A + \sin \varphi _B - \sin (\varphi _A
+ \varphi _B) \right) \,,
\tag 4.6.4
$$
$$
V(BA_1C) = R \sin ( \varphi _A /2) {\sqrt{p^2 - R^2 \sin ^2 (\varphi _A /2)}}
\tag 4.6.5
$$
and
$$
V(CB_2A) = R \sin ( \varphi _B /2) {\sqrt{p^2 - R^2 \sin ^2 (\varphi _B /2)}}\,.
\tag 4.6.6
$$
Here by 
$$
a \ge {\sqrt{c_1^2 + b_1 ^2}} \ge p{\sqrt{2}} {\text{ and similarly }} b \ge 
p{\sqrt{2}}
\tag 4.6.7
$$
and by \thetag{4.6.2} we have 
$$
\varphi _A , \varphi _B \ge \pi /3 \,.
\tag 4.6.8
$$
On the other hand, the maximum of $\alpha _1$ or $\beta _2$ occurs when $A_1$
or $B_2$ lies on the circumcircle of $\Delta ABC$, respectively (else by
\thetag{4.6.1} $A_1$ or $B_2$ would lie inside
of the circumcircle of $\Delta ABC$, 
contradicting the angular hypothesis). That is,
$$
\varphi _A , \varphi _B \le 2 \varphi _0 = 82.8192 \ldots
\tag 4.6.10
$$
Now we apply Lemma 3.9 ...

{\bf{3.}}
We turn to case (3).
Analogously as in the proof of \thetag{3.8.7}, we have also here 
$$
b_1=c_1=p\,.
\tag 4.6.1 
$$
On the other hand, by Lemma 3.1 applied to $\Delta CB_2A$ (in place of $\Delta
ABC$ in Lemma 3.1, with $\pi - \beta $ in place of $\gamma $ in Lemma 3.1), 
we have 
$$
a_2=c_2=p\,.
\tag 4.6.2
$$

----------------------------------------------------------------------------

We will copy the proof of Lemma 4.5, with the respective modifications.

We apply Lemma 3.2 to $\Delta ABC$ (with the same notations as there), with
$a_0=a$, $b_0:=b$ and $r_0:=p {\sqrt{2}}$. Observe that by application of
Lemma 3.2 the side hypotheses 
remain preserved, and also the angular hypotheses are preserved since by
strictly increasing $\gamma \,\,( \ge \pi /2)$, both angles $\alpha $ and
$\beta $ of $\Delta ABC$ strictly decrease. The supposed inequalities $\alpha
_1 , \beta _2, \gamma \ge  \pi /2$ also remain preserved. Therefore Lemma 3.2
implies that 
$$
{\text{the circumradius of }} \Delta ABC {\text{ is }}  p {\sqrt{2}}\,.
\tag 4.6.1
$$

We are going to show that 
$$
b_1=c_1=a_2=c_2=p\,.
$$
Clearly it suffices to prove 
$$
b_1=c_1=p
$$
(the proof of the remaining two equalities is analogous).
We restrict our attention to the quadrangle $ABA_1C$.

First suppose
\newline
(1) $\alpha + \alpha _1 = \pi $, 
\newline
i.e. that the quadrangle $ABA_1C$ has a
circumcircle. Like in the proof of lemma 4.5, we may suppose that $b_1=p$.
The area of the quadrangle $ABA_1C$ also equals $V(\Delta ABA_1)+V(\Delta
A_1CA)$. Like in the proof of Lemma 4.5, moving $B$ on the circumcircle of the
quadrangle $ABA_1C$ towards $A_1$, we have that $V(\Delta ABA_1)$ strictly
decreases, while $V(\Delta A_1CA)$ and $V(CB_2A)$ remain constant. Therefore
also $V(\Delta ABC) + V(\Delta BA_1C) + V(\Delta CB_2A)$ strictly decreases.
Therefore in the extremal configuration we have also $c_1=p$.
 
Second suppose
\newline 
(2) $\alpha + \alpha _1 < \pi $.
\newline
Like in the proof of Lemma 4.5, by an application of Lemma 3.1 we get
$b_1=c_1=p$.

Therefore in any of cases (1) and (2) we have
$$
b_1=c_1=p, {\text{ and analogously, also }} a_2=c_2=p\,.
\tag 4.6.2
$$

Now we apply Lemma 3.5 to the quadrangles $ABA_1C$ and $ABCB_2$, in order 
to show that 
$$
{\text{either }} \alpha _1 = \pi /2 {\text{ or }} \alpha + \alpha _1 = \pi \,,
\tag 4.6.3
$$
and
$$
{\text{either }} \beta _2 = \pi /2 {\text{ or }} \beta + \beta _2 = \pi \,.
\tag 4.6.4
$$
It suffices
to consider the quadrangle $ABA_1C$, and to show \thetag{4.6.3}. 
We use here the same notations as in
Lemma 3.5. All the hypotheses of Lemma 3.5 are satisfied by \thetag{4.6.1}
and \thetag{4.6.2}, except for $b=|AC| \le p \sqrt{7/2}$. Now we show this
last inequality. 

By the angular hypothesis and $\beta _2 \ge \pi /2$ we have that the
circumradius of $\Delta CB_2A$ is at most the circumradius of $\Delta ABC$,
i.e., $p {\sqrt{2}}$ (cf. \thetag{4.6.1}).
By \thetag{4.6.2} we have $a_2=c_2=p$. We apply Lemma 3.2 to $\Delta
CB_2A$, rather than to $\Delta ABC$ in Lemma 3.2, 
with $p$ in place of $a_0$ and $b_0$ in Lemma 3.2 and with $p{\sqrt{2}}$ in
place of $r_0$ in Lemma 3.2. Thus we will get that the circumradius of $\Delta
CB_2A$ is $p{\sqrt{2}}$. (Observe that by the proof of Lemma 3.2 
$\beta _2 \,\,(\ge \pi /2)$ and the circumradius of $\Delta CB_2A$ increase
simultaneously.)
To establish this, however we have to show that by
changing $\Delta CB_2A$ suitably we can
preserve the side hypothesis and the angular
hypothesis, and also the inequalities $\gamma , \alpha _1, \beta _2 \ge \pi
/2$.

Preservation of the side hypothesis is clear: the only changed sides are
$a_2,c_2,b_2$. However, $a_2$ and $c_2$ are fixed, and equal $p$. Then, since
$\beta _2 \ge  \pi /2$, therefore $b_2 \ge {\sqrt{a_2^2 + c_2^2}} \ge
p{\sqrt{2}} > a_2=c_2=p$.

Preservation of the angular hypothesis is necessary to verify only for $\beta
+ \beta _2$. Originally we had $\beta + \beta _2 \le \pi $. Then we increased
the common side $b_2=b$ of $\Delta ABC$ and $\Delta CB_2A$. This means
increasing of $\beta + \beta _2$. If during the changing of $\Delta CB_2A$ we
arrive to $\beta + \beta _2 = \pi $, then \thetag{4.6.4} is proved, and no
more changing of $\Delta CB_2A$ is necessary. 
 
Lastly, the inequalities $\gamma , \alpha _1, \beta _2 \ge \pi /2$ are
preserved, since $\gamma $ and $\alpha _1$ are not changed, and $\beta _2$ is
increased.

This proves \thetag{4.6.4}, and analogously, \thetag{4.6.3}. That is, either
$\alpha _1 = \beta _2 = \pi /2$, or $\alpha + \alpha _1 = \beta + \beta _2
= \pi $, or, e.g., $\alpha _1 = \pi /2$ and $\beta + \beta _2 = \pi $. In all
of these cases the figure is uniquely determined, and simple calculations show
that the total area of our three triangles equals either 
$p^2 \cdot ({\sqrt{3}}/2 +
1)= p^2 \cdot 1.8660 \ldots $, or $p^2 \cdot (5{\sqrt{7}} + 7 {\sqrt{3}}
+8)/16 = p^2 \cdot 2.0845 \ldots $, or $p^2 \cdot 
(29 {\sqrt{7}})/32 = p^2 \cdot 
2.3977 \ldots $. This finishes the proof of our lemma.
$\blacksquare $
\enddemo

%%%%%%%%%%%%%%%%%%%%%%%%%%%%%%%%%%%%%%%%%%%%%%%%%%%%%%%%%%%%%%%%%%%%%%%%%%%%%

\proclaim{Lemma 3.9}
Let a packing consisting of the triangles $\Delta ABC$, $\Delta BA_1C$, 
$\Delta CB_2A$ satisfy the side hypothesis and
let $c_1=b_1=a_2=c_2=p$
and let $\alpha + \alpha _1, \beta + \beta _2 < \pi$. 
Let the circumradius of $\Delta ABC$ be
$p{\sqrt{2}}$. Let the central angles corresponding to chords of the
circumcircle of $\Delta ABC$ of lengths
$a,b,p$ be $\varphi _A$, $\varphi _B$ and $\varphi _0 = 
\arcsin ({\sqrt{7}}/4) = 41.4096\ldots ^{\circ }$. 
Let $\pi /3 \le \varphi _A, \varphi _B \le 2 \varphi
_0$. 
Then $V(\Delta ABC) + V(\Delta BA_1C) +V(\Delta CB_2A)
\ge p^2(1+{\sqrt{3}}/2) = 1.8660 \ldots $. Equality holds if and only if 
$\varphi _A = \varphi _B = \pi /3 $.
\endproclaim

%%%%%%%%%%%%%%%%%%%%%%%%%%%%%%%%%%%%%%%%%%%%%%%%%%%%%%%%%%%%%%%%%%%%%%%%%%
%%%%%%%%%%%%%%%%%%%%%%%%%%%%%%%%%%%%%%%%%%%%%%%%%%%%%%%%%%%%%%%%%%%%%%%%%%%

\demo{Proof} 
Let $R$ be the circumradius of $\Delta ABC$. Then we have
$$
V(\Delta ABC) 
= (R^2/2)\left( \sin \varphi _A + \sin \varphi _B - \sin (\varphi _A
+ \varphi _B) \right) \,,
\tag 3.9.1
$$
$$
V(\Delta BA_1C) 
= R \sin ( \varphi _A /2) {\sqrt{p^2 - R^2 \sin ^2 (\varphi _A /2)}}
\tag 3.9.2
$$
and
$$
V(\Delta CB_2A) 
= R \sin ( \varphi _B /2) {\sqrt{p^2 - R^2 \sin ^2 (\varphi _B /2)}}\,.
\tag 3.9.3
$$
Adding these inequalities and taking in consideration $R = p{\sqrt{2}}$ and 
$\left( 1 - 2 \sin ^2 (\varphi _A / \right.$
\newline
$\left. 2) \right) ^{1/2} = (\cos \varphi _A)^{1/2}$ 
(and its
analogue for $\varphi _B$, where under the square root sign there are positive
numbers by the hypotheses of this lemma)
we get
$$
\cases
\left( V(\Delta ABC) + V(\Delta BA_1C) + V(\Delta CB_2A) \right) /p^2 =
\sin \varphi _A  + \sin \varphi _B \\
- \sin (\varphi _A + \varphi _B) + 
{\sqrt{2}} \sin (\varphi _A /2) {\sqrt{\cos \varphi _A}} +
{\sqrt{2}} \sin (\varphi _B /2) {\sqrt{\cos \varphi _B}}\,.
\endcases
\tag 3.9.4
$$
Then \thetag{3.9.4} is a function of two variables $\varphi _A$ and $\varphi
_B$.
The derivative of \thetag{3.9.4} with respect to $\varphi _A$ is
$$
\cases
\cos \varphi _A  - \cos (\varphi _A + \varphi _B) + 
{\sqrt{2}}(1/2) \cos (\varphi _A / 2) {\sqrt{\cos \varphi _A}} \\
+ {\sqrt{2}} \sin (\varphi _A
/2) \left( 1/(2{\sqrt{ \cos \varphi _A}}) \right) (-\sin \varphi _A) = \\
\left( \cos \varphi _A - \cos ( \varphi _A + \varphi _B) \right) + \\
(1/{\sqrt{2}}) \left( \cos (\varphi _A /2) \cos \varphi _A - \sin (\varphi
A/2) \sin \varphi _A \right) / {\sqrt{\cos \varphi _A}} \\
= \left( \cos \varphi _A - \cos ( \varphi _A + \varphi _B) \right) +
(1/{\sqrt{2}}) \cos (3 \varphi _A /2) / {\sqrt{\cos \varphi _A}}
\,.
\endcases
\tag 3.9.5
$$

We assert strict monotone 
decreasing of \thetag{3.9.5} as a function of $\varphi _A$,
that is, strict convexity of \thetag{3.9.4} as a function of $\varphi _A$. 

The last expression in \thetag{3.9.5} is a sum of two summands. 
Differentiating the first summand we obtain
$$
\sin ( \varphi _A + \varphi _B) - \sin \varphi _A \,. 
\tag 3.9.6
$$
This is nonpositive (in fact negative in $( \pi /3, 2 \varphi _0)$)
by 
$$
\varphi _A +
\varphi _B \in [2 \pi /3, \pi ) {\text{ and }} \varphi _A \in [\pi
/3, \pi /2)\,,
\tag 3.9.7
$$ 
implying strict monotone decreasing of the first summand in \thetag{3.9.5}.
The second summand is a fraction, with numerator nonpositive, and denumerator
positive. Its absolute value 
$$
\left( -(1/{\sqrt{2}}) \cos (3 \varphi _A /2) \right) / {\sqrt{\cos \varphi _A}}
\tag 3.9.8
$$
has nonnegative numerator and positive denumerator, the first one
strictly 
increasing, the second one strictly
decreasing, hence \thetag{3.9.8} is strictly increasing, 
implying that the second summand of \thetag{3.9.5} is strictly
decreasing as well.

This shows that also \thetag{3.9.5} is a strictly
decreasing function of $\varphi _A$,
i.e., \thetag{3.9.4} is a strictly convex function of $\varphi
_A$. Therefore \thetag{3.9.4} as a function of $\varphi _A$
attains its minimum only at an endpoint of the
interval $[\pi /3, 2 \varphi _0]$. 
Similarly, \thetag{3.9.4} as a function of $\varphi _B$
attains its minimum only at an endpoint of the
interval $[\pi /3, 2 \varphi _0]$. 
Therefore the only possible minimum points of \thetag{3.9.4} for $(\varphi _A,
\varphi _B) \in [\pi /3, 2 \varphi _0] \times [\pi /3, 2 \varphi _0]$
are $(\varphi _A, \varphi _B) = 
(\pi /3 , \pi /3), \,\,( \pi
/3, 2 \varphi _0),\,\, (2 \varphi _0, \pi /3),\,\, (2 \varphi _0, 2 \varphi _0
)$. The values of \thetag{3.9.4} in these points are $p^2(1+ {\sqrt{3}}/2
p^2 \cdot 1.8660 \ldots $,
$p^2 \cdot 2.0845 \ldots $, $p^2 \cdot 2.0845 \ldots $, 
$p^2 \cdot 2.3977\ldots $, respectively, hence the first one is the minimum.

Unicity of the minimum follows from the proof.
$\blacksquare $
\enddemo

%%%%%%%%%%%%%%%%%%%%%%%%%%%%%%%%%%%%%%%%%%%%%%%%%%%%%%%%%%%%%%%%%%%%%%%%%%%%

%, and then the angles $\angle OBA_1
%= \angle OCA_1$ decrease or increase. 

%%%%%%%%%%%%%%%%%%%%%%%%%%%%%%%%%%%%%%%%%%%%%%%%%%%%%%%%%%%%%%%%%%%%%%%%%%%%

if one unit ball lies in 2-nd, 4-th, ... neighbouring layer of the layer 
in which the other unit ball lies, the centres of the balls in one layer
forming a regular triangular lattice with edge length $2$ 
in a plane in ${\Bbb{R}}^3$.

%%%%%%%%%%%%%%%%%%%%%%%%%%%%%%%%%%%%%%%%%%%%%%%%%%%%%%%%%%%%%%%%%%%%%%%%%%%

preserving $T(2,0,0)$ and $T(0,2,0)$,
we can rotate $T(1,0,d)$
about the $z$-axis in the positive sense a bit, 
obtaining a new linear transformation $T'$, such that $\| T'(0,1,{\sqrt{2}}) -
T'(0,2,0) \| $ once more becomes $2$. By Pythagoras' theorem, for $\varepsilon
$ small, also this rotation will be small. Then $\| T'(2,0,0) - T'(0,2,0) \|
=2$ and, for $\varepsilon $ sufficiently small, all eight lattice vectors
$T'( \pm 2, \pm 2, \pm {\sqrt{2}})$ will have norms either $2$, or a bit
larger than $2$, and all other non-zero lattice vectors in $T'\Lambda $ remain
longer than some constant strictly larger than $2$. The axes of rotation of
our vertical string of unit balls meet the $xy$-plane in a point lattice whose
basic parallelogram is a rectangle of diagonal $2$, which is close to a
square, for $\varepsilon $ sufficiently small. At the same time, the Delone
decomposition of the $xy$-plane corresponding to these points of intersection
consists of translates of this rectangle, joining along entire sides. We
conjecture that for $\varepsilon >0$ sufficiently small, the above constructed
packing is the densest packing of translates of our string of unit balls, with
centres equidistant with distances $2d \in (2{\sqrt{2}},
2({\sqrt{2}} + \varepsilon )]$.

%%%%%%%%%%%%%%%%%%%%%%%%%%%%%%%%%%%%%%%%%%%%%%%%%%%%%%%%%%%%%%%%%%%%%%%%%%

%%%%%%%%%%%%%%%%%%%%%%%%%%%%%%%%%%%%%%%%%%%%%%%%%%%%%%%%%%%%%%%%%%%%%%%%%%%%

Observe that a right triangle has at least
such an area as any other triangle with two sides equal to the legs of the
right triangle, thus 
$$
2m(L)^2 \ge M(L)\sqrt{\left( 2m(L) \right) ^2 - M(L)^2}
\tag 2.9.1
$$ 
(the last
expression here being the area of a triangle of sides $2m(L)$, $2m(L)$, 
$2M(L)$). 

Now we show that, for any $\varepsilon >0$, for a suitable packing
a non-obtuse Delone triangle can have an area at most 
$M(L)\sqrt{\left( 2m(L) \right) ^2 - M(L)^2} + \varepsilon $. 
Let us consider the following
translates of $L$: $L+(M(L),0, \ldots ,0)$, $L-(M(L),0, \ldots ,0)$ (these are
touching translates in the sense of \thetag{2.5}) 
and $L+(0,x_2,x_3, \ldots , x_n)$, where $x_2$ is chosen
so that $L$ and any of $L \pm (0,x_2,x_3, \ldots , x_n)$ should be touching,
for any given $(x_3, \ldots , x_n) \in {\Bbb{R}}^{n-2}$. (This is possible by
the rotational symmetry of $L$ about the $x_3 \ldots  x_n$-coordinate plane.)
Then take the supremum of the densities of all two-dimensional lattice packings 
of translates of $L$ (in the sense of (2.5))
spanned by the above three mutually touching translates 
of $L$, when $(x_3, \ldots , x_n) \in {\Bbb{R}}^{n-2}$ is arbitrary.

The corresponding two-dimensional point lattices project orthogonally to the
$x_1x_2$-plane injectively onto two-dimensional point lattices in the
$x_1x_2$-coordinate plane, since $m(L)>0$.

\thetag{2.9.1} says that $\min \{ 2m(L)^2 , V_0 \}
= V_0$. Hence, for any packing of translates of $L$,
the number density of the intersections
of the axes of rotation of the translates of $L$ is at most 
$1/(2V_0)$,
which implies by hypothesis (2.9)
the theorem.

%%%%%%%%%%%%%%%%%%%%%%%%%%%%%%%%%%%%%%%%%%%%%%%%%%%%%%%%%%%%%%%%%%%%%%%%%%
%%%%%%%%%%%%%%%%%%%%%%%%%%%%%%%%%%%%%%%%%%%%%%%%%%%%%%%%%%%%%%%%%%%%%%%%%%%

$$
\cases
V(\Delta ABC) + V(\Delta BA_1C) + V(\Delta CB_2A) = \\
V(\Delta ABA_1) + V(\Delta A_1CA) + V(\Delta CB_2A) \,.
\endcases
\tag 4.2.7
$$
Now, replacing the vertices of the quadrangle $ABA_1C$ by the respective
vertices  of the quadrangle $BA_1CB_2A$, the same considerations as in case
(a) give, rather than \thetag{4.2.6},  
$$
\cases
V(\Delta ABC) + V(\Delta BA_1C) +V(\Delta AB_2C) = \\
V(\Delta ABA_1) + V(\Delta A_1CA) +V(\Delta AB_2C) \ge \\
p^2 ( {\sqrt{2}}/2 + {\sqrt{7}}/8) = p^2 \cdot 1.7449 \ldots > \\
p^2 \cdot (5/4) \cot {\pi /10} = p^2 \cdot 1.7204 \ldots > \\
p^2({\sqrt{7}}/4 + 1) = p^2 \cdot 1.6614 \ldots \,,
\endcases
\tag 4.2.8
$$
which proves both lower estimates of the lemma in case (b) as well.

%%%%%%%%%%%%%%%%%%%%%%%%%%%%%%%%%%%%%%%%%%%%%%%%%%%%%%%%%%%%%%%%%%%%%%%%

Observe that we have a convex
deltoid $OBA_1C$ (and analogously $OCB_2A$) where $O$
is the circumcentre of $\Delta ABC$. This deltoid has side lengths
$|OB|=|OC|=R$ and $|BA_1|=|CA_1|=p$. Then in this deltoid $\varphi _A, \alpha
_1, \alpha + \alpha _1 = \varphi _A /2 + \alpha _1, 
a=|BC|$ increase or decrease simultaneously.
%, and then the angles $\angle OBA_1
%= \angle OCA_1$ decrease or increase. 
Then the minimum  
of $\alpha _1$ or $\beta _2$ occurs when $a$ or $b$ is minimal, i.e.,
$a=p{\sqrt{2}}$ or $b=p{\sqrt{2}}$
(cf. \thetag{4.6.5}). The maximum of
$\alpha _1$ or $\beta _2$ occurs when $A_1$
or $B_2$ lies on the circumcircle of $\Delta ABC$, respectively, i.e., $\alpha
+ \alpha _1 $ or $\beta + \beta _1 $
attains its maximum $ \pi $. 

%%%%%%%%%%%%%%%%%%%%%%%%%%%%%%%%%%%%%%%%%%%%%%%%%%%%%%%%%%%%%%%%%%%%%%%%%%%%

Consider the DV-cells for
the system of centres of the balls in $\Lambda $. 
We may suppose that $E'$ lies in the DV-cell of the origin.
The minimal distance of $E$
and the centre of some ball of $\Lambda $ (which is at least $2$)
is attained exactly when the minimal
distance of $E'$ and the centre of some ball of $\Lambda $ is attained, and
this happens when we consider the distance of $E'$ and the origin. Then this
last distance is maximized when $E'$ is a vertex of the DV-cell of the origin.
Observe that we have a right triangle $0EE'$, whose hypotenuse $0E$ has length
at least $2$. Then 
$$
|EE'|={\sqrt{|0E|^2-|0E'|^2}} \ge {\sqrt{4-|0E'|^2}} \,,
\tag 2.1.4
$$
thus we can bound $|EE'|$ from below by setting $|0E|=2$ and letting $E'$ be a
vertex of the DV-cell of the origin.

{\bf{1.}} If $\gamma
= \pi /2$, then by \thetag{1.4} we have 
$2 =|BC| = |AC|$ and $d={\sqrt{2}}$. Then the lattice $\Lambda
$ is a square lattice of edge length $2$, and the Delone polygons are squares
of side length $2$, and the DV-cells also constitute a square lattice of edge
length $2$ (a translate of the Delone sundivision). 
Then the distances of the vertices of the DV-cell of $0$
to $0$ are $|0E'|={\sqrt{2}}$, hence by \thetag{2.1} 
$$
|EE'| \ge {\sqrt{2}} \,.
\tag 2.1.5
$$
Then the neighbourly horizontal layer has a distance at least ${\{sqrt{2}}$
from the $xy$-plane. And this can be realized: the densest lattice
packing of $B^3$ can be given as the union of closely packed square lattice
horizontal layers with edge length $2$. Thus we have a lattice packing, namely
the densest lattice packing of translates of $B^3$.

{\bf{2.}}
If $\beta = \pi /2$, then $ABCD$ is a rectangle and the Delone polygons are 
just the rectangles consisting the tiling ${\Cal{T}}$, 
and the DV-cells also constitute a $2 \times 2d$ rectangular lattice
(a translate of the Delone subdivision).
Then the distances of the vertices of the DV-cell of $0$
to $0$ are $(0E'|={\sqrt{d^2+1}}$, hence by \thetag{2.1} 
$$
|EE'| \ge {\sqrt{3-d^2}} \,.
\tag 2.1.6
$$
For case of equality $E'$ should be the centre of a rectangle consisting the
tiling ${\Cal{T}}$, and then $E$ should be the vertex of a rectangular pyramid 
with base a $2 \times 2d$ rectangle and with lateral edges $2$.
That
is, our lattice is spanned by the vertices of a rectangular pyramid, of base a
$2 \times 2d$ rectangle, 
and of lateral edges of length $2$. This is the first way of
giving this lattice in the theorem. The volume of the basic
parallelepiped is $2d \cdot 2 \cdot {\sqrt{2^2-(d^2+1)}}$, and the
density of our lattice arrangement is as stated in the theorem. 
We discuss {\bf{2}} further in {\bf{4}} below.

{\bf{3.}}
Since by \thetag{1.2} we have $\alpha <  \pi /2$, after {\bf{1}} and {\bf{2}}
there remains the case when $\Delta ABC$ is acute. 
Then the tiling ${\Cal{T}}$ consists of
acute triangles, such that at each vertex six triangles
meet. The DV-cell of a vertex has then as vertices the circumcentres of the six
triangles of ${\Cal{T}}$ meeting at that vertex.
Now $E$ has a distance $2$ from three vertices of some triangle $\Delta ABC$ 
of the Delone triangulation 
${\Cal{T}}$, with $E'$ being the circumcentre of $\Delta ABC$, that
now lies in the interior of the acute triangle $\Delta ABC$.  
We discuss {\bf{4}} further in {\bf{4}} below.

{\bf{4.}}
Now we investigate cases {\bf{3}} and {\bf{4}} together. Still we have
to show that we have a lattice {\it{packing}}. In a horizontal layer the
corresponding open unit balls are disjoint. Two open unit balls in neighbourly
horizontal
layers are disjoint by the considerations given above, when the minimal
distance of $E$ from the centres of our unit balls lying in the $xy$-plane was
determined. 
It remains to show that the open unit 
balls in at least second neighbour horizontal layers are disjoint. For this it
is sufficient to show that the height $h_E$
of the tetrahedron $ABCE$ corresponding
to the vertex $E$ --- i.e., when $ABCD$ is a rectangle of ${\Cal{T}}$ (case
{\bf{2}}) or
$E'$ lies in the acute triangle $\Delta ABC$ of ${\Cal{T}}$ (case {\bf{3}}), 
--- is at least $1$. Recall that in both cases {\bf{2}} and {\bf{3}}
$E'$ is the circumcentre of 
$\Delta ABC$. Let $R$ be the circumradius of $\Delta ABC$. 
Then $h_E={\sqrt{2^2-R^2}}$. Therefore we have to prove that $R \le
{\sqrt{3}}$. 
Observe that $E'$ lies on the perpendicular bisectors of the sides $AB$
and $BC$. If $|AC|$ increases from its minimum $2$ till its maximum
$2{\sqrt{d^2+1}}$ then the
angle $\beta = \angle ABC$ strictly increases, and by elementary geometric
considerations $E'$ moves on the perpendicular
bisector of side $AB$ farther from side $AB$. 
Then $E'$ will be the farthest from side 
$AB$ when $|AC|$ attains its maximum $2{\sqrt{d^2+1}}$, which happens if we
have a rectangular lattice on the $xy$-plane, when $R={\sqrt{d^2+1}} \le
{\sqrt{3}}$.

If $ABCD$ is no rectangle, then we have a triangulation of the $xy$-plane to
acute triangles, such that at each vertex of the triangulation six triangles
meet. The DV-cell of a vertex has then vertices the circumcentres of the six
triangles meeting at that vertex, and the dual tiling, the Delone
triangulation, is just the triangulation of the $xy$-plane constructed above.
Now $E$ has a distance $2$ from three vertices of some triangle $\Delta ABC$ 
of the 
Delone triangulation, with $E'$ being the circumcentre of $\Delta ABC$, that
now lies in the interior of the acute triangle $\Delta ABC$.  

%%%%%%%%%%%%%%%%%%%%%%%%%%%%%%%%%%%%%%%%%%%%%%%%%%%%%%%%%%%%%%%%%%%%%%%%%%%

Actually, by $2 =|BC| \le |AC|$ we have $\alpha < \beta $ and thus 
$$
\alpha  < \pi /2 \,.
\tag 2.1.2
$$

%%%%%%%%%%%%%%%%%%%%%%%%%%%%%%%%%%%%%%%%%%%%%%%%%%%%%%%%%%%%%%%%%%%%%%%%%

%and that some
%ball in our two-dimensional ball lattice has as centre the origin.

%%%%%%%%%%%%%%%%%%%%%%%%%%%%%%%%%%%%%%%%%%%%%%%%%%%%%%%%%%%%%%%%%%%%%%%%

%(this is different from the above choice of
%the origin, but this will cause no difficulty). 

%%%%%%%%%%%%%%%%%%%%%%%%%%%%%%%%%%%%%%%%%%%%%%%%%%%%%%%%%%%%%%%%%%%%%%%

Now we already have the formulas for the coordinates of our points in order to
prove that the system of translates of our string described in the paragraph
below the Theorem is in fact a packing.

Since a string itself forms a packing,
it suffices to show that the $0$-th string and any other string form a
packing. 
In fact, then an analogous statement holds for the $ \pm 2$'nd, $\pm
4$-th, etc. strings rather than for the $0$-th string. 
Finally, the $0$'th and the $1$'st strings can be taken
over to each other by a central symmetry with respect to the point
$(A_0+B_0)/2$, preserving our system of translates of our string,
hence also the $1$'st string 
and any other string forms a packing. Then an analogous statement holds for 
the $-1$'st, $\pm 3$-th, etc. strings, rather than for the $1$'st string.

In the $0$'th string some circle centres are
$A_0$ and $A_1$. In the $1$'st string some circle centre is
$B_0$, so the $0$'th and $1$'st strings together form a packing by definition. 

In the $2$'nd
string some circle centre is $C_1$. Its projection to the axis of the $0$'th
string lies in $[A_0,A_1]$, even actually in $[(A_0+3A_1)/4, A_1]$. Hence
the smallest distance of $C_1$ to the centre of a circle in the $0$'th string 
is $|C_1A_1|=2$. So the $0$'th and $2$'nd strings together form a packing.

In the $3$'rd string some circle centre is $C_1+(B_0-A_1)$.  
Its projection to the axis of the $0$'th
string lies in $[A_0,A_1]$, even actually in $[A_0,(A_0+A_1)/2]$. Hence
the smallest distance of $C_1+(B_0-A_1)$ 
to the centre of a circle in the $0$'th string is $\| C_1+(B_0-A_1) - A_0 \| =
\| \left( d+ {\sqrt{3}} {\sqrt{4-d^2}}, {\sqrt{3}}d+3{\sqrt{4-d^2}} \right) /2
\| $. We want to prove that this is at least $2$. We use Pythagoras theorem,
we can square the inequality to be proved, and rearranging we obtain
$$
4+2 {\sqrt{3}}d {\sqrt{4-d^2}} \ge d^2 \,.
\tag 0.2
$$
However, by $d \in ({\sqrt{3}},2)$ we have even $4 > d^2$, so the
inequality to be proven holds.
Hence the $0$'th and $3$'rd strings together form a packing.

The distance of the axes of the $0$'th string and the at least $4$-th string
is at least the distance of the axes of the $0$'th and $4$-th strings, i.e.,
the double of the distance of the $0$'th and $2$'nd strings, i.e., the
double of the $y$-coordinate of $C_1$, which double is 
${\sqrt{3}}d+{\sqrt{4-d^2}}$.
For the packing property it suffices to show that this is at least $2$. 
We can square the inequality to be proved, and then by rearranging we get 
$$
2d^2+2{\sqrt{3}}d{\sqrt{4-d^2}} \ge 0 \,,
\tag 0.3
$$
which holds, so the inequality to be proven holds.
Hence the $0$'th and at least $4$'th strings together form a packing.
This ends the proof 
that the system of translates of our string described in the paragraph
below the Theorem is in fact a packing.

%%%%%%%%%%%%%%%%%%%%%%%%%%%%%%%%%%%%%%%%%%%%%%%%%%%%%%%%%%%%%%%%%%%%%%%%%%